\documentclass[10pt]{article}
\usepackage[T1]{fontenc}    
\usepackage[utf8]{inputenc}
\usepackage[english]{babel}
\usepackage{lmodern}
\usepackage{algorithm}
\usepackage{algpseudocode}
\usepackage{authblk}

\usepackage{amsfonts}
\usepackage{amsmath}
\usepackage{amssymb}
\usepackage{amsthm}
\usepackage{bm} 			
\usepackage{esint}			
\usepackage{enumitem}

\usepackage[pdftex]{graphicx}
\usepackage{graphics}
\usepackage{float}
\usepackage{tikz}

\usepackage{lscape} 			
\usepackage{scrextend} 			

\usepackage{parskip}
\usepackage{color} 			
\usepackage{cancel}			
\usepackage{stmaryrd}			
\usepackage[sort,nocompress]{cite} 	
\usepackage{hyperref}                   
\usepackage{multirow}

\newcommand{\R}{\mathbb{R}}

\newcommand{\Z}{\mathbb{Z}} 

\newcommand{\NN}{\mathbb{N}}

\newcommand{\LL}{{\rm{L}}}
\newcommand{\CC}{{\rm{C}}}
\newcommand{\HH}{{\rm{H}}}

\newcommand{\lt}{\left}
\newcommand{\rt}{\right}
\def\[{[\![}
\def\]{]\!]}

\newcommand{\ii}{{\rm{i}}}
\newcommand{\Id}{{\rm Id}}

\newcommand{\per}{{\rm per}}

\newcommand{\Q}{{\rm Q}}

\def\XXint#1#2#3{{\setbox0=\hbox{$#1{#2#3}{\int}$}
    \vcenter{\hbox{$#2#3$}}\kern-.5\wd0}}

\renewcommand{\epsilon}{\varepsilon}
\renewcommand{\tilde}{\widetilde}
\newcommand{\langl}{\lt\langle}
\newcommand{\rangl}{\rt\rangle}

\newcommand{\et}{\quad\text{and }}
\newcommand{\dans}{\quad\text{in }}
\newcommand{\si}{\quad\text{if }}

\newcommand{\pour}{\quad\text{for }}
\newcommand{\pourtout}{\quad\text{for all }}

\newtheorem{theorem}{Theorem}[section]

\newtheorem{lemma}[theorem]{Lemma}

\newtheorem{proposition}[theorem]{Proposition}
\newtheorem*{proposition*}{Proposition}

\newtheorem*{theorem*}{Theorem}
\newtheorem*{lemma*}{Lemma}
\newtheorem*{corollary*}{Corollary}
\newtheoremstyle{TheoremNum}
      {\topsep}{\topsep}              		
      {\itshape}                      		
      {}                              		
      {\bfseries}                     		
      {.}                             		
      { }                             		
      {\thmname{#1}\thmnote{ \bfseries #3}}	
\theoremstyle{TheoremNum}

\theoremstyle{remark}
\newtheorem{remark}{Remark}

\renewcommand{\P}{P}

\newcommand{\ueps}{u^\epsilon}
\newcommand{\ubar}{\bar{u}}

\newcommand{\aeps}{a_\epsilon}
\newcommand{\abar}{\bar{a}}

\newcommand{\LLb}{{\rm{\bf{L}}}}
\newcommand{\LLpot}{{\rm{\bf{L}}}_{\rm pot}}
\newcommand{\LLsol}{{\rm{\bf{L}}}_{\rm sol}}

\newcommand{\bu}{{\bm u}}
\newcommand{\bhu}{\hat{\bm u}}
\newcommand{\bhpsi}{\hat{\bm \psi}}
\newcommand{\bhphi}{\hat{\bm \phi}}
\newcommand{\bhzeta}{\hat{\bm \zeta}}

\newcommand{\bpsi}{{\bm \psi}}
\newcommand{\bha}{\hat{\bm a}}
\newcommand{\bhA}{\hat{\bm A}}
\newcommand{\ba}{{\bm a}}
\newcommand{\bhg}{\hat{\bm g}}
\newcommand{\bg}{{\bm g}}
\newcommand{\bhP}{\hat{\bm p}}

\newcommand{\bhJ}{\widehat{\bm{\mathcal{J}}}}
\newcommand{\bl}{\bm{l}}

\newcommand{\cbhP}{\hat{\rm{\bm p}}}
\newcommand{\cbhQ}{\hat{\rm{\bm P}}}

\newcommand{\cbf}{{\rm{\bm f}}}

\newcommand{\cbhA}{\hat{\rm{\bm A}}}

\newcommand{\cbhpsi}{\hat{\Psi}}

\title{
	Stable full-field simulation of a multiscale elliptic equation by means of Quantized Tensor Trains
	}
	
\author[1]{Marc Josien}
\author[1]{Anas El Hachimi}
\author[1]{Isabelle Ramière}
\affil[1]{CEA.DES.IRESNE.DEC, Cadarache F 13108 St Paul Lez Durance}

\begin{document}
 	
 	\maketitle

\begin{abstract}
	In this article, we design an original solver based on Quantized Tensor Trains (QTT) for linear elliptic equations with heterogeneous coefficient field, that allows for extremely fine meshes.
	It can achieve full-field simulations in dimensions $d=2$ and $d=3$ with a number of Degrees of Freedom (DoFs) up to $20$ \textit{orders of magnitude} beyond the classical solvers, recovering accurately the solution as well as its gradient in the $\LL^2$ norm.
	For treating such an enormous amount of data, the solver crucially relies on the exponential compression properties of QTTs.
    This significantly improves upon the existing literature.
    The main ingredient of the proposed solver consists in the introduction of a penalization term involving the Helmholtz--Leray projector in the equation governing the gradient unknown. For practical reasons related to the expression of the Helmholtz--Leray projector, the penalized equation is solved in Fourier space. The primal solution is then obtained from the gradient via the Green operator.
        
	A core property of the solver is that it is unconditionally stable with respect to the mesh size.
    Based on numerical evidence supported by mathematical analysis, we show that reliable gradients and solutions can be obtained, and guaranteed by the proposed \textit{a posteriori} error estimator.
    As an illustration, we successfully solve an elliptic equation in a microstructured material with up to $10^{37}$ virtual degrees of freedom in dimension $d=3$.
\end{abstract}

\section{Introduction}
 
\subsection{Motivation}

	The numerical simulation of multiscale systems in materials science aims at predicting the behavior of materials featuring heterogeneities such as microstructures or defects \cite{le_bris_systemes_2005, torquato_random_2002}.
	In many interesting cases, taking into account the phenomena at the microscopic scale is essential for accurate prediction of the large-scale behavior of the structure, which the engineer is ultimately interested in.
	This motivated the development of the homogenization theory, which ramified into mathematics, applied sciences and numerics (see for example  \cite{allaire_shape_2002, bensoussan_asymptotic_1979, jikov_homogenization_1994, tartar_general_2009, torquato_random_2002, bargmann_generation_2018}).
	The fundamental objective of this theory is to average out the small scale and simulate a simpler \textit{effective} system only at the large scale.
	In general, this averaging process involves solving a problem on a volume element, which is representative of the microstructure but far smaller than the whole structure:
	this is the so-called cell problem or \textit{Representative Volume Element} (RVE) method for periodic or random microstructures, respectively.
	
	This approach works particularly well when considering periodic materials;
	still, concrete applications rarely fall into the realm of periodic homogenization.
	This motivated various numerical approaches such as the Heterogeneous Multiscale Method \cite{engquist2007heterogeneous} or Finite Element Squared (FE$^2$) \cite{feyel2003multilevel}, just to name a few popular approaches.
	These are, however, complex and costly, requiring to design small-scale and large-scale problems that communicate to each other.
	Last, they rely on an underlying scale separation, which may not always hold when considering real materials: defects, interfaces, gradients of material properties may arise.
	
	The main practical motivation for homogenization-based numerical computations is the following: Relying on the current computational capabilities, a brute-force full-field simulation that takes into account all scales in a single simulation is apparently out of reach.
	For example, consider the mechanical Finite-Element simulation of a 3D structure of characteristic size $L$, made of a heterogeneous material, where $\epsilon \ll L$ is the fine scale of heterogeneities.
	This would typically involve $(L/h)^3$ elements, where $h \ll \epsilon$ is the discretization step.
	This easily leads to a number of Degrees of Freedom (DoFs) that exceeds the memory capabilities of the most modern supercomputers --not to mention the need for highly parallel solvers.
	
	The above observation is apparently obvious; it may actually be \textit{misleading}.
	Indeed, there is no need to store \textit{all} of the DoFs of the fields of interest; in many cases, it is sufficient to have their \textit{compressed} representation.
	In order to overcome the DoFs limitation in terms of memory, the compression should not come as a post-processing step; rather the original problem should be directly expressed and solved in the compressed format.
	This is the fundamental idea behind \textit{a priori} reduced-order modeling approaches such as Proper Generalized Decomposition (PGD) \cite{Chinesta2013}, which rely on low-rank separated approximations. 
	However, since PGD is essentially based on a canonical tensor decomposition, its efficiency may deteriorate for high-dimensional problems with complex multidimensional dependencies, due to the potentially large number of modes required to achieve accurate approximations.
	
	Tensor network representations \cite{hackbusch_tensor_2012, khoromskij_tensor_2018, bachmayr_low-rank_2023} provide an alternative low-rank structure based on factorizations, which offers improved numerical robustness in high-dimensional settings.
	Indeed, tensor representations provide both compact and expressive representations of structured data.
	There are various tensor representations, the efficiency of which depend on the structure of the data.
	Among the most prominent formats are the Canonical Polyadic decomposition, the Tucker decomposition  \cite{hackbusch_tensor_2012}, and the more recent Tensor Train (TT) decomposition \cite{oseledets_tensor-train_2011} and its Quantized variant \cite{khoromskij_tensor_2018} (QTT).
	In the next two sections, we justify by reviewing the literature that, in principle, QTTs may indeed achieve full-field simulations of elliptic equations in microstructured materials.
	
\subsection{Quantized Tensor Trains}
\label{Sec:QTT1}

QTTs were first successfully used in the community of Quantum Mechanics, in combination with an algorithm for finding the first eigenvalue of an operator.
This is the well-known Density Matrix Renormalization Group (DMRG) algorithm \cite{white1993density}.
In the early 2010s, the concepts percolated in the mathematical community, notably thanks to the work of Oseledets and Khoromskij  \cite{oseledets_tensor-train_2011, khoromskij_tensor_2018}.
QTTs were perceived from the beginning as a promising tool for building PDE solvers.
In particular, Khoromskij noticed that QTTs achieve exponential compression when discretizing various analytical functions on a Cartesian grid.
The efficiency of the compression is quantified by the maximal rank $r$ of the matrices constituting the QTT (see Section \ref{Sec:QTT} for a precise definition): the required storage scales quadratically in~$r$, but only \textit{logarithmically} in the discretization step.
This exponential compression is specific to QTTs, by contrast with most of the other tensorization techniques \cite{hackbusch_tensor_2012, bachmayr_low-rank_2023}, which typically require a memory footprint scaling linearly with the inverse of the discretization step.

QTTs now come with a furnished toolbox to perform linear algebra \cite{khoromskij_tensor_2018, waintal_who_2026}.
These operations manipulate QTTs, without ever decompressing them --which is often out of reach, exceeding by far the memory storage.
The basic operations such as addition and application of operators are straightforward \cite{waintal_who_2026}.
Then, solving linear systems or finding the maximal eigenvalue of operators relies upon algorithms deriving from the DMRG algorithm, the so-called (Modified) Alternating Least Square method (MALS or ALS) in the parlance of \cite{holtz_alternating_2012}, or the Alternating Minimal Energy (AMEn) \cite{dolgov_alternating_2014}.
Next, the discrete Quantum Fourier Transform (QFT) \cite{khoromskij_tensor_2018, chen_direct_2024} can be achieved, with a complexity exponentially small with respect to the classical Fast Fourier Transform (FFT).

The maximal TT-rank $r$ plays a central role not only in memory concerns, but also determines the efficiency of the algorithms applied on the QTTs --the complexities these algorithms scale algebraically in $r$.
Yet, since all algebraic operations typically increase $r$, re-compression is achieved by appealing to a Singular Value Decomposition adapted to TTs: the TT-SVD \cite{oseledets_tensor-train_2011}.
For turning a compressible high-dimensional tensor into a QTT, there is the Tensor Cross Interpolation (TCI) algorithm \cite{oseledets_tt-cross_2010}.

\begin{remark}[On the Quantum nature of QTTs]
	More generally, QTTs' exponential compression makes them fit for solving high-dimensional problems and overcoming the curse of dimensionality.
	The word ``Quantized'' is justified by the fact that the QTT structure looks like the qubits of a Quantum Computer.
	Hence, somehow, QTTs may be seen as a practical way to emulate a Quantum Computer on a classical one, with the limitation that it works only with lowly entangled states \cite{fernandez_learning_2024} --this means that the maximal TT-rank is reasonably small.
	\\
	Yet, we emphasize that QTTs are only \textit{Quantum-inspired} algorithms, and do not rely on the future development of a Quantum Computer.
	Interested readers may consult \cite{liu_towards_2024, deiml_quantum_2025} for formulating some problems related to thermomecanics on a Quantum Computer.
\end{remark}

\subsection{Review of QTT-based strategies for solving elliptic PDEs}\label{Sec:QTTPDE}

In the past years, QTT techniques have been explored for various types of PDEs (see for example \cite{ rakhuba_robust_2021, bachmayr_stability_2020, kornev_tetrafem_2024, niedermeier2026solving}).
We report here contributions to solving linear elliptic PDEs of the type
\begin{align}\label{E*}
	-\nabla \cdot \lt( a \nabla u\rt) = f,
\end{align}
in a domain of typical size $L=1$. Here, $a$ is a heterogeneous coefficient field that may feature a microstructure of characteristic scale $\epsilon \ll 1$.
We focus on the physical dimensions $d=2$ and $d=3$, with the objective of achieving a fine mesh size.

It was early remarked that the discrete Laplacian with finite differences can be represented as a QTT operator or a \textit{Matrix Product Operator} (MPO) of low rank \cite{oseledets2010approximation, khoromskij2011d, khoromskij_tensor_2018}.
However, even if, in principle, very fine mesh size can be attained thanks to the QTT compression,
the ill-conditioning of the discrete Laplace operator prevents from solving numerically the Laplace equation$-\Delta u=f$ on fine meshes \cite{rakhuba_robust_2021} by directly approximating the Laplacian operator, \textit{e.g.} with finite differences.
Indeed, the discrete Laplace operator, is of typical condition number of $h^{-2}$, where $h$ is the mesh size: if $h<10^{-8}$, this goes beyond the machine precision.

Rakhuba \cite{rakhuba_robust_2019, rakhuba_robust_2021} made an attempt to circumvent this difficulty by building analytically the inverse of the 1D Laplacian as an MPO.
He combined it with the decomposition of the $2$-dimensional Laplacian $\Delta = \Delta_{1{\rm D}} \otimes \Id + \Id \otimes \Delta_{1{\rm D}}$ in order to build ``derivative-free'' expression of the equation $-\Delta u = f$  (this generalizes to any dimension).
This strategy achieved the impressive mesh size of $h=10^{-12}$ on a unit cube in $2$ and $3$ dimensions.
Yet, it does not allow for non-constant coefficients $a$,  and does not approximate well the gradient $\nabla u$.
We underline that, if this strategy could retrieve accurately $\nabla u$, then, combining it with the Lippman-Schwinger approach \cite{schneider_review_2021} might tackle heterogeneous coefficients as well.

A different strategy was proposed in \cite{bachmayr_stability_2020} by Bachmayr and Kazeev: a preconditioned Finite Element (FE) QTT-based solver.
The authors formalized a notion of well-conditioning related to QTTs, which are sensitive to the error amplification that may arise from the underlying multiplicative decomposition (\textit{cf.} \eqref{QTT-0} below).
Equipped with this definition, they built a provably stable well-conditioned MPO to solve elliptic equations, that may attain the fine mesh size $h=10^{-15}$ while retrieving both $u$ and $\nabla u$.\footnote{Strictly speaking, how to obtain $\nabla u$ is not explained in \cite{bachmayr_stability_2020}, but it can be inferred, and the authors established $H^1$ convergence of their method, which guarantees its accuracy.}
However, as the authors emphasized, in practice, their approach is limited because of the large maximal TT-rank $r$ of their MPO: at least $r=1152$ for the constant coefficient $a=1$ in dimension $d=2$.
Moreover, this TT-rank depends multiplicatively in the TT-rank of the discretization of the coefficient $a$ and exponentially in the dimension $d$.
This prevents its practical use beyond very simple heterogeneous coefficients in dimension $d=2$ and makes the dimension $d=3$ out of reach.
Indeed, a TT-rank of $1000$ is already large, and TT-ranks higher than $10.000$ are unaffordable without HPC strategies.

In \cite{kazeev_quantized_2022}, the authors analyzed the homogenization theory through the lens of the QTTs.
In particular, they assumed that $a$ features a special compressible microstructure, based on interleaved scale-separated periodic patterns.
Then, they used homogenization techniques for establishing that there exists an efficient QTT representation of $u$ and $\nabla u$.
Nevertheless, the expressivity of Quantized Tensor Trains extends well beyond the framework of periodic homogenization.
These theoretical findings are illustrated by numerical examples.
However, their numerical strategy is not described there, and delayed to a forthcoming article, which has not been published yet.
Thus, we cannot discuss it.

Last, but not least, there had been attempts in the community of FFT-based solvers to build QTT-based solvers.
These solvers rely on an explicit expression of the Green operator for the Laplacian, which is usually computed in the Fourier variable, for solving thermomechanics problems in microstructured media \cite{schneider_review_2021} in periodic cubes with a Cartesian mesh.
A first try was done by Risthaus and Schneider in \cite{risthaus_imposing_2024}, by appealing to QTTs for building the Green operator, but no real superiority over the usual FFT approach was observed for their strategy; this was related to the fact that the authors' problem was not expressed exclusively as QTTs --thus, decompression was unavoidable.
In \cite{hauck_sfft-based_2024}, the authors proposed to replace the FFT by the QFT.
They reported some advantages for their QTT-based strategy in terms of memory and time, but their discretization steps remain roughly of the same order as for classical FFT-based solvers.
We suspect that this might be related to their use of Lippman-Schwinger iterations, which are quite computationally expensive, since the whole iterative process involves numerous applications of MPOs followed by compressions.

\subsection{Our contribution}

In this article, we design an original general-purpose QTT-based solver for the elliptic equation \eqref{E*} on extremely fine meshes in dimensions $d=2$ and $d=3$, with the aim of obtaining accurate $\LL^2$ approximations of $u$ and $\nabla u$.
This solver allows for considering heterogeneous coefficient fields provided that these are compressible in the sense of QTTs.
The proposed solver is agnostic with regards to the (potential) homogenization process and has no prior analytical knowledge of the structure of the solution.
From a user's perspective, it looks like a brute-force full-field PDE solver for multiscale linear elliptic equations working on a desktop computer.
As such, our method represents a qualitative improvement over \cite{rakhuba_robust_2021} by allowing heterogeneous coefficients and enabling the accurate determination of the gradient~$\nabla u$.
Also, it quantitatively improves upon \cite{bachmayr_stability_2020} by significantly reducing the ranks of the involved MPOs, allowing more complex coefficient fields in dimension $d=2$ and for treating the dimension $d=3$.

After precisely settling our problem in Section \ref{Sec:Pb}, we detail in Section \ref{Sec:Algo} our numerical strategy.
It applies the QTT formalism to a discretized variational formulation of the problem~\eqref{E*} with the new unknown $\psi=\nabla u$.
To enforce the constraint that $\psi$ is a gradient, we employ the Helmholtz-Leray projector, which is naturally expressed in the Fourier variable.
The numerical method extensively uses the QTT toolbox for linear algebra.

Next, in Section \ref{Sec:NumAnalysis}, we analyze mathematically our numerical strategy.
In particular, we show that it is unconditionally stable in terms of the discretization step, enjoys an \textit{a priori} estimate, and an \textit{a posteriori} error estimator.
Stability is necessary for considering extremely fine meshes.
Furthermore, the \textit{a posteriori} error estimator is highly valuable in order to assess the quality of the output of our solver.

Last, in Section \ref{Sec:NumRes}, we perform numerical experiments with our method in dimensions $d=2$ and $d=3$.
First, we investigate the TT-rank of the Helmholtz-Leray projector.
Then, we validate the solver on two test-cases.
The first one involves a simple manufactured solution, where we benchmark our method against two other ones from the literature proposed in \cite{khoromskij_tensor_2018, bachmayr_stability_2020}.
The second test case is closer to real applications: it features a random microstructure that is periodized and modulated on large scale.
We use our solver with up to $10^{27}$ \textit{virtual} DoFs in dimension $d=2$ and up to $10^{37}$ \textit{virtual} DoFs in dimension $d=3$, and secure a relative error of $10^{-3}$ in the $\LL^2$ norm on the gradient $\nabla u$.
By \textit{virtual} DoFs, we emphasize these are indirectly represented through a compression, which adapts itself to the problem data.
To illustrate the fine discretization of this test-case: it would amount to solving a stationary heat equation inside a cube of \textit{heterogeneous} matter of $1$ m$^3$  with \textit{a discretization step below the atomic scale}.
By way of comparison, classical FFT solvers typically reach the memory limit when the number of DoFs exceeds $10^{10}$, \textit{cf.} \cite{schneider_review_2021}.

\section{Problem}\label{Sec:Pb}

We now settle a framework representative of stationary thermal diffusion in heterogeneous material --generalization to linear elasticity is straightforward.
We consider the following elliptic PDE in divergence form:
\begin{align}\label{E}
	- \nabla \cdot a \nabla u = \nabla \cdot g.
\end{align}
The coefficient $a$ stands for a the material law of a heterogeneous medium --for example, its thermal conductivity.
The equation \eqref{E} is set on the periodic cube\footnote{This is equivalent to having \eqref{E} satisfied on $[0,1]^d$, and assuming periodicity of $u$ and anti-periodicity of $a\nabla u \cdot n$ on the boundary of $[0,1]^d$, where $n$ is the outer normal.}  $\Q_{1,\per} :=\R^d / \Z^d$, for $d \in \NN$.
This domain has the beneficial property of having no boundary, allows for using the Fourier transform, and enjoys explicit formulas for the inverse Laplacian --these properties are exploited in FFT-based solvers \cite{schneider_review_2021}, and also in the QTT method \cite{rakhuba_robust_2021}.
Dealing with more complex domains and boundary conditions is beyond the scope of this article.
Notice that having $\nabla\cdot g$ in \eqref{E} as a source term is not restrictive compared to~\eqref{E*}; we refer to Remark \ref{Rk:f} for a precise explanation.
However, we stick to~\eqref{E}, which is conceptually simpler for our explanations.

We assume that:
\begin{enumerate}[label={(\roman*)}]
	\item\label{A-i} The heterogeneous coefficient $a \in \LL^\infty(\Q_{1,\per})^{d\times d}$ is elliptic, bounded, and symmetric. 
	Up to renormalizing it, we assume that there exist constants $\Lambda, \lambda>1$ such that the following inequality holds for any $x\in \Q_{1,\per}$:
	\begin{align}\label{A-i-eq}
		\lambda^{-1} \leq a(x) \leq \Lambda.
	\end{align}
	This inequality means that $\lambda^{-1}|\xi|^2 \leq \xi \cdot a(x) \xi \leq \Lambda |\xi|^2$ holds for any $\xi \in \R^d$.
	\item\label{A-ii} The source term $g$ satisfies  $g \in \LLb^2(\Q_{1,\per}) := \lt(\LL^2(\Q_{1,\per})\rt)^d$.
\end{enumerate}
Under Assumptions  \ref{A-i} and  \ref{A-ii}, thanks to the Fredholm alternative \cite[Th.\ 6.6]{brezis_functional_2011}, \eqref{E} has a unique solution $\nabla u \in \LLb^2(\Q_{1,\per})$.
The latter enjoys the following sharp estimate
\begin{align}\label{num:0001}
	\|\nabla u\|_2 \leq \lambda \|g\|_2.
\end{align}

Notice that $u$ itself is only unique up to the addition of a constant.
A usual additional constraint for restoring uniqueness on the level of $u$ is imposing the vanishing average
\begin{align*}
	\int u = 0.
\end{align*}
Hereinafter, the integrals are always implicitly on the domain $\Q_{1, \per}$.

\section{Principle of our numerical method}\label{Sec:Algo}

	Our approach involves the following main ingredients:
	\begin{itemize}
		\item the Helmholtz-Leray decomposition for reformulating \eqref{E} as a minimization problem on a new variable $\psi \in \LLb^2(\Q_{1,\per})$ generalizing $\nabla u$, under the constraint that $\psi$ is a gradient, that we relax as a penalization term,
		\item the use of the Fourier transform to express analytically the Helmholtz-Leray projector,
		\item the QTT representation for formulating the discretized problem, which allows for using extensively the available QTT toolbox.
	\end{itemize}
	Each of these ingredients is not original \textit{per se}: the first two are reminiscent of FFT-based solvers \cite{moulinec_fast_1994, schneider_review_2021}, whereas we use QTT tools and concepts from the literature \cite{khoromskij_tensor_2018, oseledets_tensor-train_2011, oseledets_tt-cross_2010, holtz_alternating_2012, chen_direct_2024, fernandez_learning_2024}.
	Nevertheless, the way they are assembled here seems particularly efficient, in the sense that it yields a well-conditioned and medium-rank QTT linear system.
	
	From an abstract point of view, a main merit of our approach is to maintain a clear correspondence between the structures of the continuous framework $\HH^1(\Q_{1, \per})$ and the discrete QTT framework, where \eqref{E} is originally set and then discretized, respectively.
	This correspondence is algebraic, by preserving the variational nature of \eqref{E} within the discretized framework; this is usual, \textit{e.g.} for Galerkin Finite Element methods \cite{ern_finite_2017}.
	But, as importantly, this correspondence is also topological: to the $\LLb^2$ norm on $\psi$ corresponds the $\ell^2$ norm on the QTT discretization of $\psi$.
	The first norm is naturally adapted to \eqref{E}, because it expresses naturally the continuity property $\|\nabla u\|_2 \lesssim \|g\|_2$.
	The second norm plays a crucial role for QTT algorithms (such as ALS); indeed, QTT approximations are naturally controlled in the Euclidean space $\ell^2$.
	\footnote{The TT-SVD algorithm \cite{khoromskij_tensor_2018}, which is at the core of most QTT approximation strategies, is almost-optimal for best-rank approximations of MPS with respect to the $\ell^2$ norm.}

	\subsection{Helmholtz-Leray reformulation}\label{Sec:HL}
	
		Equation \eqref{E} has an obvious variational structure (recall that $a$ is symmetric).
		However, we want to preserve the continuity property~\eqref{num:0001} in our numerical strategy, to accurately retrieve $\nabla u$.
		Since QTTs algorithms are naturally expressed in $\ell^2$ spaces (\textit{cf.} beginning of Section \ref{Sec:Algo}), we introduce a new variable $\psi = \nabla u \in \LLb^2(\Q_{1,\per})$.
		In this perspective, we use the Helmholtz-Leray decomposition.
		Following \cite[Chap.\ 1]{jikov_homogenization_1994}, we introduce the orthogonal decomposition
		\begin{align*}
			\LLb^2(\Q_{1,\per})  = \LLpot^2(\Q_{1,\per}) \oplus \LLsol^2(\Q_{1,\per}),
		\end{align*}
		where
		\begin{align*}
			 \LLpot^2(\Q_{1,\per}) := \{\psi \in \LLb^2(\Q_{1,\per}), \exists u \in \HH^1(\Q_{1,\per}), \psi = \nabla u \},
			 \quad \et  \quad
			 \LLsol^2(\Q_{1,\per}) := \{\psi \in \LLb^2(\Q_{1,\per}), \nabla \cdot \psi = 0\}.
		\end{align*}
		
		Equipped with this conceptual tool, we rewrite \eqref{E} as the minimization of the following functional
		\begin{align}\label{Def:J}
			\mathcal{J}(\psi) := \frac{1}{2} \int \psi \cdot a \psi + \int \psi \cdot g,
		\end{align}
		under the constraint $\psi \in \LLpot^2(\Q_{1,\per})$.
		
		To express more explicitly this constraint on $\psi$, we make use of the  Helmholtz-Leray orthogonal projector on $\LLsol^2(\Q_{1,\per})$.
		We denote it as the operator
		\begin{equation}\label{Eq:PGamma}
			P:=\Id+\Gamma, \qquad \pour \Gamma:=- \nabla \Delta^{-1}\nabla \cdot.
		\end{equation}
		In the above definition, $\Gamma$ is the Green operator that plays a central role for FFT-based solvers \cite{moulinec_fast_1994}.
		Obviously, the constraint $\psi \in \LLpot^2(\Q_{1,\per})$ is equivalent to $P\psi=0$.
		In the physical space, $P$ is represented as a convolution with a non-integrable distribution \cite{gilbarg_elliptic_2001}, which may cause numerical difficulties.
		By contrast,
		taking the Fourier transform, $\mathcal{F}(P)$ is simply a multiplier of symbol
		\begin{align}\label{Gamma}
			\hat{p}_{ij}(k) = \delta_{ij} - \frac{k_ik_j}{0^++|k|^2},
		\end{align}
		expressed in terms of the Fourier variable $k=(k_1,\dots,k_d) \in (2\pi\Z)^d$.
		Here, the notation $0^+$ only means that $\hat{p}_{ij}(k=0)$ is well-defined and equal to $\delta_{ij}$.
		Notice that the projector $P$ is self-adjoint (since its Fourier transform is real); it is thus indeed orthogonal.
		
		Taking the Fourier transform of \eqref{Def:J} and using the Parseval theorem, solving \eqref{E} amounts to minimizing the functional
		\begin{align}\label{E-2}
			\hat{\mathcal{J}}(\hat{\psi}) :=\frac{1}{2}\langl\hat\psi, \hat{a} * \hat \psi \rangl + \langl \hat\psi, \hat{g} \rangl,
		\end{align}
		for $\hat\psi \in (\ell^2((2\pi\Z)^d))^d$,
		under the constraint
		\begin{align}\label{Const}
			\hat{p}(k) \hat{\psi}(k) = 0 \pourtout k \in (2\pi\Z)^d.
		\end{align}
		In \eqref{E-2}, the symbol $*$ stands for the discrete convolution
		\begin{align*}
			(\hat u*\hat v)(k) = \sum_{k'} \hat u(k-k')\hat v(k),
		\end{align*}
		where the sum is over $k' \in (2\pi\Z)^d$,
		and where we use the usual Hermitian form:
		\begin{align*}
			 \langl \hat{u},\hat{v} \rangl = \sum_k \hat{u}^\dagger(k) \cdot \hat{v}(k),
		\end{align*}
		where $\dagger$ stands for the Hermitian transpose.
		Notice that, in \eqref{E-2}, a contraction is implicit, namely $(\hat{a}*\hat{\psi})_i = \sum_j \hat{a}_{ij}*\hat{\psi}_j$.
		
		We relax the constraint \eqref{Const} by penalizing the functional $\hat{\mathcal{J}}$ with an additional term
		\begin{align}\label{Jmu}
			\hat{\mathcal{J}}_\mu(\hat{\psi}) :=& \frac{1}{2}\langl \hat\psi, \hat{a}*\hat\psi \rangl + \langl \hat{\psi}, \hat{g} \rangl
			+ \frac{\mu}{2}\langl \hat\psi, \hat{p}\hat\psi\rangl,
		\end{align}
		for a fixed penalty parameter $\mu \gg 1$.
		(For more details on penalization, see \cite[Chapter~17]{nocedal2006numerical}.)
		Indeed, since $P$ is an orthogonal projector, we have
		$\langle \hat\psi, \hat{p}\hat\psi\rangle = 
		\|\hat{p}\hat\psi\|_2^2$,
		so that the third right-hand side term of \eqref{Jmu} is a penalization term minimizing the solenoidal part of $\psi$.
		Equivalently, we may write the Euler-Lagrange equation associated with the minimization as
		\begin{equation}\label{EL}
			\hat{a}* \hat\psi
			+ \mu \hat{p}\hat\psi
  			= -  \hat{g}.
		\end{equation}
		Equivalently, we may also write the minimization problem in the physical space on $\psi \in \LLb^2(\Q_{1,\per})$ itself as
		\begin{align}\label{Jmu_phys}
			\mathcal{J}_\mu(\psi) :=& \frac{1}{2}\int \psi\cdot a\psi + \int \psi \cdot g
			+ \mu \int \psi \cdot P \psi.
		\end{align}
		Nevertheless, the latter is less amenable to discretization because of the singularity of $P$.
		
		Formulation \eqref{Jmu} (or \eqref{EL}) has the beneficial property of combining the terms from the coefficient, the source term and the constraint --which is related to the differential operators-- in an \textit{additive} way,
		whereas the direct variational formulation of \eqref{E} naturally combines $a$ and the differential operators in a \textit{multiplicative} way --the same holds for the usual Lippman-Schwinger iterations used in FFT-based solvers \cite{schneider_review_2021}.
		Hence, formulation \eqref{Jmu} is better for QTT-based strategies: Indeed, the TT-rank of a sum of MPOs increases in an additive way with respect to the TT-ranks of each term, whereas the TT-rank of a product MPOs increases in a multiplicative way \cite{oseledets_tensor-train_2011}.

		\subsection{Recovering $u$ and $\nabla u$}\label{Sec:getu}
		
		Our numerical method aims primarily at recovering the minimizer $\hat\psi$ of \eqref{Jmu}.
		Then, the quantities of interest  $\psi=\nabla u$ and finally $u$ are easily recovered by a post-processing step.
		Clearly, $\psi$ is obtained from $\hat\psi$ by inverting the DFT.
		Once we have $\nabla u$, we may obtain $u$ as well by using the Green operator $Q : f \mapsto u$ for the Laplace equation $-\Delta u = \nabla \cdot f$, by setting $u=-Q(\nabla u)$.
		The operator $Q$ is expressed in the Fourier variable $k$ by the associated Fourier multiplier
		\begin{align}\label{Def:q}
			\hat{q}(k) = \frac{\ii k}{0^+ + |k|^2}.
		\end{align}
		Thus, we simply evaluate $\hat{u}=-\hat{q}\cdot\hat\psi$.
		
		We underline that, by contrast with most of the usual FE numerical methods, but similarly to some FFT-based methods \cite{schneider_review_2021}, we obtain $u$ from $\nabla u$, and not the other way round.
		Here, the main reason is that the latter method is unstable when reaching very fine discretization steps \cite{rakhuba_robust_2021}.

		\begin{remark}[Case of source term not in divergence form]\label{Rk:f}
			A similar procedure can be employed if, instead of \eqref{E}, we have to solve \eqref{E*}.
			Given $f$ of vanishing average in \eqref{E*}, then, by the Fredholm alternative \cite[Th.\ 6.6]{brezis_functional_2011}, we can define $g=-\nabla v$ for $-\Delta v = f$, so that $\nabla \cdot g=f$.
			Hence, we may replace $f$ in \eqref{E*} by a source term in divergence form and get back to \eqref{E}.
			In practice, we compute numerically $g$ from $f$ as a preprocessing step, using the following identity in the Fourier variable:
			\begin{align}\label{E:00}
				\hat{g} = -\hat{q}\hat{f},
			\end{align}
			where the Fourier multiplier $\hat{q}$ is defined by \eqref{Def:q}.
		\end{remark}
		
		\begin{remark}\label{Case:Id}
			When $a$ is equal to identity (and more generally for $a$ constant, up to a modification of $\Gamma$), we can directly apply $\Gamma$ and obtain $\nabla u= \Gamma g$ (the case $-\Delta u=f$ can be treated by Remark \ref{Rk:f}).
			As such, this improves upon the well-optimized QTT-based solver of \cite{rakhuba_robust_2021} for the constant-coefficient elliptic equation, which cannot retrieve accurately the gradient of the solution.
		\end{remark}
	
	\subsection{Discretization}
		
		By construction (see \cite{khoromskij_tensor_2018} and Section \ref{Sec:QTT} here), QTTs are adapted to the uniform Cartesian grid:
		\begin{align}\label{Def:G}
			G_L := \{0, 2^{-L}, \dots, 1-2^{-L}\}^d,
		\end{align}
		where $L$ stands for the discretization level; thus, the mesh size is $h=2^{-L}$. 
		This gives rise to the dual Fourier grid:
		\begin{align}\label{Def:hatG}
			\hat{G}_L := \{0, 2\pi, \dots, 2\pi (2^{L-1}-1), -2\pi2^{L-1}, \dots, ,-2\pi\}^d.
		\end{align}
		Hence, a function $u$ on the periodic cube $\Q_{1,\per}$ is naturally discretized as $\bu$ by evaluating it on each point of $G_L$.
		Similarly, its Fourier transform $\hat{u}$ is discretized as $\bhu$ by evaluating it on $\hat{G}_L$.
		The Discrete Fourier Transform (DFT) turns an object $\bu$ discretized on $G_L$ into an object $\bm{\mathcal{F}}\bu \simeq \bhu$ discretized on $\hat{G}_L$, but the commutation diagram involving the Fourier transform and the discretization is approximate:
		there are (controlled) errors.
		For simplicity, we neglect them here, and we only consider them explicitly in Section~\ref{Sec:NumAnalysis}.
		
		We discretize $\hat{\psi}$ on $\hat G_L$ by $\bhpsi$, and we define similarly $\bha$ and $\bhg$.
		Notice that $\bhpsi_k$ and $\bhg_k \in \R^d$, and $\bha_k \in \R^{d\times d}$ where $k \in \hat{G}_L$ is the discretized Fourier variable; in other words, $\bhpsi$ and $\bhg$ are vector-valued vectors, whereas $\bha$ is a matrix-valued vector.
		There are various ways of discretizing~$\hat{p}$: these techniques are called \textit{filtering} in the FFT solvers community \cite{schneider_review_2021}, and often involves the interpretation of the Fourier symbol $\ii k$ as a discrete derivative.
		The most straightforward choice is a \textit{spectral} discretization:
		\begin{align}\label{Pkc}
			\bhP_k := \hat{p}(k).
		\end{align}
		Nevertheless, we emphasize that our strategy also works with other types of filtering, such as that associated with the finite difference operator $D_h(u)=h^{-1}(u(\cdot+ h) - u(\cdot))$.
		
		Now, we can discretize \eqref{Jmu} as
		\begin{align}\label{E-disc}
			\bhJ_\mu(\bhpsi) := \frac{1}{2} \bhpsi^\dagger \cdot \bha * \bhpsi + \bhpsi^\dagger \cdot \bhg + \frac{\mu}{2} \bhpsi^\dagger \cdot \bhP \odot \bhpsi,
		\end{align}
		where we use the Hadamard product $(\bhP \odot \bhpsi)_k := \bhP_k \bhpsi_k$.
		The point of view of minimization is conceptually closer to the QTT algorithms we will use in Section \ref{Sec:MALS}, which rely on such a structure.
		Nevertheless, analogously to \eqref{EL}, minimizing the discretized functional $\bhJ_\mu(\bhpsi)$ is equivalent to solving the linear system
		\begin{align}\label{E-lin}
			(\bha * + \mu \bhP \odot) \bhpsi = - \bhg.
		\end{align}
		Here, we make use of the notation $\bha*$ and $\bhP \odot$, for the operators $\bhpsi \mapsto \bha *\bhpsi$ and $\bhpsi \mapsto \bhP \odot\bhpsi$, respectively.

	\subsection{Use of Quantized Tensor Trains}\label{Sec:QTT}
	
	\subsubsection{Basics of QTTs}\label{SimpleExample}
	
	For readers not familiar with QTTs, we briefly recall here some basics. 
	We refer to \cite{baker_methodes_2021, waintal_who_2026, khoromskij_tensor_2018} for a thorough introduction.
	
	TTs are compressed representations of high-order tensors $\bm{f}$ with $L$ indices as a product of matrices.
	This is called the Matrix Product State (or MPS) in the community of Quantum Mechanics:
	\begin{align}\label{QTT-0}
		\bm{f}_{x_1,\dots, x_L} = \lt(\cbf_1\rt)_{x_1} \cdots \lt(\cbf_L \rt)_{x_L},
	\end{align}
	where the indices $x_1,\dots, x_L$ are called \textit{physical} indices.
	In \eqref{QTT-0}, the matrices, or cores $(\cbf_{\ell})_{x_\ell}$ belong to $\R^{r_{\ell-1} \times r_{\ell}}$, where $r_\ell$ are called the TT-ranks, and, by convention $r_{0}=1=r_{L}$.
	Decomposition \eqref{QTT-0} is referred to as a Tensor Train and the word ``Quantized'' is added when the indices $x_1, \dots, x_L$ belong to $\{0, 1\}$, which we assume hereafter.
	Indeed, in this case, the entries $x_\ell$ in \eqref{QTT-0} may be seen as qubits \cite{waintal_who_2026}.
	Storing the whole tensor $F$ requires storing its $2^L$ entries, which may exceed by far the computer memory even for moderately large $L$ (for example $L > 40$).
	By contrast, the required storage of all the constituents of the right-hand side of \eqref{QTT-0} is of less than $2L r^2$, where we define the maximal TT-rank $r:=\max_\ell r_\ell$.
	The latter quantifies the efficiency of the QTT-decomposition, \textit{cf.} Section \ref{Sec:QTT1}.
	
		We show here how to represent a simple function $f:[0, 1] \to \R$ with a given discretization step $h=2^{-L}$ as an MPS.
		We use the dyadic expansion of $x \in [0, 1]$ in the form of
		\begin{align*}
			x = \sum_{\ell=1}^L 2^{-\ell} x_\ell \qquad \pour x_\ell \in \{0, 1\},
		\end{align*}
		and discretize $f$ as $\bm{f}_{x_1, \dots, x_L} = f(x)$.
		As an example \cite{khoromskij_tensor_2018}, if $f(x)=\exp(x)$, we may represent $\bm{f}$ as an MPS of TT-ranks all equal to $1$ as follows:
		\begin{align*}
			\bm{f}_{x_1,\dots, x_L} = \exp(2^{-1}x_1) \exp(2^{-2}x_2) \cdots \exp(2^{-L}x_L).
		\end{align*}
		Indeed, we identify in \eqref{QTT-0} the scalar matrices $\lt(\cbf_\ell\rt)_{x_\ell} := \exp(2^{-\ell}x_\ell) \in \R^{1\times 1}$.
	
	\subsubsection{Formats for the MPS representation of $d$-dimensional functions}\label{Sec:Format}
	
		There are many ways to generalize the example of Section \ref{SimpleExample} in dimension $d>1$, \textit{cf.} \cite{bachmayr_stability_2020, rakhuba_robust_2021}.
		We explain here how to proceed in dimension $d=2$ (see Remark \ref{Rk:3D} below for $d=3$).
		We associate the grid points $(x,y) \in G_L$ and $(k^x, k^y) \in \hat{G}_L$ with their dyadic decompositions:
		\begin{align*}
			x = \sum_{\ell=1}^L 2^{-\ell} x_\ell
			\pour x_\ell \in \{0, 1\}
			,
			\qquad 
			\et k^x = \lt\{
			\begin{aligned}
			&2\pi \sum_{\ell=2}^L 2^{\ell} k^x_\ell  && \si k^x_1 =0,
			\\
			&-2\pi2^{L-1} + 2\pi \sum_{\ell=2}^L 2^{\ell} k^x_\ell  && \si k^x_1 =1,
			\end{aligned}
			\rt.
			\pour k^x_\ell \in \{0, 1\},
		\end{align*}
		and similarly for $y$ and $k^y$.
		Thus, we need to decide how to order the indices $x_1, \dots, x_L, y_1, \dots, y_L$ in the MPS $\bpsi$, and similarly for its Fourier transform $\bhpsi$ with respect to the indices $k^x_1, \dots, k^x_L, k^y_1, \dots, k^y_L$.
		This choice is not innocuous, because, as exemplified in Section \ref{Sec:Green} for $\hat{p}$, it determines the TT-ranks of $\bpsi$ for a fixed accuracy.
		To obtain an efficient MPO representation of the QFT, we need to align the $k^x$ and $k^y$ with $x$ and $y$, respectively, but by reversing the order \cite{khoromskij_tensor_2018, chen_direct_2024}.
		Hence, taking into account the symmetries, there are four reasonable choices, which are written in Table \ref{Tab:Format2}.
		
		\begin{table}[h]
			\begin{center}
				{\renewcommand{\arraystretch}{1.5}
			\begin{tabular}{|l|l|l|l|}
				\hline
				Designation & Ordering (physical space)
				& 
				Ordering (Fourier space)
				\\
				\hline				
				\texttt{x1x2\_y1y2}
				& $x_1, \dots, x_L, y_1,\dots, y_L$
				& $k^x_L, \dots, k^x_1, k^y_L,\dots, k^y_1$
				\\
				\hline
				\texttt{x1x2\_y2y1}
				& $x_1, \dots, x_L, y_L,\dots, y_1$
				& $k^x_L, \dots, k^x_1, k^y_1,\dots, k^y_L$
				\\
				\hline
				\texttt{x2x1\_y1y2}
				& $x_L, \dots, x_1, y_1,\dots, y_L$
				& $k^x_1, \dots, k^x_L, k^y_L,\dots, k^y_1$
				\\
				\hline				
				\texttt{x1y1}
				& $x_1, y_1, \dots, x_L, y_L$
				& $k^x_L, k^y_L, \dots, k^x_1, k^y_1$
				\\
				\hline
			\end{tabular}
			}
			\end{center}
			\caption{Different QTT formats in dimension $d=2$.}
			\label{Tab:Format2}
		\end{table}
		
		The first three formats of Table \ref{Tab:Format2} keep the dyadic components of the $x$ variable stand together and concatenate them with those of $y$.
		The nuance is how to glue these two groups, either gluing them in the same order for \texttt{x1x2\_y1y2}, or on the smallest scales for \texttt{x1x2\_y2y1} and on the largest scale for \texttt{x2x1\_y1y2}, by mirroring the indices.
		These three formats are compatible with the structure of the Laplacian $\Delta = \Delta_{1{\rm D}} \otimes \Id + \Id \times \Delta_{1{\rm D}}$, \textit{cf.} \cite{rakhuba_robust_2021}.
		\\
		Of different nature is the interleaved format \texttt{x1y1}, which sorts the dyadic components by scales. This format was successfully used in \cite{bachmayr_stability_2020}, and seems \textit{a priori} more natural for problems featuring a scale separation.
		Indeed, solutions of these problem are approximated by Ansätze considering the separated scales as independent variables, which may be interpreted as tensorizations across these scales (see \cite{allaire_shape_2002} for the two-scale expansion).

		Now, recall that $\hat\psi$ is not a scalar field, but a vector field.
		Thus, we represent $\hat\psi$ as a high-order vector-valued tensor $\bhpsi_{\bl} \in \R^d$ depending on  the multi-index $\bl=(l_1,\dots, l_{dL}) \in \{0, 1\}^{dL}$ --the interpretation of $\bl$ depends on the chosen format, \textit{cf.} Table \ref{Tab:Format2}.
		We decompose it as a TT as follows:
		\begin{align}\label{MPS-form}
			\lt(\bhpsi_{\bl}\rt)_i =  
			\lt(\cbhpsi_0\rt)_i \lt(\cbhpsi_1\rt)_{l_1} \cdots \lt(\cbhpsi_{dL}\rt)_{l_{dL}}.
		\end{align}
		where we prepend a special core $\cbhpsi_0$ of physical index $i \in \{1, \dots, d\}$.
		In this representation, the first core for $\ell=0$ is related to the vector values of $\hat\psi=(\hat\psi_i)_{i \in \{1, \dots, d\}}$, while the other cores for $\ell \in \{1,\dots, dL\}$ represent the dependence on scales of the Fourier variables, \textit{e.g.} $(k^x, k^y):=k$ in dimension $d=2$.

		\begin{remark}[Generalization to $d=3$]\label{Rk:3D}
			Formats \texttt{x1x2\_y1y2} and \texttt{x1y1} have natural generalizations to higher dimensions.
			By contrast formats \texttt{x1x2\_y2y1} and \texttt{x2x1\_y1y2} are naturally generalized as Tensor Trees rather than TT in higher dimensions \cite{rakhuba_robust_2021}.
			Yet, since we seek to avoid this more complex structure, we simply mirror once more the components, that is using $x_1,\dots,x_L, y_L, \dots, y_1, z_1, \dots, z_L$ and $x_L,\dots,x_1, y_1, \dots, y_L, z_L, \dots, z_1$ conventions, respectively.
			By a slight abuse, we consider decomposition \eqref{MPS-form} as an MPS, although the cardinal of the set of the first index is now equal to $3$ in dimension $d=3$.
		\end{remark}
		
		\subsubsection{The structure of the MPOs $\bhP\odot$ and $\bha*$}

		The MPOs $\bhP\odot$ and $\bha*$ are discretized in a similar manner as \eqref{MPS-form} by doubling the physical indices, namely
		\begin{align}\label{Eq:p}
			\lt(\lt(\bhP\odot\rt)_{\bl \bl'}\rt)_{ii'} =  
			\lt(\cbhQ_0\rt)_{ii'} \lt(\cbhQ_1\rt)_{l_1 {l_1}'} \cdots \lt(\cbhQ_{dL}\rt)_{l_{dL} {l_{dL}}'}
			\quad \et \quad
			\lt((\bha *)_{\bl \bl'}\rt)_{ii'} =  
			\lt(\cbhA_0\rt)_{ii'} \lt(\cbhA_1\rt)_{l_1 {l_1}'} \cdots \lt(\cbhA_{dL}\rt)_{l_{dL} {l_{dL}}'},
		\end{align}
		where $\bl=(l_1,\dots, l_{dL}) \in \{0, 1\}^{dL}$, $\bl'=(l_1',\dots, l_{dL}') \in \{0, 1\}^{dL}$, $i, i' \in \{1, \dots, d\}$.
		By definition of the Hadamard product, we have that $\lt((\bhP\odot)_{\bl, \bl'}\rt)_{ii'}$ is diagonal in $(\bl, \bl')$, so that expression~\eqref{Eq:p} can be simplified. Indeed, we may expand for $\ell > 0$
		\begin{align*}
			\lt(\lt(\cbhQ_\ell\rt)_{l_\ell {l_\ell}'}\rt)_{\alpha\beta} = \lt((\cbhP_{\ell})_{l_\ell}\rt)_{\alpha\beta}  \delta_{l_\ell {l_\ell}'},
		\end{align*}
		where there is no implicit sum, and	where the cores $\cbhP_{\ell}$ form the MPS\footnote{Once more, by a slight abuse, we consider the decomposition \eqref{num:0010} as an MPS, although the first core is not of the desired type; still, merging the indices $i$ and $i'$ into a single one, it remains a TT, so that all the TT-algorithms can be applied on it.} $\bhP$ that discretizes the matrix-valued field $k\mapsto \hat{p}(k)$:
		\begin{align}\label{num:0010}
			\lt(\bhP_{\bl}\rt)_{ii'}
			=  
			\lt(\cbhP_0\rt)_{ii'} \lt(\cbhP_1\rt)_{l_1} \cdots \lt(\cbhP_{dL}\rt)_{l_{dL}}.
		\end{align}
		Hence the MPO $\bhP\odot$ is of same rank as the MPS $\bhP$.
		
		By contrast, generically in $\bl' \neq \bl$, we have
			$(\bha *)_{\bl, \bl'} \neq 0$.
		We implement the discrete periodic convolution by appealing to the ``magical tensor'' \cite{waintal_who_2026}.
		This gives rise to a ternary operator $\bm{T}_{kk'k''}$ such that
		\begin{align*}
			\sum_{k'} \bm{T}_{kk'k''} \bha_{k'} = (\bha*)_{k,k''}.
		\end{align*}
		The operator $\bm{T}$ is of maximal TT-rank $r =2^d$ in format \texttt{x1y1} and $r =2$ in other formats of Table \ref{Tab:Format2}.
		Hence, using the above identity, we build an MPO $\bha*$ that is of TT-rank at most $r$ times larger than the TT-rank of the MPS $\bha$.

	\subsubsection{Building the QTT representation}\label{Sec:TCI}
	
		Concerning the TT-ranks that we may expect, the situation is contrasted between  first $\bhP$, and second $\ba$ and $\bg$.
		On the one hand $\bhP$ only depends on the discretization step.	By experimental evidence, we show in Section \ref{Sec:Green} that it can be represented by a medium-rank MPS with high precision.
		This result is related to previous observations about QTT representations of some analytical functions  \cite{khoromskij_tensor_2018}.
		\\
		On the other hand, in full generality, $\ba$ and $\bg$ cannot be represented by low-rank MPS: for example, if $a$ is generated from a Gaussian field with a small correlation range, one cannot hope for representing it with an MPS of small TT-ranks.
		On the contrary, there is the need for an underlying structure with sparse information for obtaining compression.
		Nevertheless, this compression arise in many case of interest.
		For example for periodic homogenization, but also in multiscale cases generalizing this restricted framework \cite{kazeev_quantized_2022}.
		When dealing with microstructured random media, inspired by \cite{kazeev_quantized_2022}, we advocate for an approach which is somehow related to the so-called Representative Volume Element (RVE) \cite{hill_elastic_1963, anantharaman_introduction_2012}; yet, instead of using the usual upscaling strategy, we may use this RVE to directly pave the domain and perform a full-field resolution of \eqref{E} enabled by our numerical approach.
		QTTs have the advantage that they allow for modulating this RVE by various ways --we exemplify a large-scale modulation in Section \ref{Sec:NumRes}.

		We apply the TCI algorithm \cite{oseledets_tt-cross_2010} for revealing the QTT structure of $\ba$, $\bg$, and $\bhP$.\footnote{Notice that $a$ and $g$ are usually expressed in physical variables, whereas $\hat{p}$ is given in Fourier variables by \eqref{Gamma}.}
		This algorithm can be seen as a black-box tool from linear algebra which seeks for the lowest-rank MPS approximations of a high-dimensional tensor~$Q$.
		It only samples the latter on well-chosen sets of points, but does not require to build the large original tensor~$Q$ --which is not feasible if $Q$ is a large high-dimensional tensor.
		Even if TCI seeks for an MPS representation of the lowest possible TT-rank while being close to $Q$ in the Frobenius norm up to a given accuracy, we underline that there is no guarantee that the output approximates well the initial tensor $Q$ --the zone where the error is large might not be sampled.
		Yet, in Section \ref{Sec:NumRes} below, we check this \textit{a posteriori} by estimating the $\ell^2$ norm of the approximation error by the Monte-Carlo method.
		From a user's perspective, TCI does not differ much from machine learning tools, even though its principles are fairly different.
		
	\subsubsection{The QFT MPO} 
		The QFT, which approximates the one-dimensional DFT denoted as $\bm{\mathcal{F}}_{1{\rm D}}$, is obtained as in \cite{chen_direct_2024} by means of high-precision interpolation.
		For expressing the $d$-dimensional QFT denoted here $\bm{\mathcal{F}}$, we use the following identity:
		\begin{align}\label{Fourier}
			\bm{\mathcal{F}} = 
			(\bm{\mathcal{F}}_{1{\rm D}} \otimes \Id \otimes \cdots \otimes \Id)
			(\Id \otimes\bm{\mathcal{F}}_{1{\rm D}} \otimes \Id \otimes \cdots \otimes \Id) \cdots 
			(\Id \otimes \cdots \otimes \Id \otimes \bm{\mathcal{F}}_{1{\rm D}}) = \bm{\mathcal{F}}_{1{\rm D}} \otimes \bm{\mathcal{F}}_{1{\rm D}} \otimes \cdots \otimes \bm{\mathcal{F}}_{1{\rm D}},
		\end{align}
		where the $d$ terms of the above tensor products are related to the $d$ coordinates of the physical variable (versus those of the Fourier variable).
		This structure is compatible with the formats in Table \ref{Tab:Format2}. 
		For \texttt{x1x2\_y1y2}, \texttt{x1x2\_y2y1}, \texttt{x2x1\_y1y2}, the MPO $\bm{\mathcal{F}}$ can directly be obtained by concatenating (if necessary, swapping) MPO representations of $\bm{\mathcal{F}}_{1{\rm D}}$.
		For format \texttt{x1y1}, a direct application of the above formula yields a TT-rank of $\bm{\mathcal{F}}$ depending exponentially in the dimension $d$.
		In such case, the leftmost factorized expression in \eqref{Fourier} seems more interesting, with TT-rounding after applying each factor.
			
	\subsubsection{Solving the linear system}\label{Sec:MALS}
	
		Once all the terms of the linear system \eqref{E-lin} are defined, it remains to solve it; in this task, we use an ALS-based solver \cite{holtz_alternating_2012}.
		ALS-based solvers are iterative strategies that alternatively optimize each core of the unknown MPS going back and forth along core indices (this back-and-forth movement is called ``sweep'').
		Each core optimization amounts to minimizing a quadratic functional involving a symmetric positive definite operator, or equivalently to solving the associated linear system --which we call the \textit{local problem}.
		By \cite[Th.\ 4.1]{holtz_alternating_2012}, the condition number of the linear operator $(\bha * + \mu \bhP\odot)$ is an upper bound for the condition number of the matrices involved in the local problems in the ALS strategy.

		For ALS-based solvers, of paramount importance for the efficiency is the maximal TT-rank of the operator.
		Here, since we have decoupled its two constituents $\bha *$ and $\mu \bhP\odot$ in an additive way, their TT-ranks sum up to the TT-ranks of the total operator.
		This has to be compared with \cite{bachmayr_stability_2020}, where the differential operators and the coefficient $a$ are combined in a multiplicative way, resulting in a worse maximal TT-rank.
		
	\section{Numerical analysis}\label{Sec:NumAnalysis}
	
	In this section, using classical variational tools, we derive three results, that can be summarized as follow:
	\begin{itemize}
		\item The condition number of the operator involved in \eqref{E-lin} depends only on the penalty parameter $\mu$ and the local bounds on $a$, \textit{cf.} \eqref{A-i-eq}.
		\item \textit{A priori} discretization error estimates can be obtained, comparing a suitable interpolation of the minimizer $\bhpsi^*_{\mu, L}$ of $\bhJ_\mu$ defined by \eqref{E-disc}, and the solution $\nabla u$ of \eqref{E}. It scales linearly in $\mu^{-1}$, and in the $\LL^\infty$ discretization error made on $a$ and $\hat\P$, and features a spectral convergence in the discretization step $2^{-L}$.
		\item A sharp \textit{a posteriori} error estimator is derived, using the operators $P$ and $\Gamma$ given by \eqref{Eq:PGamma}.
	\end{itemize}
	
	For efficiency, and since the arguments rely on classical variational tools, we state and comment our main results, and postpone the proofs to the Appendix \ref{App:proofs}.
	
	\subsection{Condition number}\label{Sec:Cond}
	
	\begin{proposition}\label{Prop:Cond}
		The operator $\hat{a} * + \mu \hat{p}$ appearing in \eqref{Jmu} is symmetric positive definite.
		Moreover, its condition number satisfies
		\begin{align}\label{Cond}
			C \leq \lambda (\mu+ \Lambda).
		\end{align}
	\end{proposition}
	
	Similarly, up to the QTT approximation errors on $\bha$ and $\bhP$ due to our use of TCI, we have that the condition number of $(\bha * + \mu \bhP\odot)$
	is independent of the discretization step, and bounded as in \eqref{Cond}.	
	
	\subsection{\textit{A priori} error estimator}\label{Sec:DiscErr}
	
	In the sequel, we set a discretization level $L$ and we give ourselves $\bha$, $\bhP$, $\bhg$ on the grid $\hat{G}_L$.
	We want to estimate the discretization error between the minimizer $\bhpsi^*_{\mu, L}$ of $\bhJ_\mu$ defined by \eqref{E-disc}, and the solution $\nabla u$ of \eqref{E}.
	To do so, for given $\bhphi$ defined on the grid $\hat{G}_L$, we define its interpolate $\phi$ as
	\begin{align}\label{Interpol}
		\phi(x) = \frac{1}{2^{Ld/2}}\sum_{k\in \hat{G}_L} \bhphi_k \exp(-\ii k \cdot x).
	\end{align}
	
	We first study the influence of the parameter $\mu$:
	\begin{lemma}\label{Prop:Mu}
		Assume that $\psi \in \LLb^2(\Q_{1, \per})$ minimizes $\mathcal{J}_\mu$ defined in \eqref{Jmu_phys}, and that $\nabla u$ minimizes $\mathcal{J}$.
		Then, we have the following estimate:
		\begin{align}\label{num:0005}
			\|\psi - \nabla u\|_2
			\leq \frac{(\lambda+\Lambda)(1+\lambda\Lambda)}{\mu+\lambda^{-1}} \|g\|_2.
		\end{align}
	\end{lemma}
	
	Before stating the main result, we define by duality the operator $\bhA$ so that, for any $\bhzeta$, $\bhphi$ on $\hat{G}_L$, we have
	\begin{align}\label{astar}
		\bhzeta^\dagger \cdot \bhA \bhphi = \int\zeta \cdot a \phi,
	\end{align}
	where $\zeta$ and $\phi$ are built from $\bhzeta$ and $\bhphi$ following the interpolation \eqref{Interpol}.
	We define the symmetric residual operator $\delta \bhA$ by duality
	\begin{align*}
		\bhzeta^\dagger \cdot \delta \bhA \bhphi = \bhzeta^\dagger \cdot \bha* \bhphi - \int\zeta \cdot a \phi,
	\end{align*}
	and the residues $\delta \bhP$, $\delta \bhg$ on $\hat{G}_L$ as
	\begin{align*}
		\delta \bhP_k :=\bhP_k -  \hat{p}(k) \et 
		\delta \bhg_k :=\bhg_k -  \hat{g}(k),
	\end{align*}
	and we assume that $\hat{p}(k)$ is a $d\times d$ Hermitian matrix.
	We set $\nu_a$ the spectral norm $\|\delta \bhA\|_{2,2}$, and
	\begin{align*}
		\nu_P := \|\delta \bhP\|_{\ell^\infty}, \qquad \et
		\nu_g := 2^{-Ld/2}\|\delta \bhg\|_{\ell^2}.
	\end{align*}
	
	Since we are ultimately interested in oscillating coefficients (and possibly source terms), we assume that $a, g$ are regular from a fixed scale~$\epsilon>0$ downwards.
	A practical way to quantify it is to assume that $a, g \in \CC^{p}$ satisfy
	\begin{align}\label{As-0}
		\sum_{p'=0}^p \epsilon^{p'} \lt( \|\nabla^{p'} a\|_{\LL^\infty}+ \| \nabla^{p'} g \|_{\LL^2}\rt) \leq C_p,
	\end{align}
	where $p \in \mathbb{N}$ is fixed.
	
	\begin{proposition}\label{Prop:disc}
 		Let $p \in \mathbb{N}$.
		We make the following two assumptions:
		\begin{enumerate}
			\item[(i)]
			$a, g \in \CC^{p}$ satisfy \eqref{As-0}, for a given constant $C_p$.
			\item[(ii)] $\lambda\nu_a + \lambda\mu \nu_P \leq 1/2$.
		\end{enumerate}
		Then, the following error estimate holds:
		\begin{equation}\label{Disc:error2}
			\lt\|\psi^*_{\mu, L} -\nabla u\rt\|_2 \leq C
			\lt( (2^{L}\epsilon)^{-p} + \mu^{-1} + \nu_a + \mu \nu_P + \nu_g\rt),
		\end{equation}
		where the above constant $C$ depends only on $p$, $\lambda$ and $\Lambda$ in \eqref{A-i-eq}, and $C_p$ in \eqref{As-0}.
	\end{proposition}
	
	In \eqref{Disc:error2}, the ratio $\epsilon / 2^{-L}$ appears, which is typical of numerical analysis in multiscale elliptic equation.
		
	\subsection{\textit{A posteriori} error estimator}
	
	\begin{proposition}\label{Prop:Aposteriori}
		Let $\nabla u$ be the minimizer of $\mathcal{J}$ defined by \eqref{Def:J}.
		Let $\phi \in \LLb^2(\Q_{1,\per})$ be given, and define
		\begin{align}
			\label{Def:EP}
			&E_P[\phi] := \frac{\|P\phi \|_2}{\|\nabla u\|_2},
			\qquad \et \quad E_\Gamma[\phi] := \frac{\|\Gamma \big(a \phi + g\big)\|_2}{\|\nabla u\|_2}.
		\end{align}
		Then, we have the upper and lower estimates
		\begin{align}\label{Error-total}
			\frac{\|\nabla u - \phi\|_2}{\|\nabla u\|_2}
			\leq
			(1+\lambda \Lambda ) E_P[\phi] + \lambda E_\Gamma[\phi] =:E_{\max}[\phi],
		\end{align}
		and
		\begin{align}\label{Error-total2}
			\frac{\|\nabla u - \phi\|_2}{\|\nabla u\|_2} \geq \max\lt( E_P[\phi]; \frac{1}{ \Lambda \sqrt{2}}E_\Gamma[\phi]\rt)
			=:E_{\min}[\phi] \geq  \frac{1}{2\sqrt{2}(1+ \lambda \Lambda)}E_{\max}[\phi]
		\end{align}
	\end{proposition}
	
	Up to a renormalization, error $E_P[\phi]$ is the distance of $\phi$ to $\LLpot^2$, and $E_\Gamma[\phi]$ is the distance of $(a\phi + g)$ to $\LLsol^2$.
	The first one quantifies how $\phi$ is close to being a gradient, whereas the second one measures the deviation to the identity
	$-\nabla\cdot(a\phi+g)=0$.
	
	Notice that the estimator $E_{\max}$ is sharp in the sense it is, up to a constant, a lower bound of the error, \textit{cf.} \eqref{Error-total2}.
	Naturally, the expressions of the \textit{a posteriori} estimators do not involve the penalty parameter $\mu$.
	(However, the dependence on $\mu$ appears when evaluating the estimators on the minimizer $\psi$ of $\mathcal{J}_\mu$.)
	Notice that these above estimators do not cover the discretization errors, and may themselves suffer from numerical error when used in practice --because they require approximating $P$, $\Gamma$, and $a$.
	
	\section{Numerical results}\label{Sec:NumRes}
	
	We now expose our numerical results.
	First, we build a medium-rank accurate MPS $\bhP$ that approximates the symbol of the Helmholtz-Leray projector $\hat{p}$.
	We study the impact of the QTT-format (\textit{cf.} Table \ref{Tab:Format2}) on the TT-rank of $\bhP$.
	Then, we apply our strategy for retrieving a manufactured solution and challenge it with two other QTT-based methods.
	Last, we validate our strategy on a ``realistic'' problem that can as well be solved by the homogenization approach, in dimensions $d=2$ and $d=3$.
	This final problem features a random microstructure that is periodized on small scales and then modulated on large scales.

	In this task, we use our code Sisyphe, an in-house code written in Julia and based on the libraries ITensors, ITensorMPS, and Tensor4All.
	For QTT basic manipulations (additions, MPO applications, TT-rounding), we mostly rely on ITensors and ITensorMPS  \cite{fishman_itensor_2022}.
	We also use their 1-site version of the ALS algorithm \cite{holtz_alternating_2012} --but other versions such as 2-site or AMEn \cite{dolgov_alternating_2014} might be used as well.
	Then, we rely on Tensor4All \cite{fernandez_learning_2024}, which implements the QFT of \cite{chen_direct_2024}
	for computing $\bm{\mathcal{F}}_{1{\rm D}}$.
	Last, we use the TCI implementation of Tensor4All.
	
	Sisyphe contains all the PDE-related features; it not only encompasses the numerical strategy explained in the previous sections, but also the direct finite-difference based strategy \cite{oseledets2010approximation,khoromskij2011d, khoromskij_tensor_2018} and the method from \cite{bachmayr_stability_2020}, which we use for comparison.
	Sisyphe supports limited parallelism with shared memory; hence, we use it on a single CPU node, using up to $20$ CPU cores.
	Since the code is not carefully optimized, we do not insist much on the computation time, which ranges from minutes to tens of hours in the longest tests.

	\subsection{QTT approximation for $\bhP$}\label{Sec:Green}
	
	\subsubsection{Methodology}
	
	We investigate the question of how to numerically approximate $\bhP$ defined by \eqref{Pkc} by an MPS.
	Obviously, the aim is to obtain the best relative accuracy with the smallest TT-rank.
	We restrict ourselves to the physically relevant dimensions $d=2$ and $d=3$.
	To build the MPS, we appeal to TCI.
	To assess the quality of approximation, we measure the relative accuracy through the relative $\ell^2$ norms of the difference between the output of TCI and the original analytic one $\hat{p}$.
	This measure is empirical, and performed by a Monte-Carlo method with $1000$ sample points.\footnote{Interestingly, this quantity is roughly equal to the estimated TCI error returned by the algorithm. This emphasizes the reliability of this TCI implementation for approximating such functions.}
	
	\subsubsection{First tests with fixed discretization step in dimensions $d=2$ and $d=3$}
	
\begin{figure}[h]
		\begin{center}
			\begin{tabular}{c}
				$d=2$\\
				\includegraphics[width=0.45\textwidth]{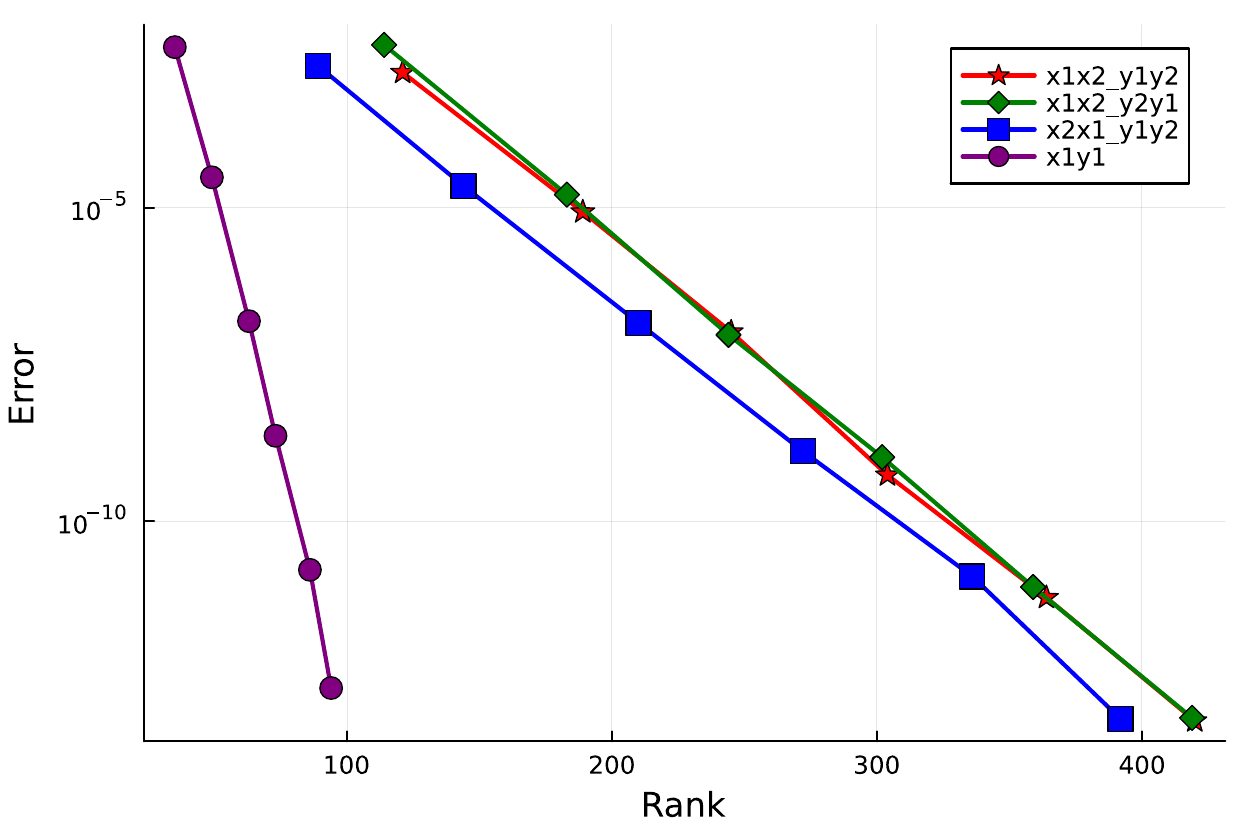}
			\end{tabular}
			\begin{tabular}{c}
				$d=3$\\
				\includegraphics[width=0.45\textwidth]{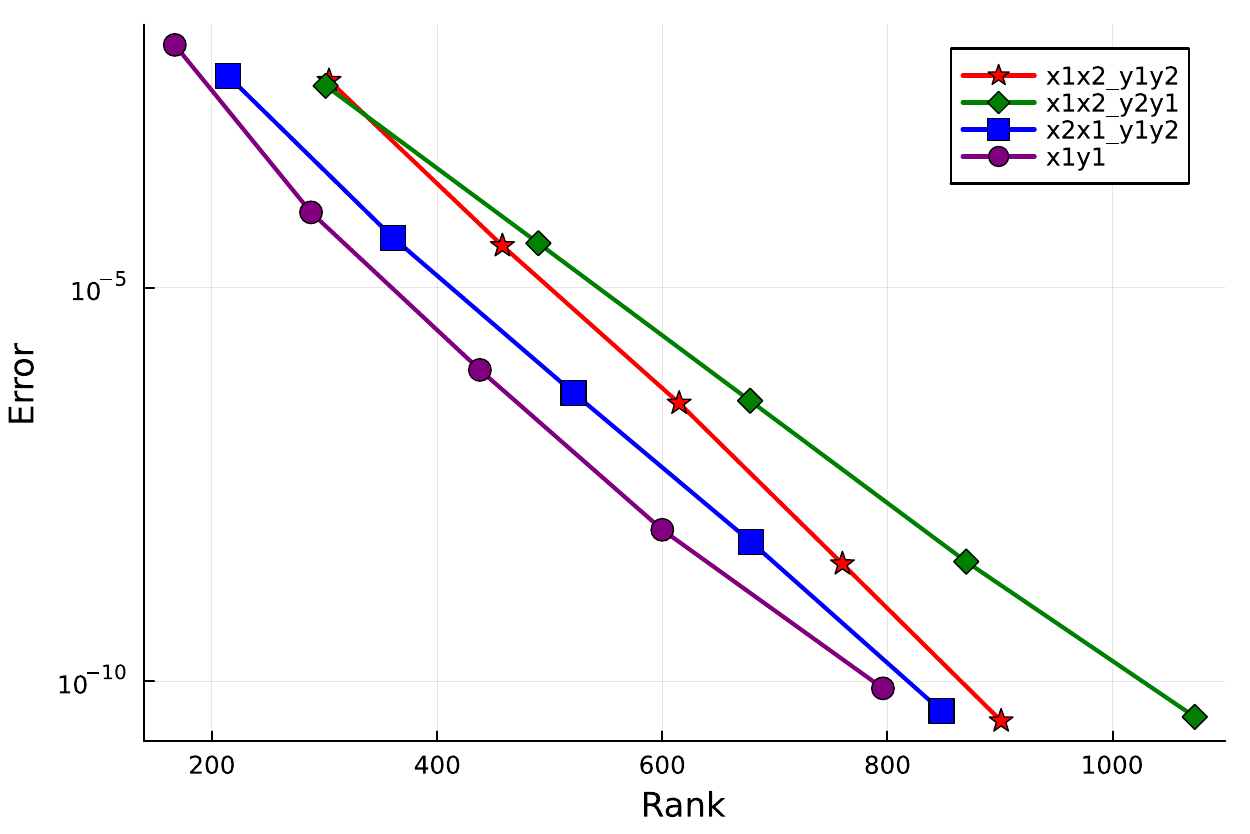}
			\end{tabular}
		\end{center}
		\caption{Empirical relative $\ell^2$ accuracy on the MPS $\bhP$ as a function of its maximal TT-rank, for $L=25$.
			On the left, $d=2$, and on the right $d=3$.
			}
			\label{Fig:d2}
	\end{figure}

	We set a quite high number of dyadic scales $L=25$, which yields a discretization step $h \simeq 3 \cdot10^{-8}$.
	Then, we build MPS approximations for $\bhP$ and we draw their relative accuracies as a function of the TT-rank, in dimensions $d=2$ and $d=3$ in Figure~\ref{Fig:d2}.
	First, in both cases $d=2$ and $d=3$, for all formats and discretizations, we observe that the error approximations scale roughly exponentially with the TT-rank, which is highly favorable.
	We also observe that a relative accuracy of $10^{-10}$ can be reached for all formats, with operators of TT-ranks of at most $350$ in dimension $d=2$ and $1100$ in dimension $d=3$.
	Such a behavior was previously observed for similar functions in \cite{khoromskij2011d}.

	Now, we come to a finer interpretation and consider more precisely the actual TT-ranks.

	In dimension $d=2$, we observe that the format \texttt{x1y1} is the most favorable.
	With this format, our approximation of $\bhP$ enjoys a relative accuracy of $10^{-10}$ and features a TT-rank of $80$.
	All other formats behave roughly the same, and have a TT-rank of approximately $330$ for this discretization step and accuracy.
	As a consequence, \textit{all} the considered formats can be employed in practice.\\
	As a reminder, \cite{bachmayr_stability_2020} provides an MPO of TT-rank at least $1152$, scaling multiplicatively in the maximal TT-rank of $a$, \textit{cf.} \cite[Rk.\ 4]{bachmayr_stability_2020}, with however the property of not depending of $L$.
	Here, by contrast with \cite{bachmayr_stability_2020}, the TT-rank of the final MPO $\bha*+\mu\bhP\odot$ is much smaller: with format \texttt{x1y1}, it is of about $80$ for the homogeneous coefficient $a=1$, and it scales with that of the Fourier transform of the coefficient field $a$ in an affine way, but with a slope of $4$.

	In dimension $d=3$, the cost of approximating $\bhP$ is obviously higher, with a TT-rank being at least of $800$ for the most favorable format.
	Once more, the format \texttt{x1y1} is the best for obtaining a relative accuracy of $10^{-10}$, but it is not that well separated from other formats, the TT-ranks of which range between $800$ and $1050$.
	We underline that, although a TT-rank of order $1000$ is not low, it may be seen as moderately high in the sense that it can still be used in practical computations on a desktop computer.
	In particular, unlike \cite{bachmayr_stability_2020}, our strategy can be used in dimension $d=3$.
	We exemplify this later in Section~\ref{Sec:2D-3D}.

	\subsubsection{Dependence on the discretization step}

	Now, we probe which TT-rank is necessary to obtain a fixed relative accuracy on $\bhP$ in $\ell^2$ norm (which is evaluated \textit{a posteriori} by the Monte-Carlo method), with a variable number $L$ of dyadic scales (the discretization step is $h=2^{-L}$).
	This prescribed relative accuracy is of $10^{-9}$ in dimension $d=2$, and $2\cdot 10^{-7}$ in dimension $d=3$.
	
	\begin{figure}[h]
		\begin{center}
			\begin{tabular}{c}
				$d=2$\\
				\includegraphics[width=0.4\textwidth]{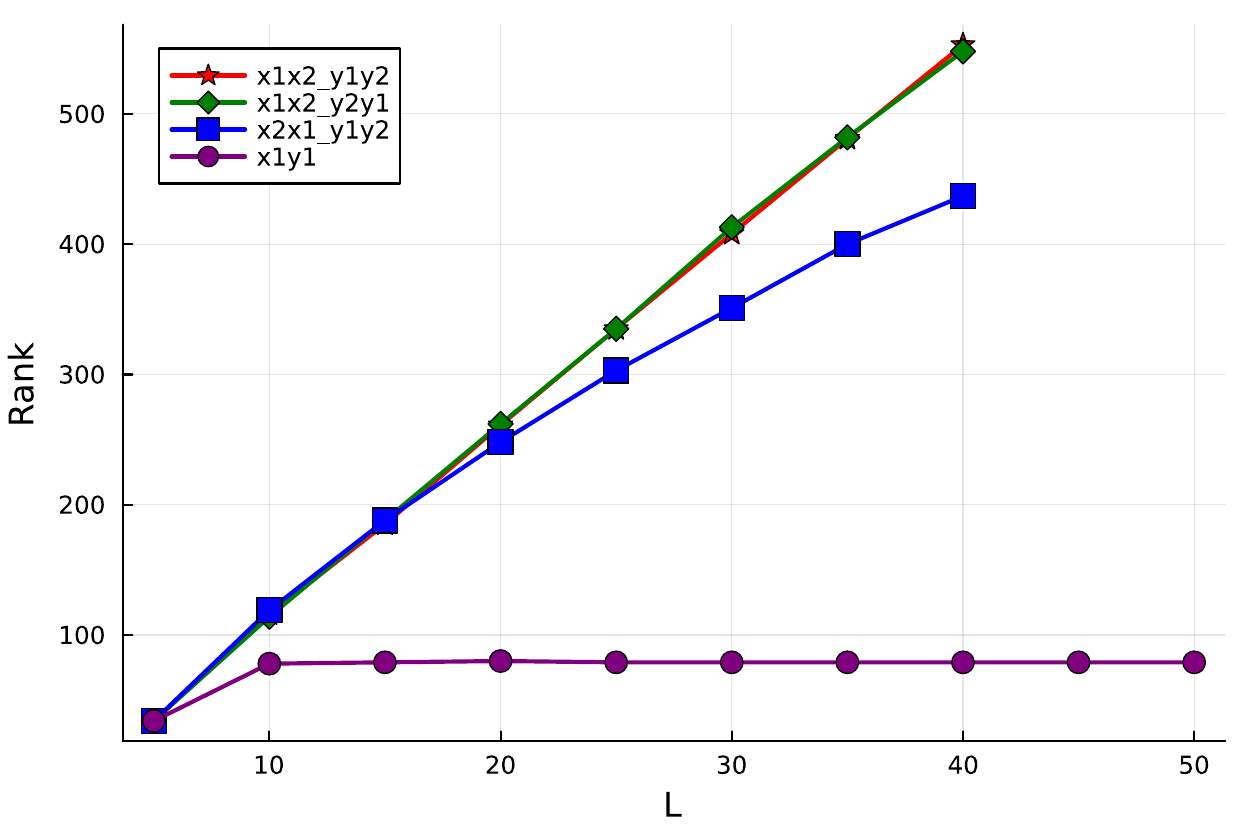}
			\end{tabular}
			\begin{tabular}{c}
				$d=3$\\
				\includegraphics[width=0.4\textwidth]{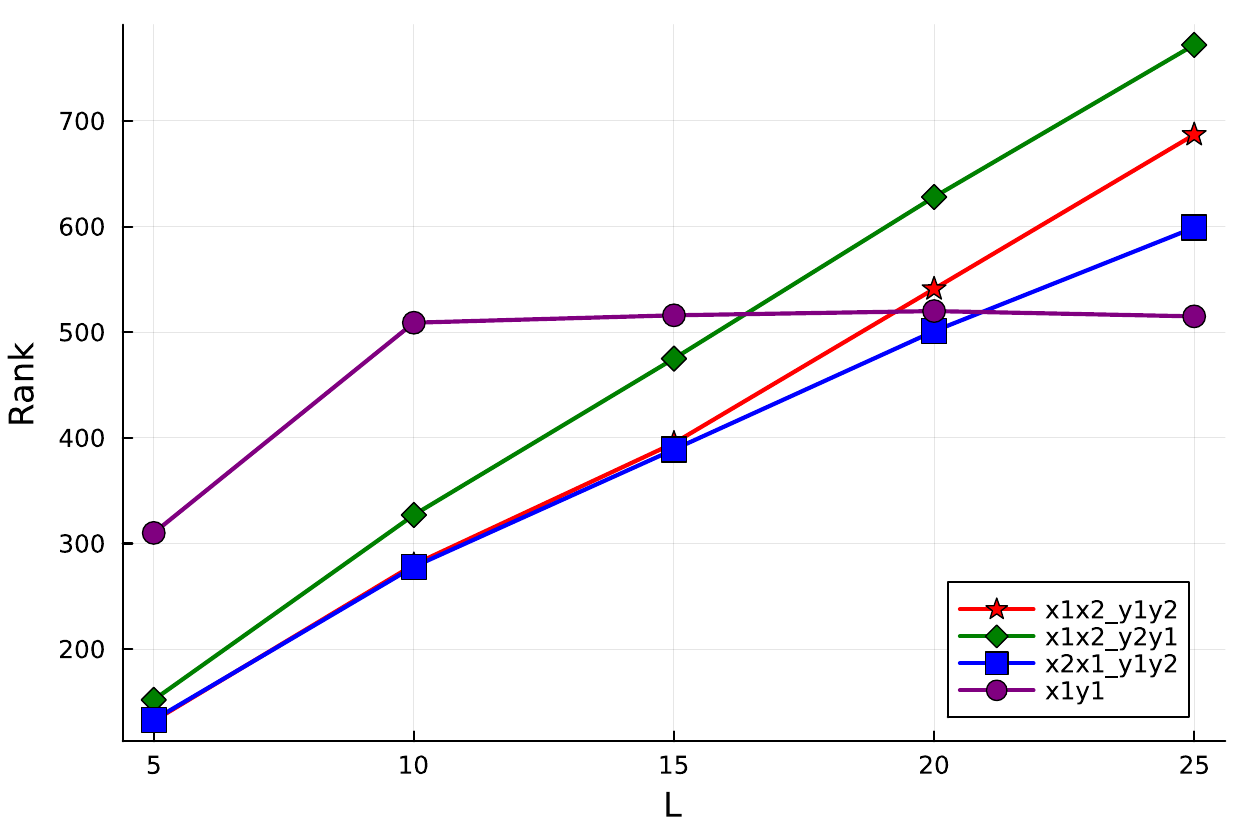}
			\end{tabular}
		\end{center}
		\caption{
			TT-rank of $\bhP$ for a fixed relative accuracy, as a function of the number $L$ of dyadic scales.
			On the left, in dimension $d=2$, for a relative $\ell^2$ accuracy of $10^{-9}$ and, on the right, in dimension $d=3$ for a relative $\ell^2$ accuracy of $2\cdot 10^{-7}$.}
		\label{Fig:rankP}
	\end{figure}

	The result is shown in Figure \ref{Fig:rankP}a) for the dimension $d=2$.
	Very interestingly, we observe that the maximal TT-rank obeys a scaling depending on the format:
		for formats \texttt{x1x2\_y1y2} and \texttt{x1x2\_y2y1}, the TT-rank scales \textit{linearly} with $L$,
		for format \texttt{x2x1\_y1y2}, the TT-rank scales slightly \textit{sublinearly} with $L$, 
		for format \texttt{x1y1}, for $L\geq 10$, the TT-rank is \textit{constant} (equal to $80$) when $L$ increases.
	This emphasizes the superiority of the format \texttt{x1y1} over the other formats, when the mesh size decreases.
	In Figure \ref{Fig:rankP}, we test the format \texttt{x1y1} down to a discretization step $h \simeq 10^{-15}$ (stopping there is arbitrary).
	
	The above qualitative observations remain true in dimension $d=3$, \textit{cf.} Figure \ref{Fig:rankP}b).
	As expected, by comparison with $d=2$, whatever the format is, the TT-ranks are obviously larger for a slightly worse relative accuracy.
	For example, the best format for large $L$ is still \texttt{x1y1}, but a relative accuracy of $2\cdot 10^{-7}$ can be achieved with a TT-rank of $520$, almost constant in $L\geq 10$, to be compared with the previous values of $10^{-9}$ and $80$ for $d=2$.
	As in \cite{bachmayr_stability_2020}, we expect that the TT-rank of $\bhP$ is exponential in dimension $d$.
	Hence, our strategy may only work in moderate dimensions $d=2$ and $d=3$.
	
	Notice that the value of the TT-rank of $\bhP$ for fixed discretization step $2^{-L}$ and fixed relative accuracy should be in principle the same for formats \texttt{x1x2\_y2y1} and \texttt{x2x1\_y1y2} in dimension $d=3$, because these are symmetric of each other, exchanging the variables $x \leftrightarrow z$.
	Thus, we think that the differences between the green and blue curves in Figure \ref{Fig:rankP}b) are non significant.
	By contrast, these format are definitely different in dimension $d=2$, where it seems reasonable to think that gluing the indices on the largest scale is more favorable than gluing them on the smallest scale.
	
\subsection{Comparison of different QTT-based strategies for a manufactured heterogeneous problem}\label{Sec:Comparison}
	
	This subsection showcases a comparison between three QTT-based solvers: 
	\begin{itemize}
		\item QTT-FD for the QTT Finite Difference scheme \cite{oseledets2010approximation, khoromskij2011d, khoromskij_tensor_2018};
		\item QTT-BPX for the Finite Element solver with multilevel Bramble–Pasciak–Xu \cite{bramble1990parallel} preconditioner introduced in \cite{bachmayr_stability_2020};
		\item QTT-HL for the solver based on the Helmholtz-Leray projector introduced here.
	\end{itemize}
	Each of these schemes requires solving a linear system $Ax=b$, where $A$ is an MPO, $b$ is a given MPS, and $x$ is the unknown MPS.
	Here, we appeal to the same 1-site ALS solver for inverting this linear system.
	
	In principle, the three strategies benefit from the QTT compression; but they actually enjoy quite different numerical properties.
	To probe these, we evaluate the three methods on a manufactured problem of the form \eqref{E} in the case $d=2$.
	The comparison is based on several quantitative indicators, including the relative error over the solution and its gradient, 
	the TT-rank of the involved MPOs, and the maximal condition numbers of the local linear systems in the ALS algorithm.
	Indeed, as discussed in \cite{holtz_alternating_2012}, each ALS sweep involves the solution to local linear systems along the TT.
	Since the errors propagate through the sweeping procedure, high condition numbers of these local linear systems may deteriorate the numerical stability and accuracy of the global method.

\subsubsection{Protocol}

Since the QTT-FD and QTT-BPX schemes aim at recovering $u$, they cannot be directly applied to \eqref{E}, for which $u$ is not unique (\textit{cf.} Section \ref{Sec:Pb}).
We restore the uniqueness of $u$ by adding a massive term to \eqref{E}. Hence, we consider the following equation:
\begin{align}\label{EE}
	\gamma u - \nabla \cdot (a \nabla u) = f + \nabla \cdot g,
\end{align}
where $\gamma$ is a positive mass parameter.\footnote{We may have chosen to penalize the integral $\int u$, but adding a massive term proves easier with our numerical tools.}
For QTT-FD and QTT-BPX, we set a moderately small mass parameter $\gamma = 1\cdot 10^{-2}$, whereas we naturally take $\gamma=0$ for QTT-HL.
In all cases, we prescribe both the coefficient field $a$ and the solution $u$, from which we deduce the right-hand side terms $f$ and $g$.
More precisely, we define
\begin{align*}
	&a(x,y):= \big(1 + \nu \sin(4\pi x)\big)\big(1 + \nu \sin(8\pi y)\big),
	\qquad u(x,y):= \frac{1}{6\pi \sqrt{5}} \, \sin^2(6\pi x)\sin^2(12\pi y),
	\\
	&
	f(x,y) = \gamma u(x,y) \qquad \and
	g(x,y) := -a(x,y) \nabla u(x,y) - \tilde{g}(x,y),
	\qquad \pour
	\tilde{g}(x,y) := \lt( \begin{array}{c} 1 + \sin(2\pi y) \\ 2 + \sin(4\pi x)\end{array} \rt),
\end{align*}
where $\nu:=\frac{3}{7}$.
Notice that $\tilde{g}(x,y)$ is of vanishing divergence; we merely use it to trigger the penalization for QTT-HL; otherwise, we may find exactly the solution even when using small penalty parameter.
From the definition of $a$, we derive the following explicit bounds:
$$
\lambda^{-1} := (1-\nu)^2 \leq a(x,y) \leq (1+\nu)^2=:\Lambda.
$$

All considered functions $a$, $u$, and $g$ admit exact low-rank QTT representations.
As a consequence, the present benchmark is particularly fit for QTT-based methods, allowing for fast resolution regardless of the mesh size $h:=2^{-L}$, where $L$ denotes the refinement level.
Nevertheless, the goal here is not to exploit this structure artificially, but to provide a controlled setting for a systematic comparison of the considered solvers.

In the ALS procedure, the initial guess is generated randomly and subsequently normalized.
Its TT-rank is chosen slightly larger than the TT-rank of the analytical solution.
The iterative procedure is stopped when either of the following criteria is satisfied: when the relative change between successive iterates
$$
\delta_{u^k}=\frac{\|u^{k} - u^{k-1}\|_2}{\|u^{k}\|_2},
$$
where $u^k$ denotes the approximation at sweep $k$
falls below ${\rm tol}=1\cdot10^{-6}$, or when the maximum number $S=30$ of ALS sweeps is reached.

The QTT-FD and QTT-HL schemes are evaluated across the QTT formats listed in Table~\ref{Tab:Format2}, whereas QTT-BPX is restricted to its native format \texttt{x1y1}.\footnote{In \cite{bachmayr_stability_2020}, the BPX structure is implemented by contracting virtual indices between cores at the same scale, leading to physical indices of size $2^d$. In dimension $d=2$, this corresponds to merging the two cores corresponding to indices $x_\ell$ and $y_\ell$ into a single core of index $\overline{x_\ell y_\ell}$.
For simplicity, we disregard this subtlety here.
}

\subsubsection{Numerical tests for QTT-FD}

For the QTT-FD scheme, we discretize the operator $-\nabla\cdot(a\nabla u)$ as
${\bm D}^* (\ba {\bm D} \bu)$,
where ${\bm D}$ is the classical first-order forward finite difference operator, and ${\bm D}^*$ its associated transposed operator.
The condition number of the operator ${\bm D}^*\ba {\bm D}$ is bounded by
\begin{align}\label{Condth}
	{\rm Cond}^{\rm th}_{\max} := \frac{8\Lambda h^{-2} + \gamma}{\gamma}.
\end{align}
We explicitly express ${\bm D}$ by means of an MPO of TT-rank $2d+1$ (see \textit{e.g}, \cite{khoromskij_tensor_2018}).

\medskip

We report the results in Table~\ref{tab:QTT-FD1}, where we consider the discretization levels $L \in \{10,20, 30\}$. 

\begin{table}[h]
	\begin{center}
		\scriptsize
		
		{\renewcommand{\arraystretch}{1.5}
			\begin{tabular}{cccccccc}
				\hline
				L & Format & Err$_u$ & Err$_{\nabla u}$ & Sweeps & Meas. Max. Cond. & ${\rm Cond}^{\rm th}_{\max}$ & Maximal relative error for local systems \\
				\hline
				\hline
				\multirow{4}{*}{$10$} & \texttt{x1x2\_y1y2} & $3.83\cdot 10^{-3}$ & $3.39\cdot 10^{-2}$ & 3 & $1.07\cdot 10^{9}$ & \multirow{4}{*}{$1.71\cdot 10^{9}$} & $4.55\cdot 10^{-13}$ \\
				& \texttt{x1x2\_y2y1} & $3.83\cdot 10^{-3}$ & $3.29\cdot 10^{-2}$ & 3 & $1.41\cdot 10^{9}$ &  & $3.41\cdot 10^{-13}$ \\
				& \texttt{x2x1\_y1y2} & $3.83\cdot 10^{-3}$ & $3.35\cdot 10^{-2}$ & 3 & $7.50\cdot 10^{8}$ &  & $2.34\cdot 10^{-13}$ \\
				& \texttt{x1y1} & $3.83\cdot 10^{-3}$ & $3.47\cdot 10^{-2}$ & 8 & $3.83\cdot 10^{8}$ &  & $2.32\cdot 10^{-13}$ \\
				\hline
				\hline
				\multirow{4}{*}{$20$} & \texttt{x1x2\_y1y2} & $5.74\cdot 10^{-5}$ & $1.22\cdot 10^{-4}$ & 30 & $1.77\cdot 10^{15}$ & \multirow{4}{*}{$1.80\cdot 10^{15}$} & $5.77\cdot 10^{-7}$ \\
				& \texttt{x1x2\_y2y1} & $7.14\cdot 10^{-5}$ & $1.11\cdot 10^{-4}$ & 30 & $4.02\cdot 10^{15}$ &  & $3.62\cdot 10^{-7}$ \\
				& \texttt{x2x1\_y1y2} & $1.20\cdot 10^{-5}$ & $6.01\cdot 10^{-5}$ & 30 & $5.90\cdot 10^{13}$ &  & $2.80\cdot 10^{-7}$ \\
				& \texttt{x1y1} & $2.35\cdot 10^{-5}$ & $3.88\cdot 10^{-5}$ & 30 & $3.55\cdot 10^{12}$ &  & $2.45\cdot 10^{-7}$ \\
				\hline
				\hline
				\multirow{4}{*}{$30$} & \texttt{x1x2\_y1y2} & $1.83\cdot 10^{1}$ & $1.01$ & 30 & $3.41\cdot 10^{19}$ & \multirow{4}{*}{$1.88\cdot 10^{21}$} & $6.13\cdot 10^{17}$ \\
				& \texttt{x1x2\_y2y1} & $1.09$ & $1.06$ & 30 & $1.68\cdot 10^{18}$ &  & $1.48\cdot 10^{7}$ \\
				& \texttt{x2x1\_y1y2} & $2.93$ & $1.05$ & 30 & $1.58\cdot 10^{18}$ &  & $2.51\cdot 10^{1}$ \\
				& \texttt{x1y1} & $7.05\cdot 10^{-1}$ & $5.61\cdot 10^{-1}$ & 30 & $2.37\cdot 10^{16}$ &  & $2.19\cdot 10^{-1}$ \\
				\hline
			\end{tabular}
		}
	\end{center}
	\caption{Results for QTT-FD across QTT formats and discretization levels: relative errors and TT-ranks of $u$ and $\nabla u$, number of sweeps for ALS, measured maximum condition numbers of local linear systems and its theoretical bound \eqref{Condth}, and maximal relative error in local linear systems for ALS, \textit{cf.} \eqref{Rel:Err}.
	}
	\label{tab:QTT-FD1}
\end{table}

First, we observe that, for $L=10$, all formats yield comparable errors for both $u$ and $\nabla u$, of order $1$ with respect to the mesh size $h=2^{-L}$ (linear scaling in $h$ is observed down to $L=15$ --not shown here).
For $L=20$, there are some discrepancies among the formats, but still the errors decays satisfactorily.
Then, for $L=30$ and for all QTT formats of Table \ref{Tab:Format2}, the errors are so high that the approximations of both $u$ and $\nabla u$ become irrelevant.

Second, we see that the measured condition numbers of local linear systems in ALS deteriorate rapidly with mesh refinement.
Indeed, as shown in \cite{holtz_alternating_2012} for general MPOs, the bound on the global condition number \eqref{Cond} is also a bound on the condition number of the local systems.
This observation for QTT-FD is well-known in the literature, see, \textit{e.g.} \cite{rakhuba_robust_2021, bachmayr_stability_2020}-- and is the main motivation for more robust strategies.
In Table \ref{tab:QTT-FD1}, the observed values remain below the theoretical upper bound for $L=10$.
For $L \geq 20$, the values are so high that they are beyond double precision\footnote{We think that numerical errors are the main reason for which we measure few condition numbers beyond  ${\rm Cond}^{\rm th}_{\max}$ for some formats for $L= 20$.}; as a consequence, the local systems become harder to invert.

To validate this hypothesis, we report in Table \ref{tab:QTT-FD1} the maximal relative error when solving local linear systems $Ax=b$ in the ALS procedure, that is the quantity
\begin{align}\label{Rel:Err}
	\frac{\|Ax - b\|_2}{\|b\|_2}.
\end{align}
As the condition numbers of these local linear systems explode when $L$ increases, we observe that the maximal relative error \eqref{Rel:Err} increases gradually from order $10^{-13}$ for $L=10$ to orders ranging between $10^{-1}$ and $10^{17}$ for $L=30$, depending on the formats.
Obviously, the accuracy of the whole ALS algorithm is impacted by the high residues in the local linear systems; this explains the instability of the method for extremely small mesh size.
Hence, we interpret the failure of QTT-FD for $L=30$ as a consequence of the poor condition numbers of the local systems.

We underline that we have done these computations with a relative low TT-rank of the initial guess, namely $r=20$.
Should the TT-rank increase, the observed error would be amplified, as the sizes of the local matrices of the ALS procedure will grow.

\subsubsection{Numerical tests for QTT-BPX}

We next turn to QTT-BPX.
The only difference in our implementation with respect to the original construction \cite{bachmayr_stability_2020} is that we adapt it to periodic boundary conditions --this does not change the structure of the MPOs nor their TT-ranks.
Although \cite{bachmayr_stability_2020} does not explicitly treat this case, we believe that their results and proofs remain valid.
To be more specific, we think that the periodic QTT-BPX is also stable with bounded \textit{representation condition numbers}, \textit{cf.} for a precise definition \cite[Sec.\ 5.5]{bachmayr_stability_2020}, resulting in the ability to tackle extremely fine mesh sizes.

\medskip

We show the results in Table \ref{tab:BPX}, where we use the discretization levels $L\in \{10, 20, 30, 40\}$.

\begin{table}[h]
	\begin{center}
	\scriptsize
	
	{\renewcommand{\arraystretch}{1.5}
	\begin{tabular}{c c c c c c c c}
		\hline
		$L$ & Format & Err$_u$ & Err$_{\nabla u}$ & Sweeps & Meas. Max Cond. \\
		\hline
		\hline
		
		$10$& x1xy1 & $4.65\cdot10^{-3}$ & $4.42\cdot10^{-3}$ & 15 & $7.64\cdot10^{2}$ \\
		\hline
		\hline	
				
		$20$&
		x1xy1 & $3.83\cdot10^{-5}$ & $5.89\cdot10^{-4}$ & 28 & $8.57\cdot10^{2}$ \\
		\hline
		\hline
		
		$30$ &
		x1xy1 & $3.81\cdot10^{-5}$ & $5.90\cdot10^{-4}$ & 30 & $8.57\cdot10^{2}$ \\
		\hline
		\hline
		
		$40$
		&x1xy1 & $3.81\cdot10^{-5}$ & $5.90\cdot10^{-4}$ & 25 & $8.57\cdot10^{2}$ \\		
		\hline
	\end{tabular}
	}
	\end{center}
	\caption{Results for QTT-BPX across discretization levels: numerical relative errors on $u$ and~$\nabla u$, number of sweeps for the ALS, and measured maximum condition numbers of local linear systems.}
	\label{tab:BPX}
\end{table}

In accordance with \cite{bachmayr_stability_2020}, we observe that the method is robust regardless of the discretization step down to $h=2^{-40} \simeq 10^{-12}$.
This is in relationship with the fact that the maximum condition number of the local linear systems stays roughly constant across all discretization levels, staying below $10^{3}$.
This property is the consequence of a special construction guaranteeing the control of the representation condition number of the involved MPOs --only preconditioning is not sufficient, see \cite{bachmayr_stability_2020} for more details.
When $L$ increases, the approximation errors on $u$ and $\nabla u$ decrease and finally reach a plateau not far from the prescribed tolerance.
The latter is matched by ALS for all experiments but one within $30$ sweeps.

\subsubsection{Numerical tests for QTT-HL}\label{Sec:QTTHL}

\begin{figure}[h]
	{
		\begin{center}
			\includegraphics[width=0.4\textwidth]{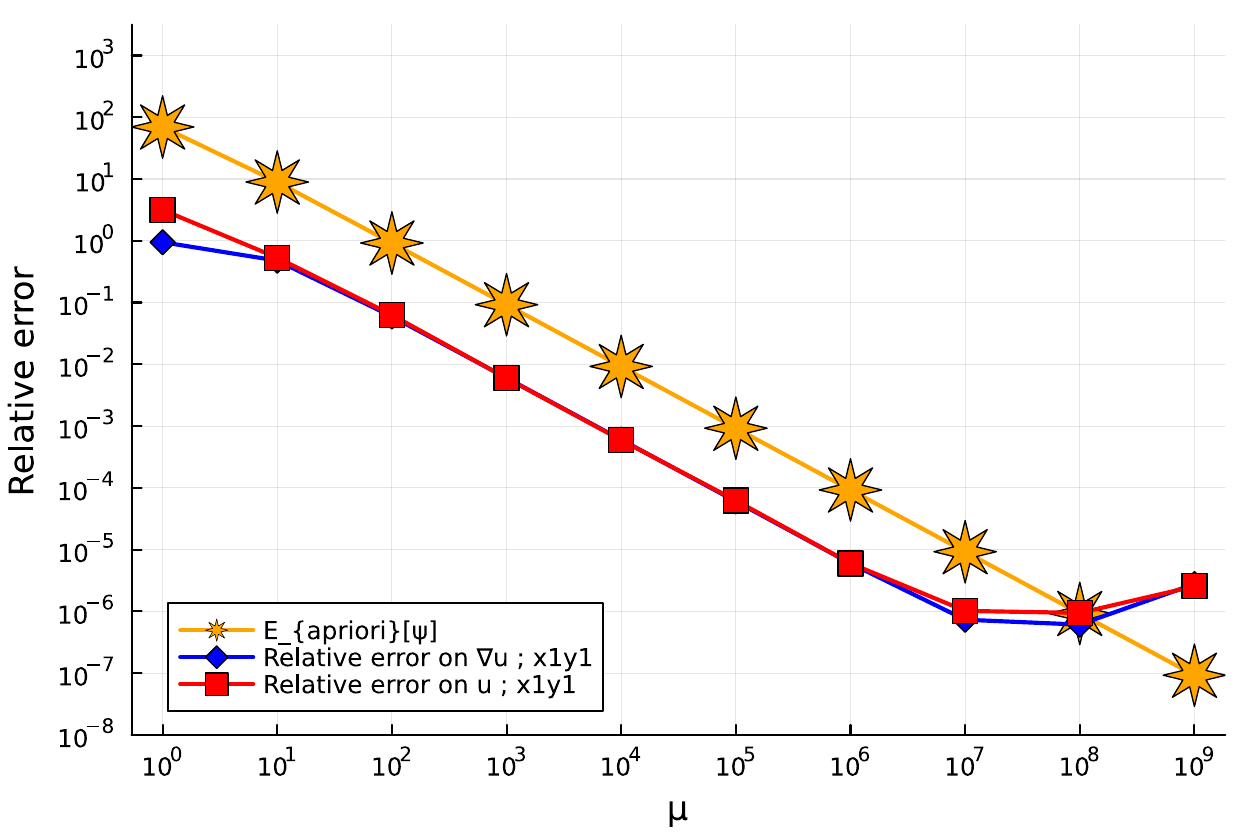}
			\includegraphics[width=0.4\textwidth]{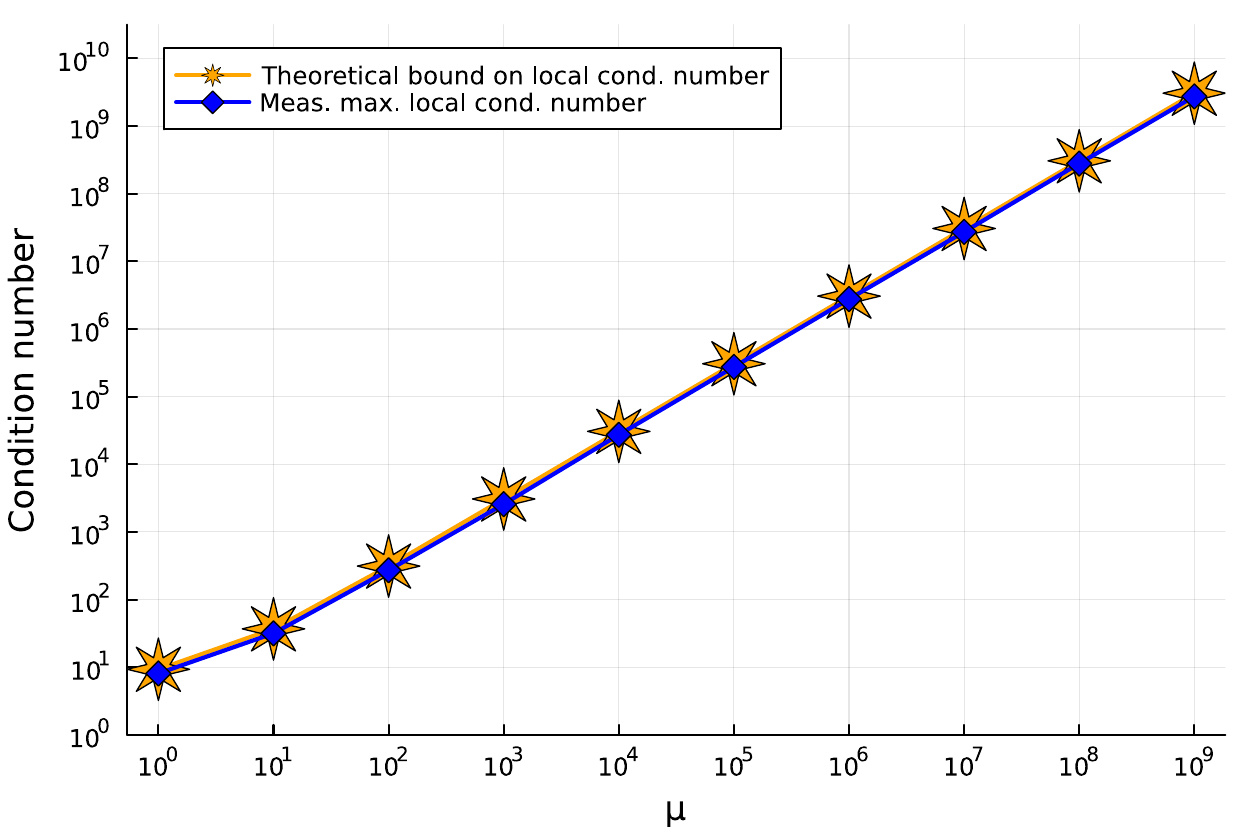}
		\end{center}
	}
	\caption{
		On the left, relative errors on $u$ and $\nabla u$ as functions of $\mu$ for the QTT-HL strategy, compared with the theoretical bound on $\nabla u$; on the right, measured condition number of local systems and associated theoretical bound as a function of $\mu$.
		We use the following parameters: QTT format is \texttt{x1y1}, $L=20$, ${\rm tol}=1\cdot 10^{-6}$, $S=30$.
	}
	\label{Fig:com5}
\end{figure}

First, we study the impact of the penalty parameter $\mu$.
In Figure \ref{Fig:com5}, we pick the QTT format \texttt{x1y1} and set $L=20$, and we make $\mu$ vary between $1$ and $10^9$.
(Similar experiments for other QTT formats of Table \ref{Tab:Format2} yield similar observations --but are not shown here.)

As expected, when $\mu$ increases, we observe that the relative error on $\nabla u$ between the numerical approximation of the method and the analytical value decreases linearly in $\mu^{-1}$ until $\mu=10^8$.
Since the ALS iterations stop when $\delta_{u^k} \leq 10^{-6}$, the relative error naturally saturates before $10^{-6}$, which explains the deviation from the linear scaling at $\mu=10^{8}$.
Also, we see that the theoretical \textit{a priori} error estimate derived from Lemma \ref{Prop:Mu} is satisfied until $\mu=10^{8}$.

Next, we observe that the measured maximal condition number of local systems increases linearly in $\mu$, and virtually matches by below the estimate of Proposition \ref{Prop:Cond}.
(A precise inspection of the values yields that the measured maximal local condition number is $9\%$ to $14\%$ smaller than the theoretical bound.)
A higher condition number makes local linear systems harder to invert.
As a trade-off between accuracy and condition number, we henceforth set $\mu=10^4$.

\medskip

Now, we make the QTT formats and $L$ vary, and report the results in Table~\ref{tab:QTT-HL1}.

\begin{table}[h]
	\begin{center}
	\scriptsize
	{\renewcommand{\arraystretch}{1.5}
		
\begin{tabular}{cccccccccc}
	\hline
	L & Format & Err$_u$ & Err$_{\nabla u}$ & $E_{\min}[\psi]$ & $E_{\max}[\psi]$ & $E_{\rm apriori}$ & Sweeps & Meas. Max. Cond. & ${\rm Cond}^{\rm th}_{\max}$ \\
	\hline
	\hline
	\multirow{4}{*}{$10$} & \texttt{x1x2\_y1y2} & $6.02\cdot 10^{-4}$ & $7.93\cdot 10^{-4}$ & $4.62\cdot 10^{-4}$ & $3.35\cdot 10^{-3}$ & \multirow{4}{*}{$9.23\cdot 10^{-3}$} & 4 & $2.82\cdot 10^{4}$ & \multirow{4}{*}{$3.06\cdot 10^{4}$} \\
	& \texttt{x1x2\_y2y1} & $6.02\cdot 10^{-4}$ & $7.86\cdot 10^{-4}$ & $4.62\cdot 10^{-4}$ & $3.35\cdot 10^{-3}$ &  & 4 & $2.81\cdot 10^{4}$ &  \\
	& \texttt{x2x1\_y1y2} & $6.02\cdot 10^{-4}$ & $7.54\cdot 10^{-4}$ & $4.62\cdot 10^{-4}$ & $3.35\cdot 10^{-3}$ &  & 3 & $2.78\cdot 10^{4}$ &  \\
	& \texttt{x1y1} & $6.02\cdot 10^{-4}$ & $7.86\cdot 10^{-4}$ & $4.62\cdot 10^{-4}$ & $3.35\cdot 10^{-3}$ &  & 4 & $2.55\cdot 10^{4}$ &  \\
	\hline
	\hline
	\multirow{4}{*}{$20$} & \texttt{x1x2\_y1y2} & $6.02\cdot 10^{-4}$ & $5.99\cdot 10^{-4}$ & $4.62\cdot 10^{-4}$ & $3.35\cdot 10^{-3}$ & \multirow{4}{*}{$9.23\cdot 10^{-3}$} & 4 & $2.83\cdot 10^{4}$ & \multirow{4}{*}{$3.06\cdot 10^{4}$} \\
	& \texttt{x1x2\_y2y1} & $6.02\cdot 10^{-4}$ & $5.90\cdot 10^{-4}$ & $4.62\cdot 10^{-4}$ & $3.35\cdot 10^{-3}$ &  & 5 & $2.89\cdot 10^{4}$ &  \\
	& \texttt{x2x1\_y1y2} & $6.02\cdot 10^{-4}$ & $5.75\cdot 10^{-4}$ & $4.62\cdot 10^{-4}$ & $3.35\cdot 10^{-3}$ &  & 3 & $2.80\cdot 10^{4}$ &  \\
	& \texttt{x1y1} & $6.02\cdot 10^{-4}$ & $5.82\cdot 10^{-4}$ & $4.62\cdot 10^{-4}$ & $3.35\cdot 10^{-3}$ &  & 5 & $2.74\cdot 10^{4}$ &  \\
	\hline
	\hline
	\multirow{4}{*}{$30$} & \texttt{x1x2\_y1y2} & $6.02\cdot 10^{-4}$ & $5.74\cdot 10^{-4}$ & $4.62\cdot 10^{-4}$ & $3.35\cdot 10^{-3}$ & \multirow{4}{*}{$9.23\cdot 10^{-3}$} & 4 & $2.83\cdot 10^{4}$ & \multirow{4}{*}{$3.06\cdot 10^{4}$} \\
	& \texttt{x1x2\_y2y1} & $6.02\cdot 10^{-4}$ & $6.04\cdot 10^{-4}$ & $4.62\cdot 10^{-4}$ & $3.35\cdot 10^{-3}$ &  & 4 & $2.84\cdot 10^{4}$ &  \\
	& \texttt{x2x1\_y1y2} & $6.02\cdot 10^{-4}$ & $6.23\cdot 10^{-4}$ & $4.62\cdot 10^{-4}$ & $3.35\cdot 10^{-3}$ &  & 3 & $2.85\cdot 10^{4}$ &  \\
	& \texttt{x1y1} & $6.02\cdot 10^{-4}$ & $6.29\cdot 10^{-4}$ & $4.62\cdot 10^{-4}$ & $3.35\cdot 10^{-3}$ &  & 5 & $2.77\cdot 10^{4}$ &  \\
	\hline
	\hline
	\multirow{4}{*}{$40$} & \texttt{x1x2\_y1y2} & $6.02\cdot 10^{-4}$ & $5.76\cdot 10^{-4}$ & $4.62\cdot 10^{-4}$ & $3.35\cdot 10^{-3}$ & \multirow{4}{*}{$9.23\cdot 10^{-3}$} & 4 & $2.84\cdot 10^{4}$ & \multirow{4}{*}{$3.06\cdot 10^{4}$} \\
	& \texttt{x1x2\_y2y1} & $6.02\cdot 10^{-4}$ & $6.05\cdot 10^{-4}$ & $4.62\cdot 10^{-4}$ & $3.35\cdot 10^{-3}$ &  & 4 & $2.86\cdot 10^{4}$ &  \\
	& \texttt{x2x1\_y1y2} & $6.02\cdot 10^{-4}$ & $5.77\cdot 10^{-4}$ & $4.62\cdot 10^{-4}$ & $3.35\cdot 10^{-3}$ &  & 3 & $2.77\cdot 10^{4}$ &  \\
	& \texttt{x1y1} & $6.02\cdot 10^{-4}$ & $6.03\cdot 10^{-4}$ & $4.62\cdot 10^{-4}$ & $3.35\cdot 10^{-3}$ &  & 4 & $2.68\cdot 10^{4}$ &  \\
	\hline
\end{tabular}

	}
	\end{center}
	\caption{Results for QTT-HL across QTT formats and discretization levels: relative errors on $u$ and $\nabla u$, \textit{a posteriori} estimators $E_{\min}[\psi]$ and $E_{\max}[\psi]$ given by \eqref{Error-total} and \eqref{Error-total2}, \textit{a priori} estimator $E_{\rm apriori}$ given by \eqref{num:0005}, number of sweeps for ALS, maximal measured condition number on local linear systems of ALS, and its theoretical bound given by \eqref{Cond}.
	(Results concerning Err$_u$, $E_{\min}[\psi]$, $E_{\max}[\psi]$ are not constant across the tests, but only differ when considering more digits --which are irrelevant for our discussion.)}
	\label{tab:QTT-HL1}
\end{table}

Across all discretization levels, and regardless of the QTT format, QTT-HL exhibits remarkably stable results.
The relative errors on $u$ and $\nabla u$ range between $5\cdot 10^{-4}$ and $8\cdot 10^{-4}$.
In particular, the relative error on $\nabla u$ follows the expected behavior: it is indeed between the two \textit{a posteriori} error estimates, and below the \textit{a priori} error estimate.
The relative error fails to decay with respect to the discretization level $L$, but saturates: this is due to our use of a constant finite penalty parameter $\mu$, which is reflected in the constant \textit{a priori} error.
The number of sweeps before convergence of the ALS procedure is between $3$ and $5$, which is low.
This can be explained by the fact that the condition number of the local linear systems is well-controlled, slightly below the theoretical bound ${\rm Cond}^{\rm th}_{\max}$ of \eqref{Cond}, which is of $3.06 \cdot 10^4$.
All these results are apparently unconditionally stable with respect to the mesh refinement.\footnote{We also tested $L=50$ and obtain similar results.}

Our numerical results concerning the error approximations on $\nabla u$ are not only consequences of Lemma \ref{Prop:Mu} and Proposition \ref{Prop:Aposteriori}, but they also hint that our strategy is robust against perturbations.
We emphasize the nuance here: as is usual, our numerical analysis is blind to the computational errors that are due to machine precision and approximations errors when computing, \textit{e.g.} the QR decomposition of MPS cores (this is required by the ALS procedure).
As explained in \cite{bachmayr_stability_2020}, as a consequence of the TT structure, this may lead to amplified errors, even when the MPOs under consideration approximate well-conditioned operators.
Based on these numerical evidences, it seems that our construction does not suffer from this possible drawback.

\subsubsection{Performance evaluation}\label{Sec:Perf}

Computation times are not very meaningful with our code; hence, in order to evaluate the performance of each strategy, we prefer to report in Table \ref{Tab:ComparisonEnd} the TT-rank of the associated MPO $A$ of the linear system $Ax=b$.
This TT-rank may depend on the QTT format as well as the mesh size.
Notice however that this indicator itself is not a perfect mirror of the computational complexity since further MPO decompositions may be used to accelerate the resolution of the linear system within the ALS scheme, such as a factorization as a product of MPOs.

\begin{table}[h]
	\begin{center}
		\scriptsize
		{\renewcommand{\arraystretch}{1.5}
			\begin{tabular}{c c c c c}
				\hline
				Solver & QTT format & MPO rank  of $A$, $L=15$ & MPO rank  of $A$, $L=20$
				\\
				\hline
				\hline
				\multirow{4}{*}{QTT-FD}
				&\texttt{x1x2\_y1y2} & $79$ & $79$
				\\
				&\texttt{x1x2\_y2y1} & $79$ & $79$
				\\
				
				&\texttt{x2x1\_y1y2} & $79$ & $79$
				\\
				
				&\texttt{x1y1} & $163$ & $163$
				\\
				\hline
				\hline
				QTT-BPX & \texttt{x1y1} &  $11392$ & $11392$
				\\
				\hline
				\hline
				\multirow{4}{*}{QTT-HL}
				&\texttt{x1x2\_y1y2} & $222$ & $310$
				\\
				
				&\texttt{x1x2\_y2y1} & $222$ & $309$
				\\
				
				&\texttt{x2x1\_y1y2} & $218$ & $295$
				\\
				
				&\texttt{x1y1} &  $90$ & $90$
				\\
				\hline
			\end{tabular}
		}
	\end{center}
	\caption{TT-rank of the MPO $A$ of the linear system $Ax=b$ involved in the QTT-FD, QTT-BPX and QTT-HL strategies for $L=15$ and $L=20$. (The error on $\Gamma$ is of order $10^{-12}$ for QTT-HL.)}
	\label{Tab:ComparisonEnd}
\end{table}

We observe on Table \ref{Tab:ComparisonEnd} that QTT-FD with all QTT formats of Table \ref{Tab:Format2}, QTT-BPX with  \texttt{x1y1}, and QTT-HL with \texttt{x1y1} all feature MPOs of TT-rank independent of the discretization fineness.
By comparison, the MPO in QTT-HL is of increasing TT-rank when $L$ increases for other formats.
This is in accordance with the theory and the observations of Section \ref{Sec:Green}.

The lowest MPO ranks are observed for QTT-FD with formats \texttt{x1x2\_y1y2}, \texttt{x1x2\_y2y1}, and \texttt{x2x1\_y1y2}, and for QTT-HL with format \texttt{x1y1} ($79$ and $90$, respectively).
QTT-FD with format \texttt{x1y1}, as well as QTT-HL with other formats shows moderate TT-ranks ($163$ for any $L$, and about $300$ for $L=20$, respectively).
By contrast, QTT-BPX involves an MPO of TT-rank more than one order of magnitude above that of the others methods (about $11000$), which makes it barely tractable on a single computation node.\footnote{Our implementation actually uses an MPO factorization for QTT-BPX, \textit{cf.} \cite{bachmayr_stability_2020}; otherwise the method saturates the computer memory. Yet, even with this optimization, the method proves far slower than the others.}

\subsubsection{General comparison of the three QTT-based schemes}

We now summarize the global behavior of the three considered solvers: QTT-FD, QTT-BPX, and QTT-HL. Although all methods are built within the same tensor-train framework, they exhibit quite different numerical properties in terms of stability, compression efficiency, and scalability.

\medskip

The QTT-FD scheme represents the most direct tensorization of a classical finite-difference discretization. Its main limitation lies in the strong deterioration of the conditioning of the local ALS systems as the mesh is refined, in $h^{-2}$. This leads to rapidly growing condition numbers and, consequently, to instability for fine mesh size $h$.
QTT formats do not alter the instability, which is intrinsic to this method.
Nevertheless, reaching $L=20$ is already slightly beyond the capabilities of classical solvers on a supercomputer (due to memory restrictions) in dimension $d=2$, and well beyond in dimension $d=3$.
QTT-FD remains thus appealing when the mesh size $h$ is not too small; for example, very recently, a multigrid-inspired strategy \cite{li_tailoring_2026} was proposed for replacing ALS-like algorithms in order to further reduce the computational cost for this strategy.

\medskip

The QTT-BPX method significantly improves the conditioning behavior by means of a multilevel BPX-type preconditioner. This results in uniformly bounded condition numbers across discretization levels and a robust convergence.
Furthermore, the numerical errors on $u$ and $\nabla u$ decay linearly in the discretization step $h$.
However, this remarkable (and mathematically proven) stability comes at the cost of substantially increased TT-ranks for the involved MPO, which is exponential in the dimension $d$ and multiplicative with respect to the rank of the coefficient field.
This rank growth severely impacts both memory consumption and computational cost, and becomes a major limitation.
Indeed, without resorting to HPC strategies, dimension $d=3$ cannot be tackled, neither complex heterogeneous coefficient fields in dimension $d=2$.

\medskip

Last, the proposed QTT-HL scheme combined with the QTT format \texttt{x1y1} features satisfactory quantitative properties in terms of rank, stability and conditioning.
In terms of MPO ranks, it is similar to QTT-FD and far better than QTT-BPX. 
By comparison with QTT-BPX, it allows for considering more complex $2$-dimensional coefficient fields, and can tackle $3$-dimension problems (see Section \ref{Sec:homog} below).
In terms of conditioning, like QTT-BPX, it features a global condition number that is bounded uniformly with respect to the mesh size $h$.
Hence, as QTT-BPX, it allows for extremely fine $h$.
These properties are slightly deteriorated for other formats of Table \ref{Tab:Format2}; for these, the MPO rank increases linearly with $L$ (\textit{cf.} Sections \ref{Sec:Green} and \ref{Sec:Perf}).
However, even with these degraded QTT formats, it compares advantageously with QTT-FD and QTT-BPX.

The price to pay is a penalty parameter $\mu \gg 1$, which offers a balance between solution accuracy and condition number.
Although our penalization strategy makes the approximation errors saturate in $\mu^{-1}$, preventing the linear scaling thereof in the mesh size $h$ for $h \ll \mu^{-1}$, this drawback is virtually irrelevant for our multiscale problems.
Indeed, in our cases, we do not seek for an extremely small mesh size for having an extremely precise solution.
Rather, we want to obtain the full field solution to \eqref{E} with a prescribed error approximation (typically below $1\%$ on both the solution and its gradient), which requires to set $h\ll \epsilon$, where $\epsilon$ is the characteristic size of the heterogeneities.
	
	\subsection{Validation of the proposed strategy in a multiscale setting}\label{Sec:homog}
	
	In this section, we study a multiscale problem.
	There, we show that our strategy can tackle a brute-force resolution of the problem, and assess the quality of the solution by comparing it to the Ansatz produced by a homogenization strategy.
	We select the format \texttt{x1y1}, which proves superior according to the insights of the previous sections.
	Nevertheless, other formats can be considered, at the price of larger TT-ranks or poorer error approximation: we discuss this at Section \ref{Sec:OtherFormats} on a single example.
	
	\subsubsection{Problem setting}
	
	We assume that $a=\aeps$ is decomposed as
	\begin{align}\label{Ansatz}
		\aeps(x) = b(x) c(x/\epsilon),
	\end{align}
	where $\epsilon \ll 1$, $b$ is a scalar field and $c$ is a $\Q_1$-periodic scalar field.
	Then, according to the homogenization theory \cite{allaire_shape_2002, tartar_general_2009}, the solution $\ueps$ to
	\begin{align}\label{Eeps}
		- \nabla \cdot \lt(\aeps \nabla \ueps\rt) = \nabla \cdot g \qquad \text{in} \quad \Q_{1,\per},
	\end{align}
	can be approximated by the two-scale expansion
	\begin{align}\label{2-scale}
		\ueps \simeq u^{\epsilon, 1} := \ubar + \epsilon \sum_i \partial_i \bar{u} \phi_i(\cdot/\epsilon),
	\end{align}
	with $\phi_i$ the first-order corrector and $\bar{u}$ the homogenized solution, both defined after.
	Its gradient is accordingly approximated as
	\begin{align}\label{Uncoupling}
		\nabla \ueps \simeq \nabla  u^{\epsilon, 1} := \sum_i \partial_i \bar{u} (e_i + \nabla \phi_i(\cdot/\epsilon)),
	\end{align}
	where $e_i$  denotes the $i$-th canonical  basis vector.
	By a straightforward generalization of the classical periodic homogenization results \cite{jikov_homogenization_1994}, the convergence rates scale as
	\begin{align}\label{CvRate}
		\lt\|\ueps - \ubar\rt\|_{\LL^2(\Q_1)} \lesssim \epsilon \qquad \et \quad
		\big\|\nabla\ueps - \sum_i \partial_i \bar{u} (e_i + \nabla \phi_i(\cdot/\epsilon))\big\|_{\LL^2(\Q_1)} \lesssim \epsilon.
	\end{align}

	In \eqref{2-scale}, the corrector $\phi_i$ is of zero mean and solves the cell problem
	\begin{align}\label{E:corr}
		-\nabla \cdot (c (e_i + \nabla \phi_i)) = 0 \quad \dans \Q_{1,\per},
	\end{align}
	where $e_i$ is a vector of the canonical basis of $\R^d$, and $\ubar$ solves the homogenized equation
	\begin{align}\label{E-hom}
		-\nabla\cdot(\abar \nabla \ubar) = \nabla\cdot g \dans \Q_{1,\per}.
	\end{align}
	The homogenized coefficient field $\abar(x)$ is defined as
	\begin{align}\label{Def:abar}
		\abar(x) \cdot e_i = b(x) \int c(e_i+\nabla\phi_i).
	\end{align}
	Unlike in the classical periodic homogenization setting, here, $\abar(x)$ does depend on $x$, but only through the slowly varying coefficient $b$.
	Yet, here, the specific algebraic structure \eqref{Ansatz} that we choose allows for explicit expression of $\abar(x)$ by solving correctors problems depending on a single microstructure (given by $c$).
	The coefficient $c$ is generated by a simple random process that we periodize: somehow, from a homogenization perspective, $c$ encodes for a Representative Volume Element.
	However, unlike in the usual strategy of scale separation \cite{jikov_homogenization_1994}, we choose here to \textit{pave} the actual material with the RVE and solve the full-field problem by our QTT solver.
	We underline that \eqref{2-scale} and \eqref{Uncoupling} are used only for validating the solution.

	\subsubsection{Specific choices of microstructure}\label{Sec:microEmpirical}
	
	We make the specific choices:
	\begin{align*}
		b(x) = \lt\{
		\begin{aligned}
			&1 + 0.5 \cos(2 \pi x_1)\cos(2 \pi x_2) &&\si d=2, \\
			&1 + 0.5 \cos(2 \pi x_1) && \si d=3,
		\end{aligned}
		\rt.
	 	\qquad \et \quad
		g_i(x) = \cos((4+2i)\pi x_i) \pour i \in \{1, \dots, d\},
	\end{align*}
	where, here, $x_i$ denotes the $i$th spatial coordinates of $x \in [0, 1)^d$, and $g_i$ the $i$th component of $g\in \LLb^2$.
	These functions $b$ and $g$ are quite simple: they can be analytically represented by MPS of TT-ranks lower than $10$ (for any QTT format and any dimension), making virtually no error on these quantities.\footnote{Analytical MPS expressions exist for trigonometric functions \cite[Sec.\ 4.1.1]{khoromskij_tensor_2018}.}
	The complexity lies in $c$, where we choose
	\begin{align}\label{Def:c}
		c(x) = 1 + \sum_{l \in \Z^d} c_0(x-l)
		\pour 
		c_0(x)
		 =  \nu \sum_{n=1}^N G(x-X_i),
	\end{align}
	for a given $\nu >0$ and where $G$ is a Gaussian function
	\begin{align}\label{Def:g}
		G(x) = \exp(-x^2 / 2 r^2),
	\end{align}
	for a given $r>0$.
	Here $X_i$ is generated by a random process, the so-called RSA \cite{torquato_random_2002}.
	This generates $N$ non-intersecting spheres centered at $X_i$ with radius $2r$ in the periodic cube $\Q_{1,\per}$.
	In \eqref{Def:c}, we clearly decoupled the unit cell microstructure $c_0$ and the periodization made by the sum over $l$.
	
	We remark that $a$ satisfies Assumption \ref{A-i} with $\lambda=2$ and $\Lambda \leq 1.5(1+1.5\nu)$.
	
\begin{figure}[h]
		\begin{center}
			\includegraphics[width=0.45\textwidth]{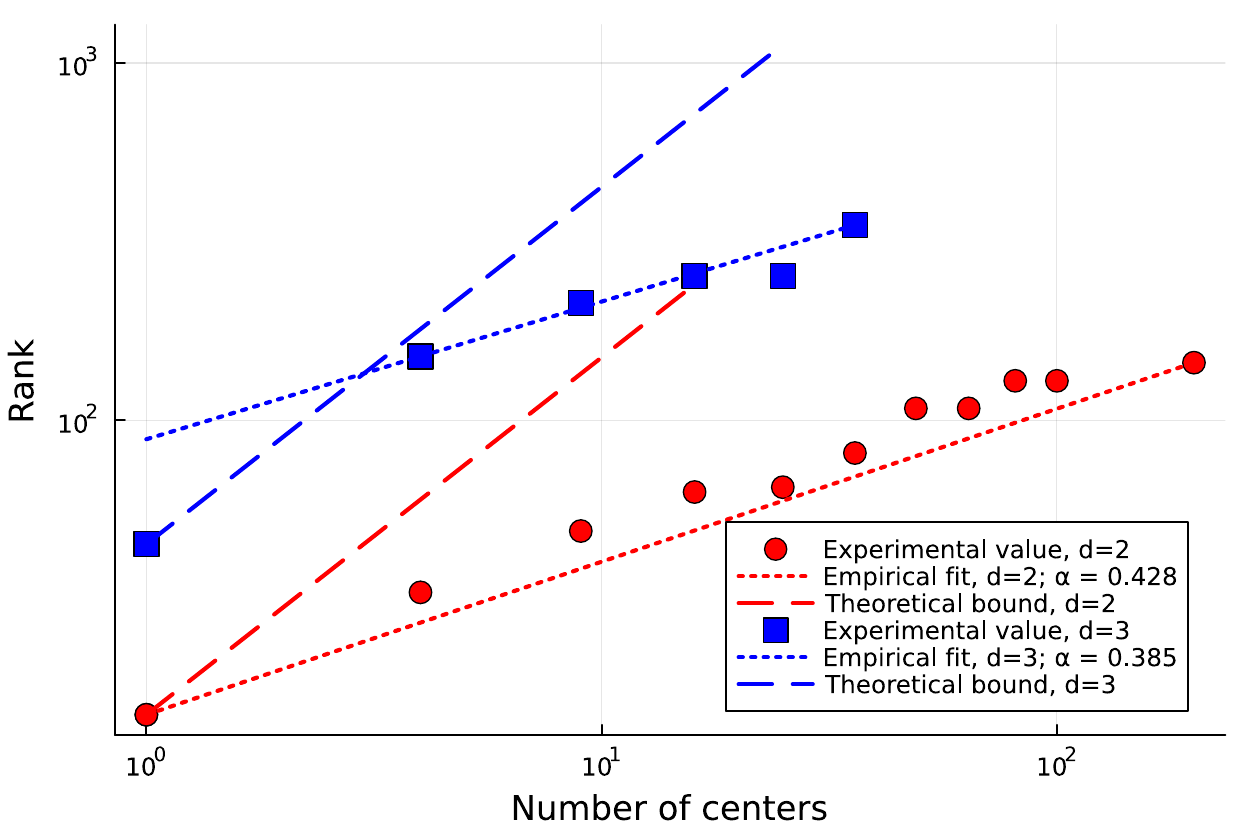}
			\includegraphics[width=0.45\textwidth]{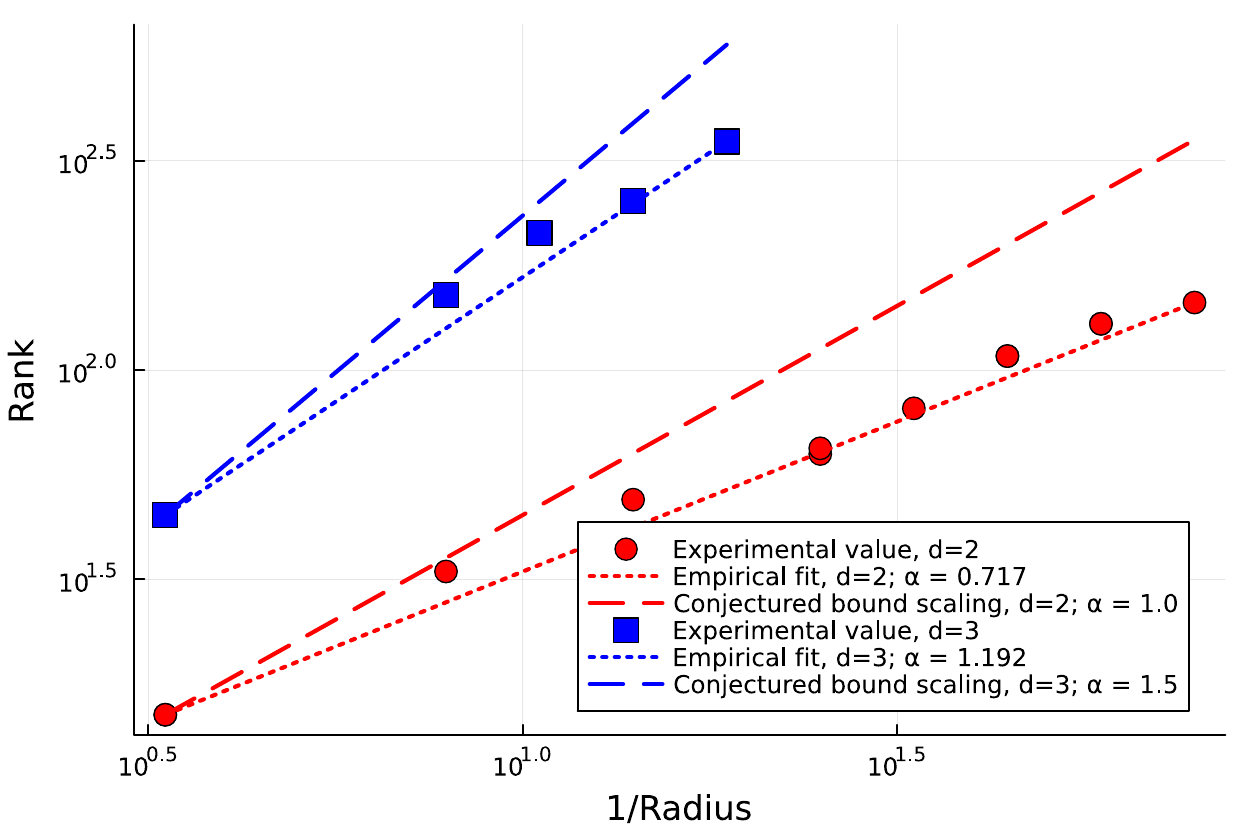}
		\end{center}
		\caption{The maximal TT-rank of $\bm{c}$ as a function of $N$ on the left and $r$ on the right, with Monte-Carlo error of order $10^{-6}$, in dimensions $d=2$ (red) and $d=3$ (blue).
		Here, $\nu=5$, with a small-scale discretization level $L_0=15$, and format \texttt{x1y1}.
		On the left, the dashed lines feature the rigorous theoretical bound linear in $N$; on the right, the dashed lines feature the conjectured scaling from the bound in \cite{lindsey_multiscale_2024}.
		}\label{Fig:RankC}
	\end{figure}

	We investigate on the TT-rank of an MPS ${\bm c}$ representing $c$.
	We set $\nu=5$ resulting in a moderate contrast $\lambda\Lambda$ of about $20$ and choose adaptively $r$ as a function of $N$ so that the spheres of center $X_i$ and radius $2r$ generated by RSA are close to a maximal configuration (\textit{i.e.} one cannot add a new sphere of such radius without intersecting already placed spheres).
	We appeal to TCI to build an MPS $\bm{G}$ representing $G$, that we convolve periodically with $\sum_{i=1}^N \delta_{X_i}$.
	Hence, the TT-rank of the resulting MPS ${\bm c}$ is bounded by $N$ multiplied by the TT-rank of $\bm{G}$.
	However, we observe that truncating the MPS ${\bm c}$ yields a TT-rank smaller than this pessimistic bound\footnote{This observation was first made by Nicolas Jolly.}, \textit{cf.} Figure \ref{Fig:RankC}.
	Empirically, we see that the rank scales like $N^\alpha$, with $\alpha \simeq 0.4$ in dimensions $d=2$ and $d=3$, with a relative error in $\ell^2$ norm of order $10^{-6}$.
	Such a scaling is unclear to us, but very favorable for simulating multiscale systems.
	
	We also remark that the TT-rank of ${\bm c}$ is correlated with the radius $r$ --the latter being approximately proportional to $N^{-1/d}$.
	We observe empirically that it scales like $r^{\alpha}$ with $\alpha = 0.7$ in dimension $d=2$ and $\alpha =  1.2$ in dimension $d=3$.
	This might be related to the result \cite[Cor. 10]{lindsey_multiscale_2024}, where the author proves that the TT-rank of 1-dimensional functions with Fourier support limited to $[-\Omega, \Omega]$ is of TT-rank at most $C\sqrt{\Omega}$.
	We expect this result to generalize to higher dimensions $d>1$; heuristically, one could argue that the scaling would then be replaced by $\Omega^{d/2}$.
	Since the Gaussian $G$ defined by \eqref{Def:G} is almost band-limited in $[C/r, C/r]^d$ for $C =10$, we would then obtain that the TT-rank should be limited by $r^{-d/2}$.
	This is still pessimistic, but definitely closer to the reality than the immediate bound in $r^d$ (which is related to linearity in $N$ of the pessimistic estimate on the TT-rank).

	\subsubsection{Protocol}\label{Sec:Protocol}
	
	Our protocol is the following: 
	We generate a coefficient $c$ as in Section \ref{Sec:microEmpirical}.
	Then, we choose $\epsilon=2^{-L_\epsilon}$, and a discretization step $h:=2^{-L}$ with $L=L_0+L_\epsilon$.
	In the sequel, $L_0$ is fixed, whereas $L_\epsilon$ or equivalently $\epsilon$ may vary.
	We consider two different grids: the smallest one represents only the microscopic periodic pattern and meshes the microscopic cube $[0,\epsilon]^d$ with discretization step $h=2^{-L}$, while the largest one represents the full coefficient field and meshes the macroscopic cube $[0, 1]^d$ with the same discretization step $h=2^{-L}$.
	We approximate $c(x/\epsilon)$ by an MPS ${\bm c}$ on the smallest grid by means of TCI (see Section \ref{Sec:microEmpirical}).
	Then, we periodize\footnote{This is an analytical manipulation, just setting trivial cores to scales higher than $\epsilon$.} ${\bm c}$ to get it on the largest grid, retrieving an approximation of $c(\cdot/\epsilon)$ on $[0, 1]^d$.
	Next, we analytically build $b$ and $g$ as MPS on the largest grid.
	
	We are now in position to solve the equation \eqref{Eeps} on the largest grid.
	To do so, we build by TCI the MPS $\bhP$, choose a penalty parameter $\mu$, and use the algorithm described in Section~\ref{Sec:Algo}.
	For definiteness, we chose the 1-site ALS numerical method \cite{holtz_alternating_2012} as a QTT linear solver, with a fixed largest TT-rank $R_0$, stopping the algorithm when the relative $\ell^2$ difference between states $\psi^{n}$ and $\psi^{n+1}$ after a sweep reaches a given relative tolerance ${\rm tol}$, or after the solver reaches a maximal number of sweeps denoted as $S$.
	With the same parameters, we solve the corrector equation \eqref{E:corr} on the smallest grid, deduce the homogenized matrix $\abar$ by evaluating \eqref{Def:abar}.
	Thus, we may solve the homogenized equation \eqref{E-hom} on the largest grid.\footnote{We do not need the fine discretization $h=2^{-L}$ there because we \textit{a priori} know that $\ubar$ is regular.
	But, since the computational overhead is not high, for simplicity, we solve \eqref{E-hom} as well on the largest grid.}
	
	We validate our numerical resolutions by two different ways.
	First, we measure the \textit{a posteriori} errors defined by \eqref{Def:EP}, \eqref{Error-total}, and \eqref{Error-total2}.
	Second, we also evaluate the numerical errors related to homogenization corresponding to
	\begin{align}\label{Def:Ehomnabla}
		&E_{\rm hom}[\nabla \ueps] :=\frac{\|\nabla \ueps - \nabla  u^{\epsilon, 1}\|_2}{\|\nabla \ueps\|_2}
		\qquad \et \quad
		E_{\rm hom}[\ueps] :=\frac{\| \ueps - \ubar\|_2}{\|\ueps\|_2}.
	\end{align}
	As explained in \eqref{CvRate}, we expect these to scale as $\epsilon$.
	We cannot directly evaluate these, since we do not have access to the actual $\ueps$ and $\nabla \ueps$; thus, we rather evaluate their numerical approximations, that we denote $E_{\rm hom}[\psi^\epsilon]$ and $E_{\rm hom}[u_{\psi^\epsilon}]$, where $u_{\psi^\epsilon}$ is deduced from $\psi^\epsilon$ as explained in Section \ref{Sec:getu}.

\subsubsection{Settings}

We make use of different parameters in dimensions $d=2$ and $d=3$; obviously, those of the $3$D case are less ambitious, since the problem is more complex.
We refer to Table \ref{Tab:Param2D} for parameter values that are fixed across all the subsequent numerical experiments.

\begin{table}[h]
	\begin{center}
		{\renewcommand{\arraystretch}{1.5}
			\begin{tabular}{|l|l|l|l|l|l|l|l|l|l|l|}
				\hline
				Parameter symbol & $d$ & $r$ & $N$  & $\nu$ & Format & $L_0$& $\mu$ & Initial maximal rank of $\bhpsi$ & $S$ & {\rm tol}
				\\
				\hline
				Value & $2$ & $0.05$ & $10$ & $5$ & \texttt{x1y1} & $15$ &  $1\cdot10^4$ & $70$ & $40$ & $1\cdot10^{-6}$
				\\
				\hline
				Value & $3$ & $0.2$ & $2$ & $5$ & \texttt{x1y1} & $12$ &  $1\cdot10^4$ & $100$ & $30$ &$1\cdot10^{-6}$
				\\
				\hline
				
			\end{tabular}
		}
	\end{center}
	\caption{Parameters for multiscale experiments in dimensions $d=2$ and $d=3$.}\label{Tab:Param2D}
\end{table}

\begin{figure}[h]
	\begin{center}
		\begin{tabular}{l l l}
			a) & b) & c)
			\\
			\includegraphics[width=0.3\textwidth, clip, trim=80 0 0 0]{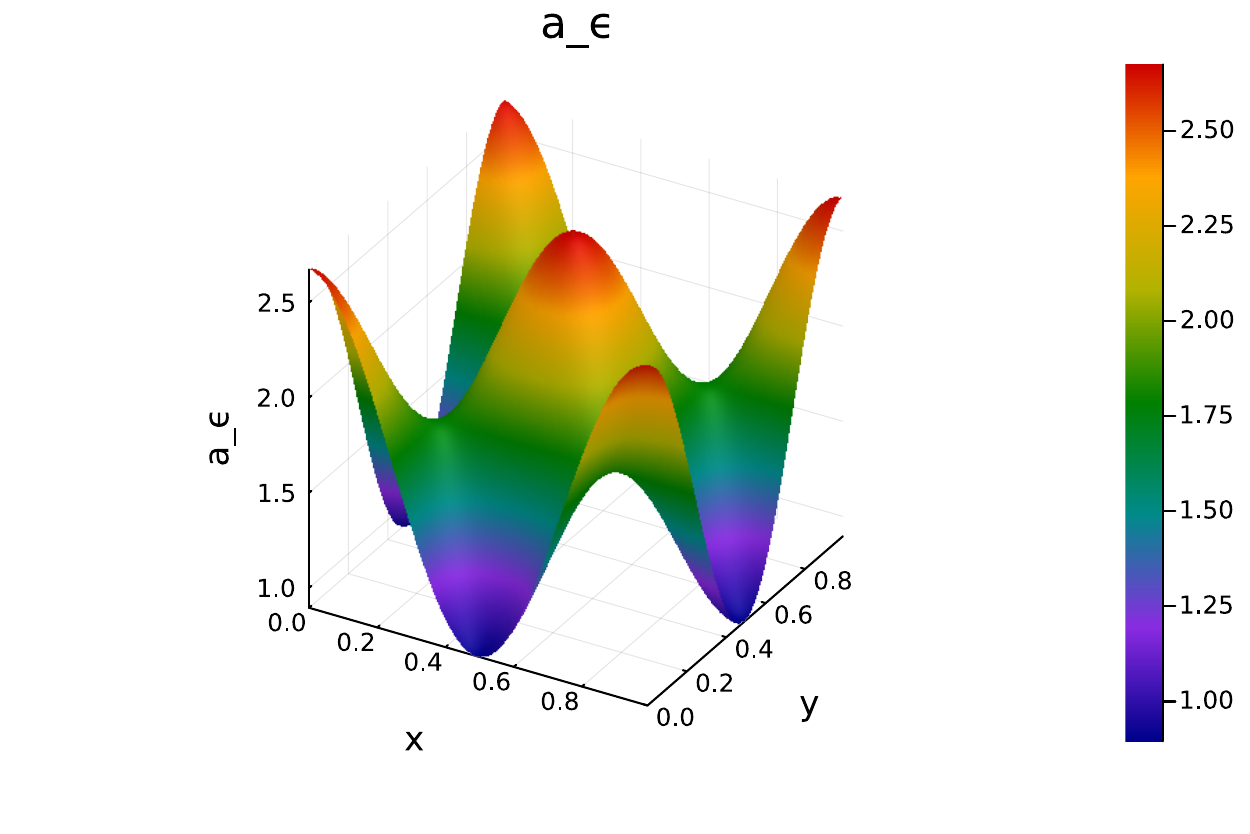}
			&\includegraphics[width=0.3\textwidth, clip, trim=80 0 0 0]{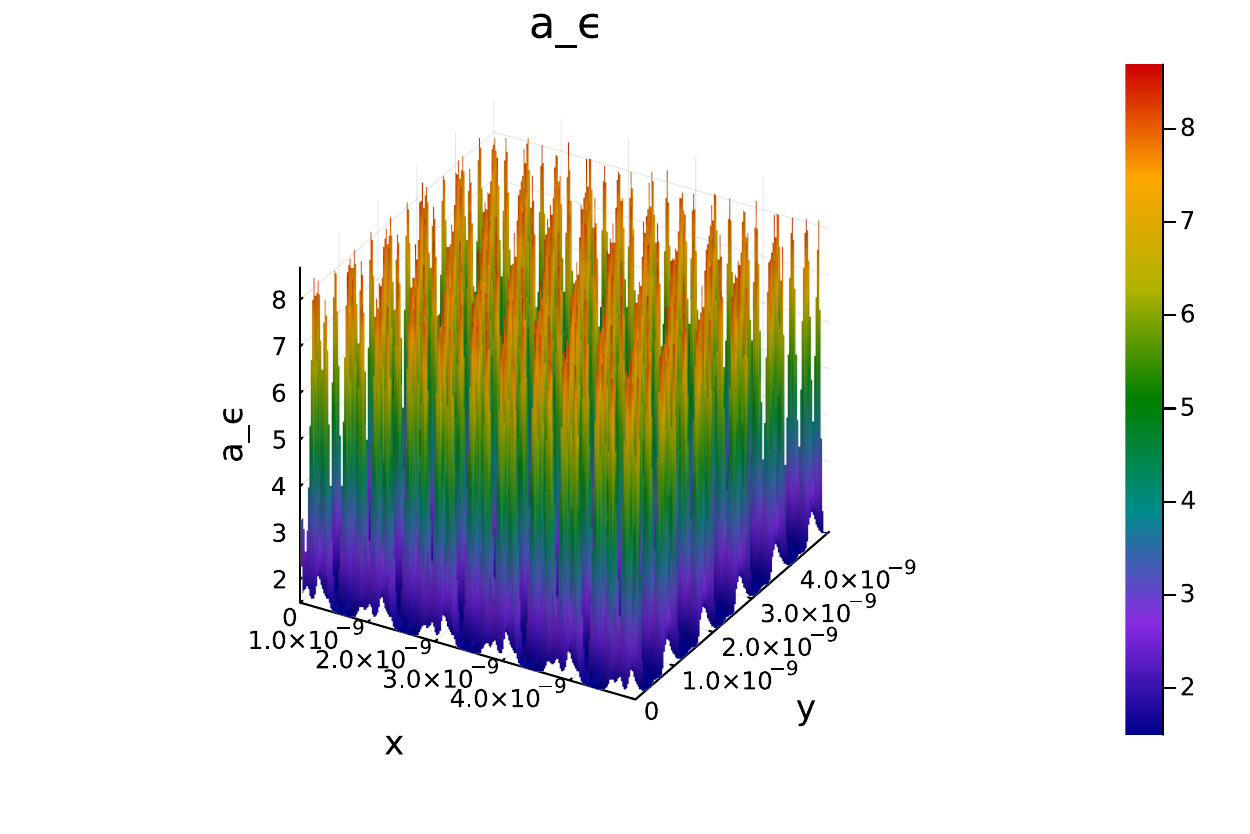}
			&\includegraphics[width=0.3\textwidth, clip, trim=80 0 0 0]{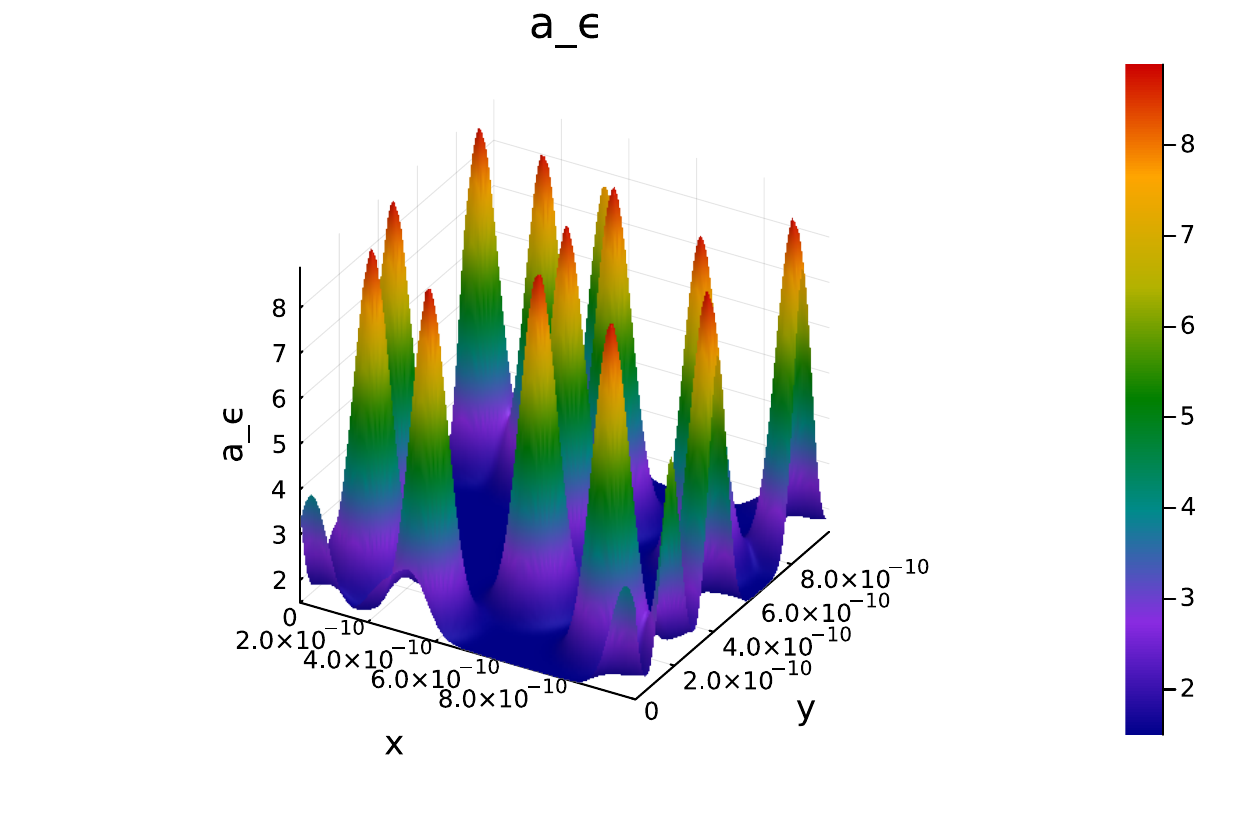}
		\end{tabular}
	\end{center}
	\caption{For $\epsilon \simeq 10^{-9}$: a), b), c) represent  $a_\epsilon$ at scales $1$,  $5\epsilon$, and $\epsilon$, \textit{i.e.} in $[0,1]^2$, $[0,5\epsilon]^d$, $[0,\epsilon]^2$, respectively.
	}\label{Fig:ag}
\end{figure}

\paragraph{Settings in $2$D}
	In dimension $d=2$, we choose a Volume Element containing $N=10$ different centers $X_i$, with radius parameter $r=0.05$.
	We make $\epsilon$ vary between $2^{-30} \simeq  1\cdot 10^{-9}$ and $1/4$.
	We underline that the smallest discretization step is of $h  = 2^{-45} \simeq 2 \cdot 10^{-14}$, which corresponds to $2^{90} \simeq 1\cdot 10^{27}$ virtual discretization points.
	Choosing $\epsilon=1\cdot 10^{-9}$, we show the coefficient $a_\epsilon$ at scales $1$, $5\epsilon$, and $\epsilon$
	in Figure \ref{Fig:ag}.
	Notice that, when showing $a$ at scale $1$  in Figure \ref{Fig:ag}a), we average the coefficient under the pixel scale, so that its microscopic features cannot be seen.
	Yet, we emphasize that $a_\epsilon$ involves through $c_\epsilon$ no fewer than $10^{18}$ periodic copies of the pattern shown in Figure \ref{Fig:ag}c) modulated by $b$.

\paragraph{Settings in $3$D}
	In dimension $d=3$, we choose a Volume Element containing $N=2$ different Gaussian centers, with radius parameter $r=0.2$.
	We make $\epsilon$ vary between $2^{-30} \simeq  1\cdot 10^{-9}$ and $1/4$, and apply the protocol of Section \ref{Sec:Protocol}.
	The smallest discretization step is of $h  = 2^{-42} \simeq 2\cdot 10^{-13}$, which corresponds to $2^{126} \simeq 8\cdot 10^{37}$ virtual discretization points.

\subsubsection{Discussion on QTT results} \label{Sec:2D-3D}
	
	First, both in dimensions $d=2$ and $d=3$, down to $\epsilon=10^{-9}$, we obtain meaningful results.
	The algorithm converges for the prescribed tolerance only for $\phi_i$ and $\ubar$ in dimension $d=2$, but it hits the maximum number of sweeps for other cases.
	This is presumably due to an insufficient maximal rank with respect to the prescribed tolerance.
	However, even though the latter is not reached, we may still evaluate the quality of the output by means of the \textit{a posteriori} estimators.
	
	Obviously, the algorithm is slower and more costly in terms of memory in dimension $d=3$ rather than in dimension $d=2$:
	indeed, we need to increase the maximal TT-rank of the solution from $70$ to $100$ to get accurate results, and increase as well the length of the TT for fixed discretization step (because of the new cores corresponding to indices $z_1$, $\dots$, $z_L$).
	Nevertheless, as shown in the subsequent results, for the chosen parameters, we do not see much differences between the dimensions $d=2$ and $d=3$ for our test case.
	Hence, we showcase our results for both dimensions in parallel.
	
	\medskip
	
\begin{figure}[h]
		\begin{center}
			\begin{tabular}{l l}
				\begin{tabular}{c}
					$d=2$
					\\
					\includegraphics[width=0.4\textwidth]{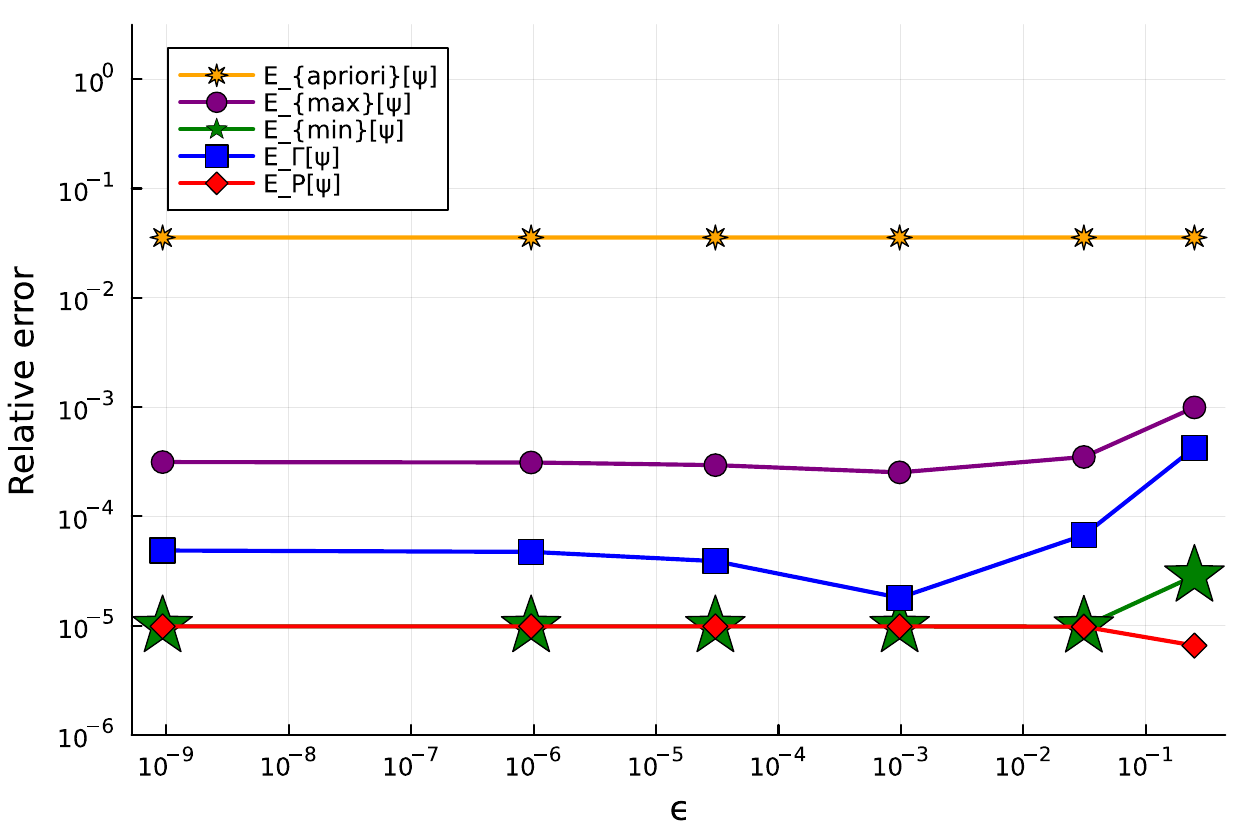}
				\end{tabular}
				\begin{tabular}{c}
					$d=3$
					\\
					\includegraphics[width=0.4\textwidth]{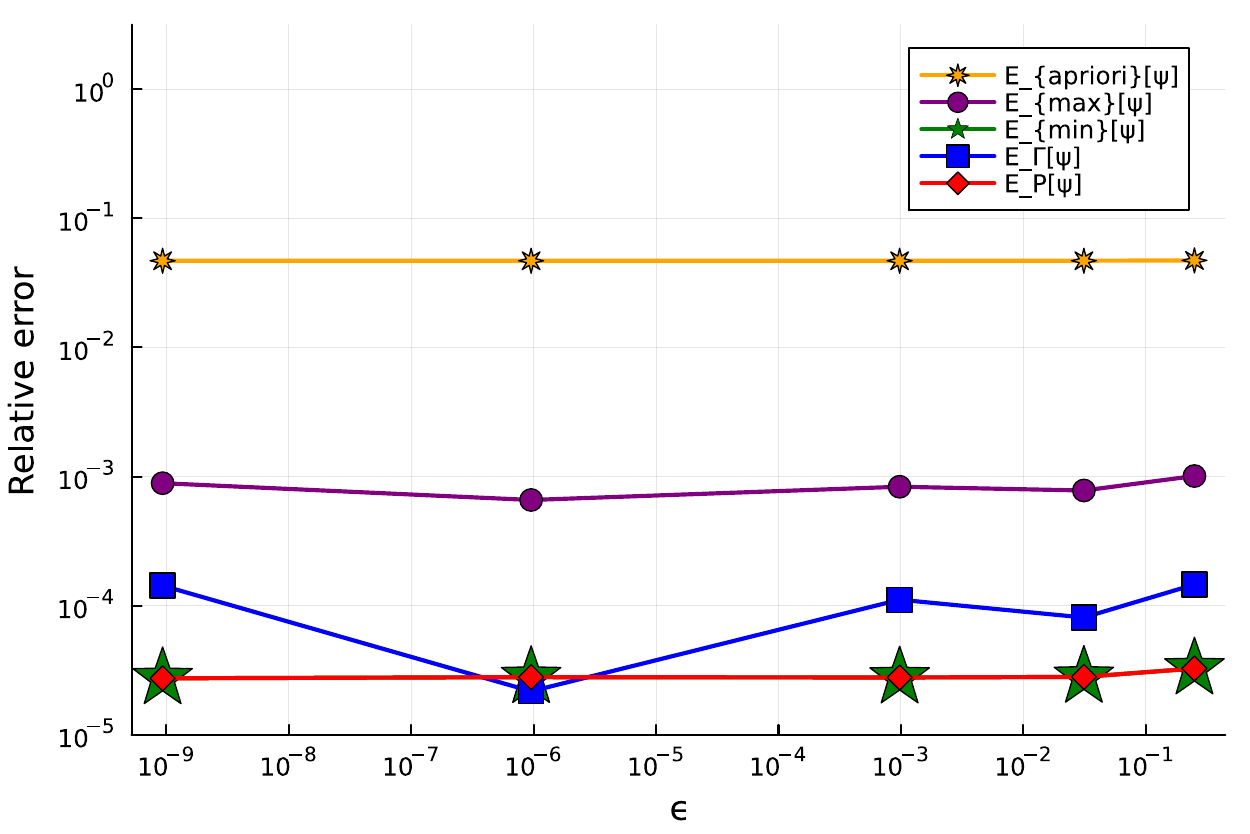}
				\end{tabular}
			\end{tabular}
		\end{center}
		\caption{
			Cases $d=2$ on the left and $d=3$ on the right of the multiscale setting of Section \ref{Sec:microEmpirical}, where parameters are given by Table \ref{Tab:Param2D}.
			On both pictures: Numerical errors $E_P[\psi^\epsilon]$, $E_\Gamma[\psi^\epsilon]$, \textit{a priori} estimator $E_{\rm apriori}[\psi^\epsilon]$ given by \eqref{E_apriori}, and \textit{a posteriori} estimators $E_{\max}[\psi^\epsilon]$ and $E_{\min}[\psi^\epsilon]$, given by  \eqref{Def:EP}, \eqref{Error-total}, and \eqref{Error-total2}, as a function of $\epsilon$.}
			\label{Fig:Errors:2D3D}
	\end{figure}

	In both dimensions, the \textit{a posteriori} error estimates safeguard the accuracy of our numerical solutions.
	These results are displayed in Figure \ref{Fig:Errors:2D3D}.
	We observe that the overall \textit{a posteriori} error estimator $E_{\max}[\psi^\epsilon]$ is between $10^{-4}$ and $10^{-3}$ in dimension $d=2$ and slightly below than $10^{-3}$ in dimension $d=3$, which is satisfactory for most applications in thermomechanics.
	In both cases, it reaches a plateau for $\epsilon \leq 10^{-2}$.
	This plateau error for $E_{\max}[\psi^\epsilon]$ is of $3\cdot 10^{-4}$ for $d=2$ and $1\cdot10^{-4}$ for $d=3$, which is better than that of the \textit{a priori} error estimator
	\begin{align}\label{E_apriori}
		E_{\rm apriori}[\psi^\epsilon] := 
		\frac{(\lambda+\Lambda)(1+\lambda\Lambda)}{\mu+\lambda^{-1}}\frac{\|g\|_2}{\|\nabla \ueps\|_2} \geq \frac{\|\psi^\epsilon - \nabla \ueps\|_2}{\|\nabla \ueps\|_2}
	\end{align}
	deduced from \eqref{num:0005}, which evaluates at approximately $3\cdot 10^{-2}$ and $5\cdot 10^{-2}$ for the cases $d=2$ and $d=3$, respectively.
	(In practice, $\psi^\epsilon$ is used instead of $\nabla \ueps$ in order to evaluate the upper bound of \eqref{E_apriori} --which implies a slight discrepancy in the bound value.)
	Also, we remark that the \textit{a posteriori} error estimator $E_{\min}[\psi^\epsilon]$ is not very far from $E_{\max}[\psi^\epsilon]$, so that the whole pair gives a useful hint on the real error approximation on $\nabla\ueps$ of the solver.
	
	Taking into account \eqref{Error-total}, we remark that the contributions of $E_P[\psi^\epsilon]$ and $E_\Gamma[\psi^\epsilon]$ to the total error estimator $E_{\max}[\psi^\epsilon]$ are roughly of the same order: 
	the fraction of $E_{\max}[\psi^\epsilon]$ due to $E_P[\psi^\epsilon]$ is between $15\%$ and $85\%$ in dimension $d=2$, and between $67\%$ and $93\%$ in dimension $d=3$.
	The contribution $E_P[\psi^\epsilon]$ is remarkably stable when $\mu$ is fixed, and we checked that it decays when $\mu$ increases (see Section \ref{Sec:QTTHL} for a discussion on the effect of $\mu$ on the error, and the balance between accuracy and condition number).
	On the contrary, we checked that $E_\Gamma[\psi^\epsilon]$ decays when we increase the prescribed maximal TT-rank of $\bhpsi$.

	Hence, for these special experiments, the chosen maximal TT-ranks are close to be optimal with respect to the fixed $\mu$ for most values of $\epsilon$, since $E_P[\psi^\epsilon]$, $E_\Gamma[\psi^\epsilon]$ are roughly of the same order.
	However, interestingly enough, in dimension $d=2$, the estimator $E_{\max}[\psi^\epsilon]$ decays with $\epsilon$ until it reaches the plateau.
	For the largest $\epsilon$, the contribution of $E_\Gamma[\psi^\epsilon]$ to $E_{\max}[\psi^\epsilon]$ dominates (from $35\%$ for $\epsilon=10^{-9}$, it reaches $85\%$ for $\epsilon=1/4$).
	We interpret this as a consequence of the fact that scale separation barely happens when $\epsilon=1/4$, so that the fixed rank of $70$ might be insufficient in that case, while penalization maintains a small $E_P[\psi^\epsilon]$.
		
	\subsubsection{Validation and comparison with the homogenization results}
	
\begin{figure}[h]
		\begin{center}
			\begin{tabular}{l l l}
			a) & b) & c)
			\\
			\includegraphics[width=0.3\textwidth, clip, trim=80 0 0 0]{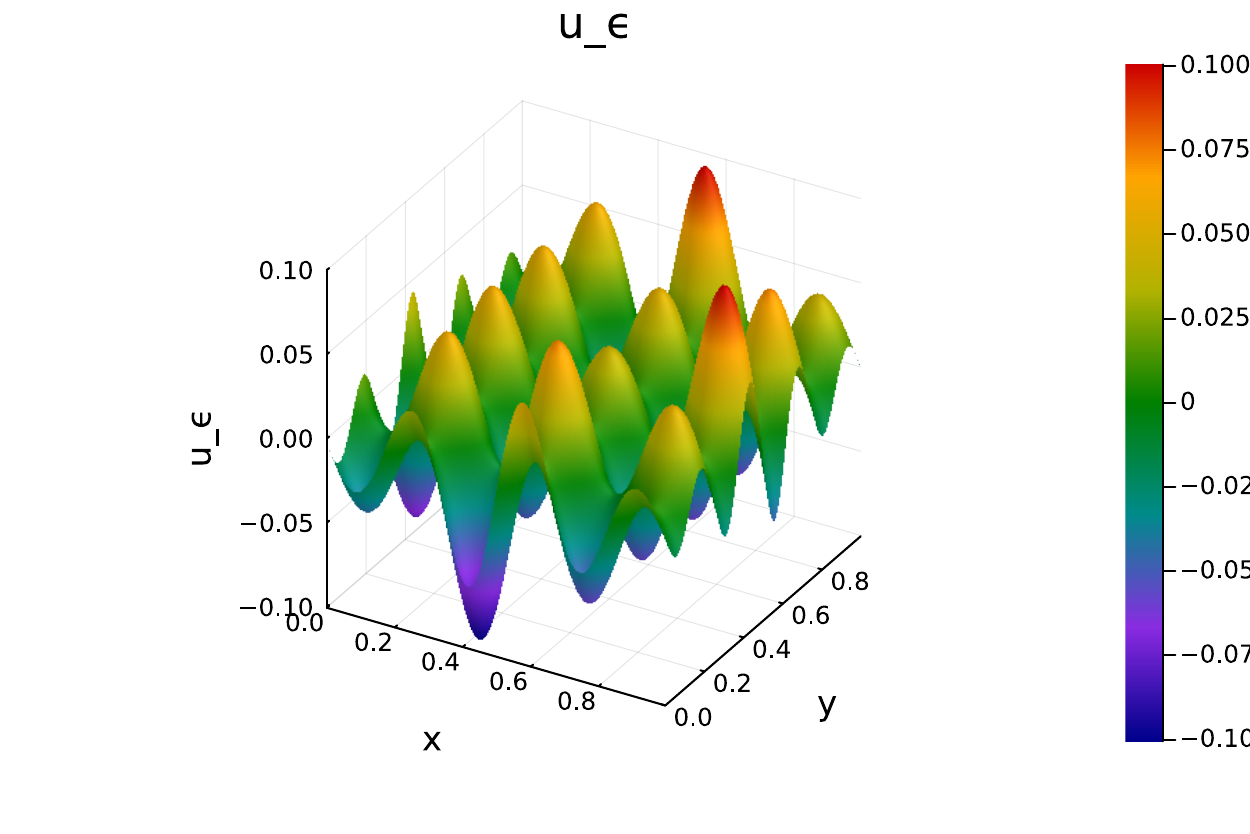} &
			\includegraphics[width=0.3\textwidth, clip, trim=80 0 0 0]{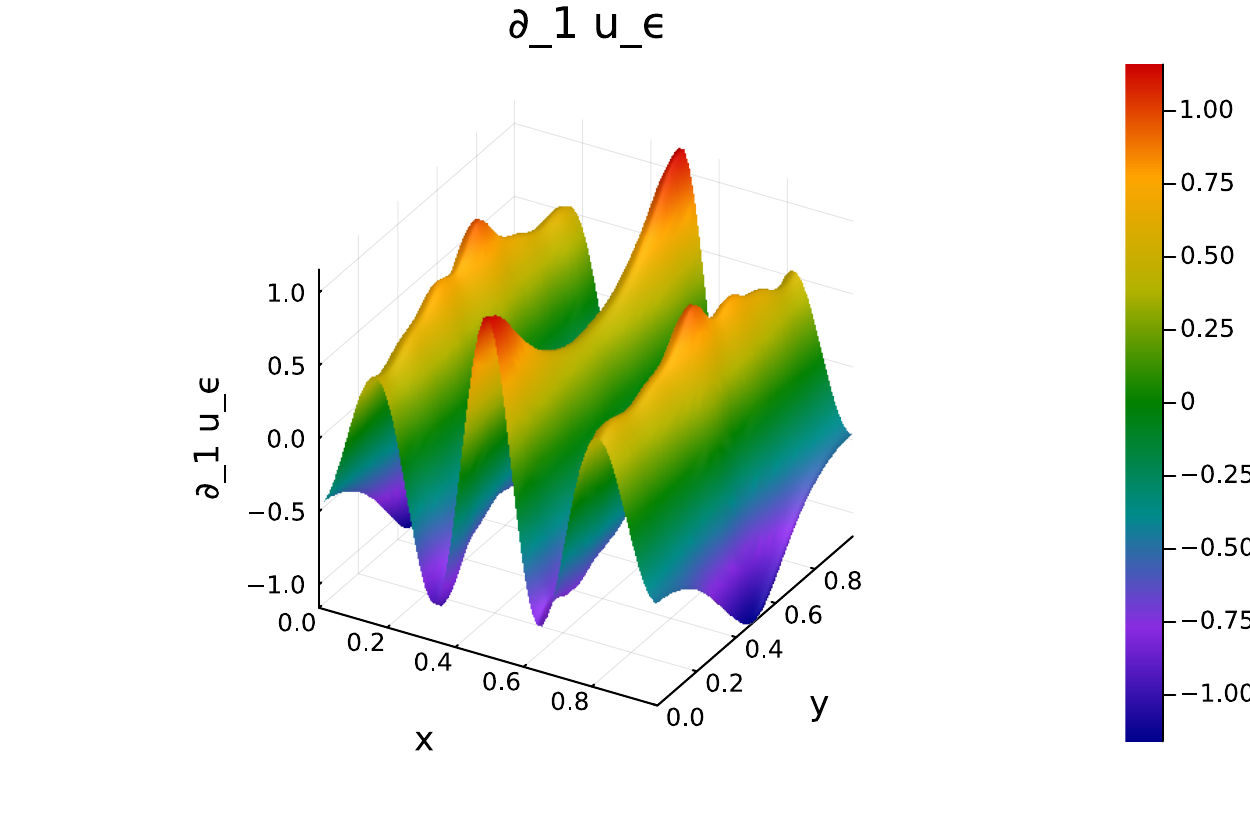}&
			\includegraphics[width=0.3\textwidth, clip, trim=80 0 0 0]{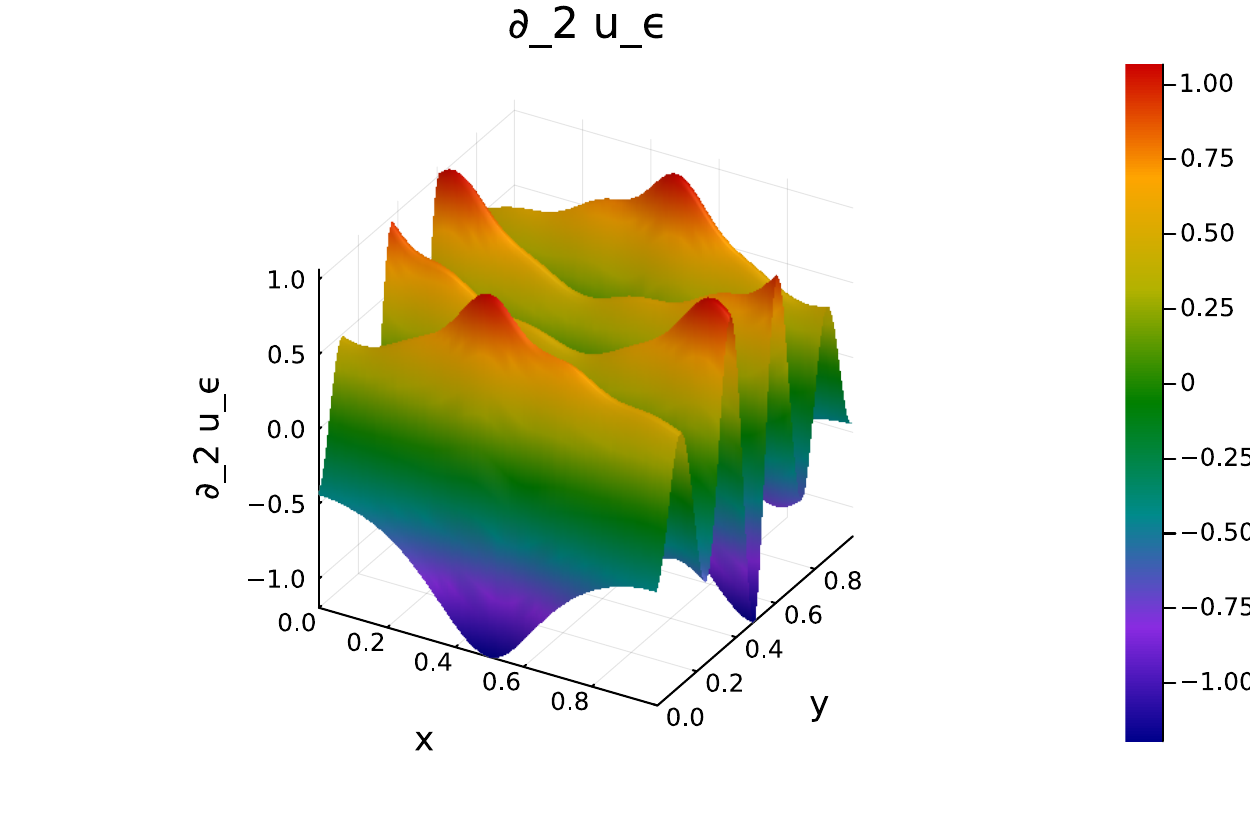}
			\\
			d) & e)
			\\
			\includegraphics[width=0.3\textwidth, clip, trim=80 0 0 0]{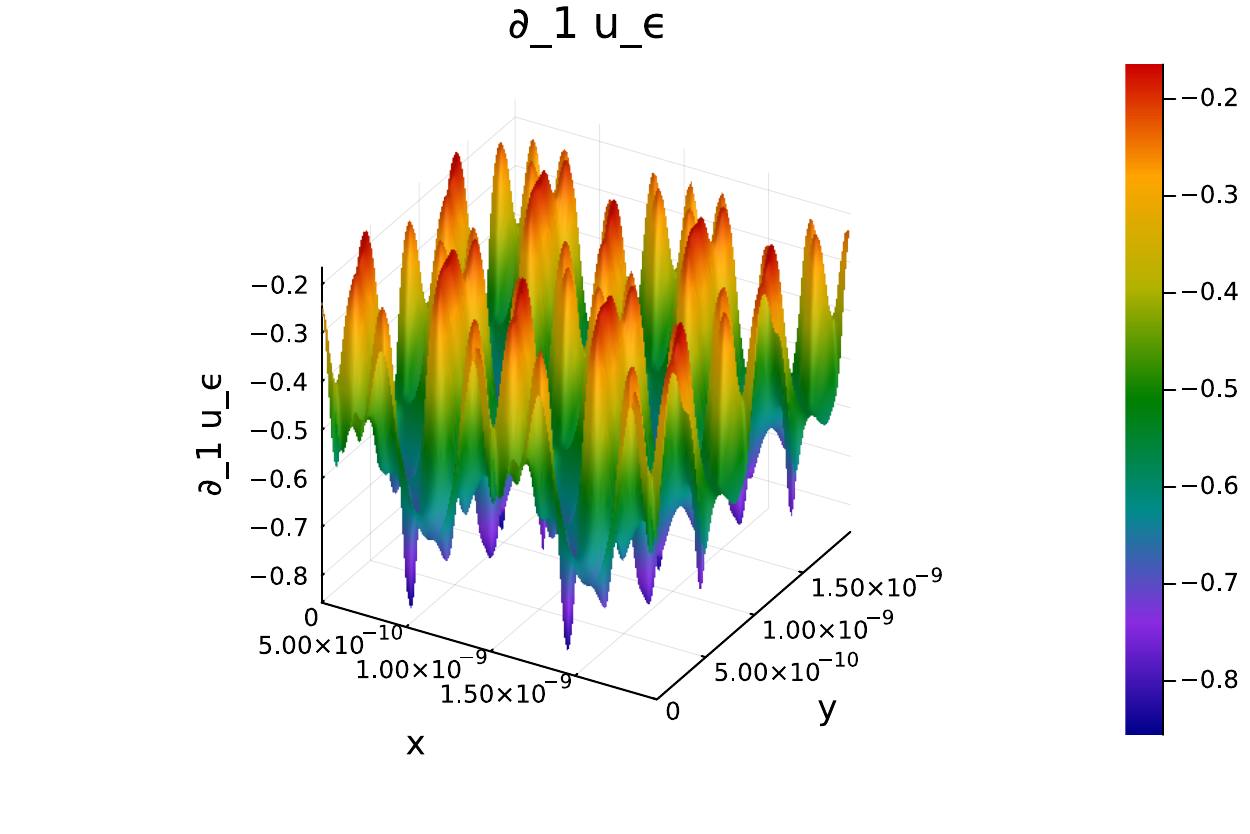}&
			\includegraphics[width=0.3\textwidth, clip, trim=80 0 0 0]{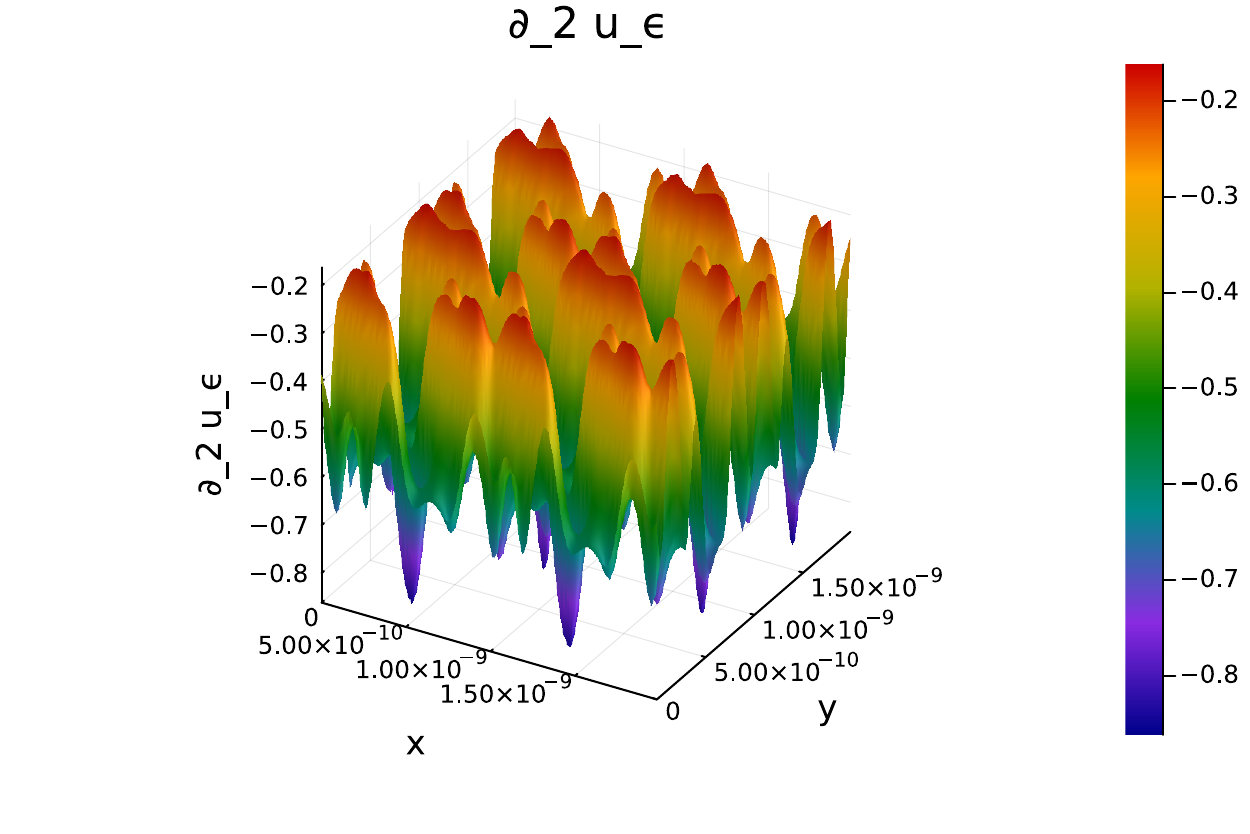}
			\end{tabular}
		\end{center}
		\caption{
		Case $d=2$ of the multiscale setting of Section \ref{Sec:microEmpirical}, where parameters are given by Table \ref{Tab:Param2D}, for $\epsilon \simeq 10^{-9}$ and $h\simeq 2\cdot 10^{-14}$.
		a), b) and c) represent $\ueps$, $\partial_1 \ueps$ and $\partial_2 \ueps$ at scale $1$;
		d) and e)  represent $\partial_1 \ueps$ and $\partial_2 \ueps$ at scale $2\epsilon$.}\label{Fig:res}
	\end{figure}

	\paragraph{Observing the solution}
	In dimension $d=2$, we draw $\ueps$ and $\nabla\ueps$ at different scales in Figure \ref{Fig:res}, for $\epsilon \simeq 10^{-9}$ and $h\simeq 2\cdot 10^{-14}$.
	As expected, the function $\ueps$ is smooth --indeed, it is close to the homogeneous solution $\ubar$ with error $\epsilon \simeq 10^{-9}$.
	On the contrary, at scale $\epsilon$, $\nabla\ueps$ features oscillations of bounded amplitude and of wavelengths smaller than $\epsilon$ correlated with the variations of $\aeps$ at that scale.
	As expected, we also observe that it shows a periodic pattern of size $[0, \epsilon]^2$ --this can be observed in Figure \ref{Fig:res} d) and e), where we plotted it on the small domain $[0,2\epsilon]^2$.
	This is in accordance with \eqref{Uncoupling}, where the $\Q_{\epsilon,\per}$-periodic functions $\nabla\phi_i(\cdot/\epsilon)$ are modulated by the macroscopic gradient $\nabla\ubar$ to give rise to an approximation of $\nabla\ueps$.
		
\begin{figure}[h]
		\begin{center}
			\begin{tabular}{c}
				$d=2$
				\\
				\includegraphics[width=0.4\textwidth]{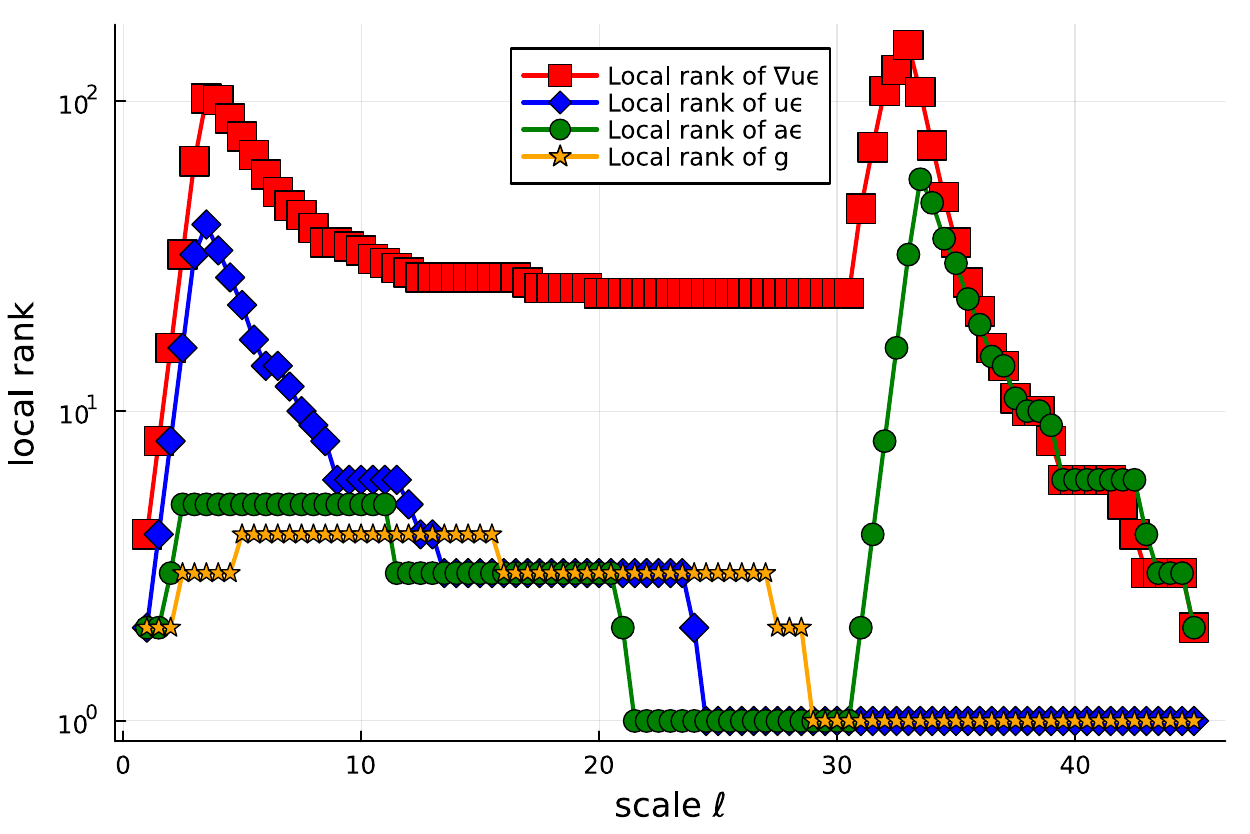}
			\end{tabular}
			\begin{tabular}{c}
				$d=3$
				\\
				\includegraphics[width=0.4\textwidth]{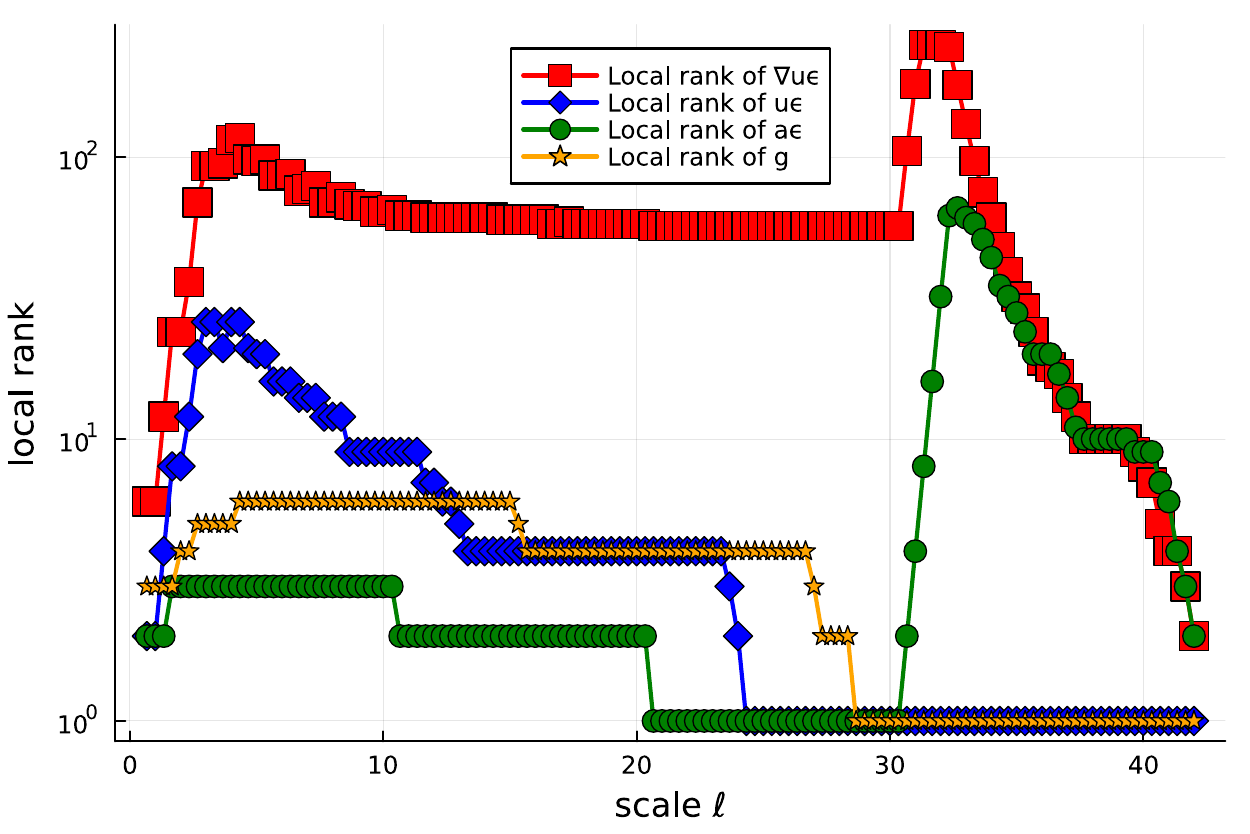}
			\end{tabular}
		\end{center}
		\caption{
			Cases $d=2$ on the left and $d=3$ on the right of the multiscale setting of Section \ref{Sec:microEmpirical}, where parameters are given by Table \ref{Tab:Param2D},  for $\epsilon = 2^{-30} \simeq 10^{-9}$.
			On both pictures: 
			Local ranks $r_\ell$ of the MPS representations of $\nabla\ueps$, $\ueps$, $a_\epsilon$, $g$ in red, blue, green, and orange as a function of the dyadic scale $\ell$.
		}\label{Fig:Rank}
	\end{figure}

	\paragraph{TT-ranks and scale separation}	
	We remark that the TT-ranks of the MPS representations of $\nabla\ueps$, $\ueps$, $a_\epsilon$, and $g$ reflect the underlying homogenization process.
	Indeed, in Figure~\ref{Fig:Rank}, we show the local TT-ranks $r_\ell$ of the MPS $\nabla\ueps$, $\ueps$, $a_\epsilon$, and $g$ in the case $\epsilon=10^{-9}$ as a function of the dyadic scale parameter $\ell$ (the larger $\ell$, the smaller the scale $2^{-\ell}$).
	As expected, for both dimensions $d=2$ and $d=3$, we observe that
	\begin{itemize}
		\item the local ranks $r_\ell$ of the MPS representing $g$ are small and concentrate in the largest scales (the latter feature is related to the fact that $g$ is regular, so that dependence on small scales is removed by TT-rounding);
		\item $\aeps$ is represented by an MPS with two peaks in terms of local TT-ranks $r_\ell$, for $\ell=5$ and $\ell = 35$, separated by a well of smaller ranks $r_\ell$ for $\ell \in\{10,\dots, 30\}$.
		The first is for the caracteristic scales of $b$, located in the largest dyadic scales $\ell \leq 20$, and the second is related to the microstructure $c(\cdot/\epsilon)$, with local ranks located in $\ell \geq \log_2(\epsilon)=30$;
		\item the MPS representing $\nabla\ueps$, similarly to $\aeps$, has local TT-ranks $r_\ell$ with two peaks near $\ell=5$ and $\ell = 35$, separated by a well of smaller ranks $r_\ell$ for $\ell \in\{10,\dots, 30\}$. This reflects the two-scale expansion \eqref{Uncoupling};
		\item the MPS representing $\ueps$ has TT-ranks $r_\ell>1$ only for $\ell \lesssim 25$, and TT-ranks $r_\ell=1$ for $\ell > 25$. This reflects that $\ueps$ is well-approximated by a regular function $\ubar$.
	\end{itemize}
	We underline that the small-scale oscillations in $\nabla\ueps$ are not present in the MPS representing $\ueps$ which is of TT-ranks $r_\ell=1$ for small dyadic scales $\ell > 25$.
	Indeed, these are suppressed by the truncation process, because their $\ell^2$ norm is too small.
	This advocates for our strategy of computing first $\nabla\ueps$ and then deducing $\ueps$ from the latter --since the other way round is impossible, because of the loss of small-scale informations in $\ueps$.
	
	\begin{figure}[h]
		\begin{center}
			\begin{tabular}{c}
				$d=2$
				\\
				\includegraphics[width=0.4\textwidth]{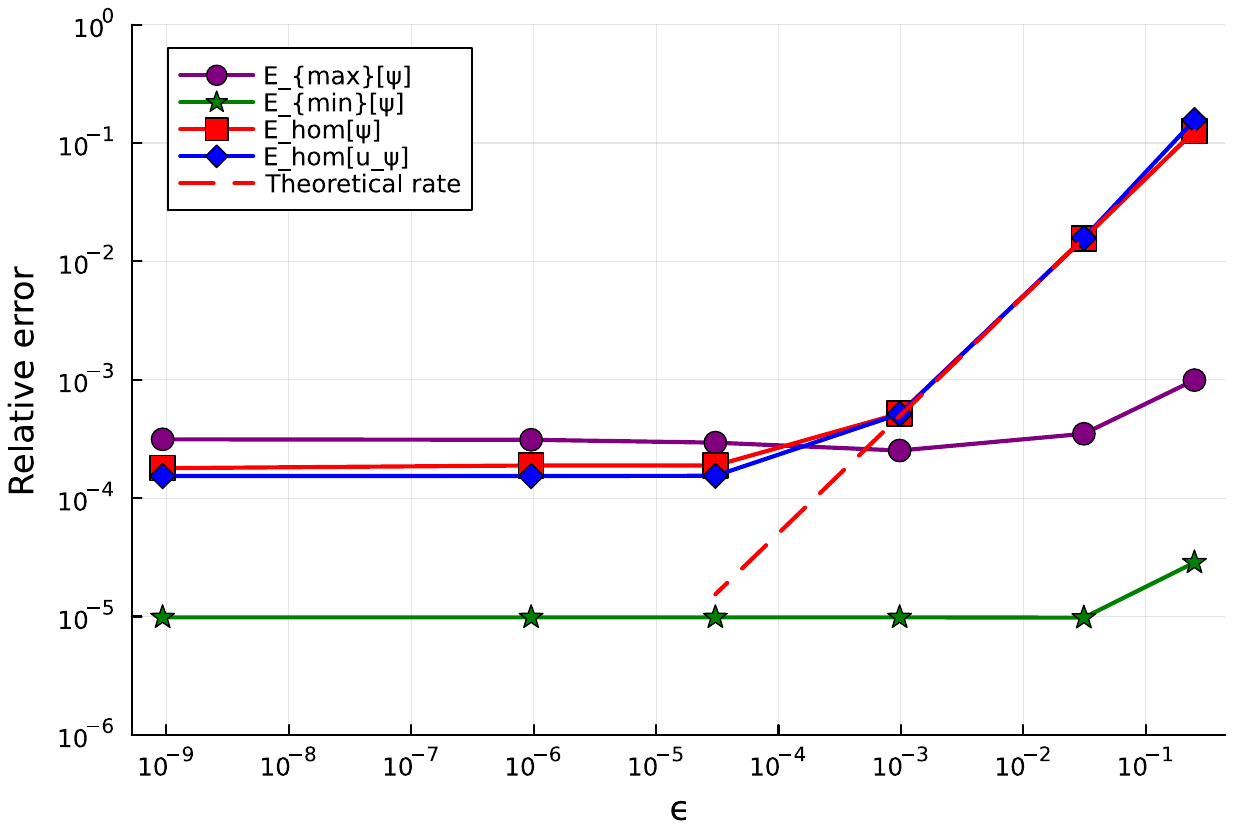}
			\end{tabular}
			\begin{tabular}{c}
				$d=3$
				\\
				\includegraphics[width=0.4\textwidth]{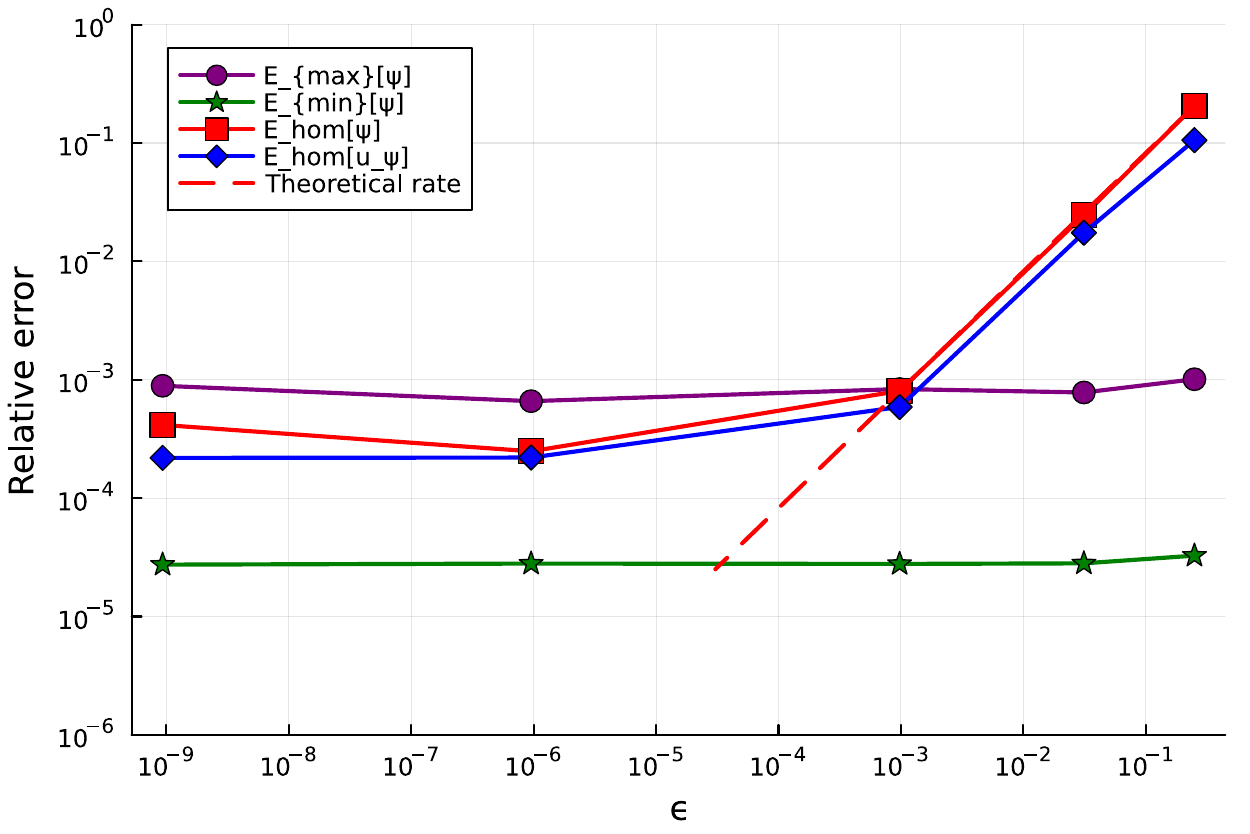}
			\end{tabular}
		\end{center}
		\caption{
			Cases $d=2$ on the left and $d=3$ on the right of the multiscale setting of Section \ref{Sec:microEmpirical}, where parameters are given by Table \ref{Tab:Param2D}.
			On both pictures: Numerical homogenization errors $E_{\rm hom}[\psi^\epsilon]$, $E_{\rm hom}[u_{\psi^\epsilon}]$, and \textit{a posteriori} estimates $E_{\max}[\psi^\epsilon]$ and $E_{\min}[\psi^\epsilon]$, given by \eqref{Def:Ehomnabla}, \eqref{Error-total} and \eqref{Error-total2} as a function of $\epsilon$.}\label{Fig:Homog}
	\end{figure}
	
	\paragraph{Convergence rates in homogenization}	
	We now validate our results with the quantitative predictions of the homogenization theory.
	We plot $E_{\rm hom}[\psi^\epsilon]$ and $E_{\rm hom}[u_{\psi^\epsilon}]$ defined by \eqref{Def:Ehomnabla} in Figure \ref{Fig:Homog}.
	We observe that the result is coherent with the homogenization theory: as long as the precision on the $\ueps$ and $\nabla \ueps$ is sufficient, we observe convergence rates $E_{\rm hom}[\psi^\epsilon]$ and $E_{\rm hom}[u_{\psi^\epsilon}]$ that are linear in $\epsilon$.
	Then, for $\epsilon$ sufficiently small, these quantities saturate at a value that is coherent with the \textit{a posteriori} error.
	(Notice that similar numerical errors are also made when computing for $\nabla \phi_i$ and $\nabla \ubar$, but not shown here.)

	\subsubsection{Comparing the different formats}\label{Sec:OtherFormats}
	
	All the above results in Section \ref{Sec:homog} were obtained with the specific format \texttt{x1y1}, which is particularly fit for representing the solution to this multiscale test case, as explained below.
	Yet, we could have chosen a different format in Table \ref{Tab:Format2}, and obtain same qualitative properties at the price of poorer error approximations.
	We exemplify this by considering the most demanding case for $d=2$ --that is $L=45$ and $\epsilon=2^{-30}$-- and considering all the formats of Table \ref{Tab:Format2}.
	
	We show the results in Table \ref{Tab:ComparisonHomog}.
	We observe that it is indeed possible to get an accurate approximation of the solution, given the maximal rank $r=70$ of the solution $\bhpsi$.
	Indeed, for all the formats, the \textit{a posteriori} estimator $E_{\rm max}[\psi^\epsilon]$ remains below $4\%$.
	This is still an acceptable accuracy for many practical problems.
	Accordingly, the validation against the homogenization theory yields an error $E_{\rm hom}[\psi^\epsilon]$ below $5\%$ --we recall that the numerical evaluation of the homogenization error is affected by the approximation errors of its constituents.
	
	\begin{table}[h]
		\begin{center}
			{\renewcommand{\arraystretch}{1.5}
				\begin{tabular}{|c|c|c|c|}
					\hline
					Format & $E_{\max}[\psi^\epsilon]$ & $E_{\min}[\psi^\epsilon]$ & $E_{\rm hom}[\psi^\epsilon]$
					\\
					\hline 
					\texttt{x1y1} & $ 3.14 \times 10^{-4} $ & $ 9.83 \times 10^{-6} $ & $ 1.79 \times 10^{-4} $ \\ 
					\hline 
					\texttt{x1x2\_y2y1} & $ 2.20 \times 10^{-2} $ & $ 7.30 \times 10^{-4} $ & $ 4.02 \times 10^{-2} $ \\ 
					\hline 
					\texttt{x1x2\_y1y2} & $ 3.24 \times 10^{-2} $ & $ 1.08 \times 10^{-3} $ & $ 4.44 \times 10^{-2} $ \\ 
					\hline 
					\texttt{x2x1\_y1y2} & $ 1.05 \times 10^{-2} $ & $ 3.46 \times 10^{-4} $ & $ 1.92 \times 10^{-2} $ \\ 
					\hline
				\end{tabular}
			}
		\end{center}
		\caption{
			Errors $E_{\max}[\psi^\epsilon]$, $E_{\min}[\psi^\epsilon]$ and $E_{\rm hom}[\psi^\epsilon]$ for various formats, for the case $d=2$ of the multiscale setting of Section \ref{Sec:microEmpirical}, where parameters are given by Table \ref{Tab:Param2D}, and with $\epsilon=2^{-30} \simeq 10^{-9}$, $L=45$.}\label{Tab:ComparisonHomog}
	\end{table}
	
	However, we also see that the format \texttt{x1y1} is quantitatively superior to the others, gaining $2$ orders of magnitude in terms of both the \textit{a posteriori} estimator and accordingly the homogenization error.
	This can be explained by the fact that the approximation \eqref{Uncoupling}, which is well-adapted to the format \texttt{x1y1}, is less fit for other formats of Table \ref{Tab:Format2}.
	On the one hand, for all formats, exponential compression holds.
	On the other hand, for the three first formats of Table \ref{Tab:Format2}, this decomposition results in \textit{multiplying} the ranks for the approximation of $\nabla \ubar$ and $\nabla \phi_i$ --whereas, when scales are separated, this results in taking the maximum of the ranks of  $\nabla \ubar$ and $\nabla \phi_i(\cdot/\epsilon)$ for the format \texttt{x1y1}.
	
	This validates \textit{a posteriori} our choice of favoring the format \texttt{x1y1} for our special problem.
	Nevertheless, this also highlights that our approach is versatile in the sense that it does not rely on a single format, unlike the approaches of \cite{bachmayr_stability_2020} or \cite{rakhuba_robust_2021}, which depend on the specificities of the format they choose.
	
	\section{Conclusions}
        
        In this article, we have introduced a novel QTT-based numerical method for solving elliptic equations. As recently highlighted in the literature, the QTT structure enables discretization levels far beyond those reachable with classical solvers, making this formalism particularly attractive, in particular for multiscale problems.
        The proposed QTT solver overcomes several limitations of state-of-the-art QTT solvers for solving elliptic equations, notably regarding coefficient heterogeneity, problem dimension, and solver stability.

        The key ingredients of the solver are as follows. First, the introduction of the Helmholtz--Leray decomposition to reformulate the elliptic problem as a minimization problem posed on a new variable generalizing the solution gradient, under the constraint that this new variable remains a gradient. This constraint is enforced through a penalization term. Second, the Fourier transform is used to express the Helmholtz--Leray projector analytically. Finally, the QTT representation is employed to formulate the discretized problem, allowing for extensive use of the available QTT toolbox, among which the ALS solver, the QFT algorithm and the TCI algorithm.
        Setting aside the QTT part, our solver looks like a variant of the celebrated FFT-based solvers.
        From a more abstract point of view, its discretized QTT formulation preserves important algebraic and topological structures involved in the variational formulation of the elliptic equation.
        
        The mathematical properties of this QTT Helmholtz--Leray (QTT-HL) solver have been investigated, showing in particular that its unconditional stability with respect to the mesh size and that it benefits from \textit{a posteriori} error estimates guaranteeing the quality of the solution.

        To assess the capabilities of the QTT-HL solver, we have conducted a benchmark against other QTT-based methods from the literature, namely the QTT Finite Difference scheme (QTT-FD) of Oseledets and Khoromskij and the Bramble–Pasciak–Xu preconditioned QTT-BPX solver of Bachmayr and Kazeev.
        In particular, the QTT-HL solver proves to be more stable with respect to the mesh size than the QTT-FD solver and compares favorably with the QTT-BPX solver in terms of controlled TT-ranks of the involved operators.
        Also, our solver capabilities extend the well-optimized QTT-based ADI solver of Rakhuba.

        Finally, we illustrate the capabilities of the QTT-HL solver on a challenging multiscale problem. In dimension $d=3$, the solver is able to perform full-field simulations involving a small scale of $\epsilon=10^{-9}\,\mathrm{m}$ over a cubic domain of characteristic size $1\,\mathrm{m}$ with a mesh size below the atomic scale, corresponding to $10^{37}$ degrees of freedom for classical solvers.
		
        This work opens new perspectives for the simulation of realistic multiscale problems for which classical solvers are limited by the prohibitive number of degrees of freedom required, and for which homogenization procedures are not sufficiently generic.
        In particular, problems where there is no clear scale separation (the worse being fractal-like domains, but large-scale modulated microstructures are already challenging), or situations with scale separation into which the large scale and the small scale interact non-trivially, such as corners and interfaces in heterogeneous materials.
        To achieve such a long-term goal, many challenges have to be overcome; to name a few: generalization of the solver to a general domain --not only a cubic one-- with at least the classical boundary conditions;
        combination of the QTT methods with parallel solution strategies in order to treat richer microstructures by increasing the TT-ranks.

        Finally, we emphasize that the efficiency of QTT-based strategies strongly depends on the QTT compressibility of the problem coefficients and solution.
        In analogy with the Kolmogorov $n$-width barriers for linear reduced-order models, QTT approximations may also suffer from intrinsic rank-complexity barriers when the solution manifold exhibits strong non-separable interactions.
        
	\section{Acknowledgment}
	
	Marc Josien thanks Xavier Waintal for introducing him to QTTs, and for many inspiring discussions.
	His student Nicolas Jolly did first experiments concerning the TT-rank of microstructured media made of sums of Gaussians, which motivated the example in Section \ref{Sec:homog}.
	
	We also thank Francesca Cuteri, researcher at the CEA Cadarache for discussions and feed-back, and Sioban Nieradzik-Kozic, who was in internship at the CEA Cadarache, for interesting discussions, and for testing our code Sisyphe and contributing to it.
	
	Finally, we acknowledge the support of \textit{La Mission Numérique} at the CEA, who funded the post-doctoral fellowship of Anas El Hachimi.
	
	\section*{CRediT authorship contribution statement}
	\textbf{Marc Josien}: Conceptualization, Formal Analysis, Software, Methodology, Investigation, Writing – original draft, Writing – review \& editing, Supervision.
	\textbf{Anas El Hachimi}: Investigation, Writing – original draft,Writing – review \& editing.
	\textbf{Isabelle Ramière}: Methodology, Writing – original draft,Writing – review \& editing.

	\section*{Use of AI tools}
		A large language model (Mistral by Mistral AI) was used as a coding aid for formatting arrays and figures. AI-based tools were partially used to assist with language editing.
		The authors performed all the mathematical reasoning and scientific interpretations.
		They take full responsibility for the scientific content of this article.
		
	\section*{Declaration of competing interest}
		The authors declare that they have no known competing financial interests or personal relationships that could have appeared to influence the work reported in this paper.
		
	\section*{License}
	
		For the purpose of Open Access, the authors have applied a CC-BY public copyright license to the present document and to all subsequent versions to the Author Accepted Manuscrit arising from this submission.

\bibliographystyle{plain}
\bibliography{0_Bib_unifie_zot}
  
\appendix

\section{Proofs of Propositions in Section \ref{Sec:NumAnalysis}}\label{App:proofs}

\subsection{Proof of Propsition \ref{Prop:Cond}}

	Both $ \hat{a} *$ and $\mu \hat{p}$ are symmetric and positive, so that their sum is symmetric positive definite.
	Next, we observe that
	\begin{align*}
		0 \leq  \mu \langl\hat\phi, \hat{p} \hat\phi \rangl \leq \mu \|\hat\phi\|_{2}^2,
	\end{align*}
	for $\phi \in \LLb^2(\Q_{1, \per})$.
	Then, by the isometric properties of the Fourier transform, we have
	\begin{align*}
		\langl \hat\phi, \hat{a} * \hat\phi\rangl
		= \int \phi\cdot a \phi.
	\end{align*}
	By our Assumption \ref{A-i}, we obtain
	\begin{align*}
		\lambda^{-1} \|\phi\|_2^2 \leq \int \phi\cdot a \phi \leq  \Lambda \|\phi\|_2^2.
	\end{align*}
	As a consequence, we obtain
	\begin{align*}
		\lambda^{-1} \|\hat\phi\|_2^2 \leq \langl \hat\phi, \lt(\hat a *  + \mu \hat P\rt)\hat\phi \rangl \leq (\mu +  \Lambda) \|\hat\phi\|_2^2,
	\end{align*}
	so that the condition number of the operator  $\hat{a} * + \mu \hat{p}$ satisfies \eqref{Cond}.

\subsection{Proof of Lemma \ref{Prop:Mu}}

	We express the Euler-Lagrange equation associated to \eqref{Jmu_phys} as
	\begin{align}\label{num:001-1}
		\int \phi \cdot a \psi + \phi \cdot g
		+ \mu \phi \cdot P \psi=0,
	\end{align}	
	for any $\phi \in \LLb^2(\Q_{1,\per})$.
	Then, we decompose $\psi$ as $\psi = \psi_{\rm pot} + \psi_{\rm sol}$, where $\psi_{\rm pot} \in \LLpot^2(\Q_{1,\per})$ and $\psi_{\rm sol} \in \LLsol^2(\Q_{1,\per})$.
	Naturally, we may as well decompose $\phi$.
	we obtain from the previous identity \eqref{num:001-1} the following:
	\begin{align}\label{num:0002}
		\int \phi_{\rm pot} \cdot a \psi_{\rm pot} + \phi_{\rm pot}\cdot g =  - \int \phi_{\rm pot} \cdot a \psi_{\rm sol},
	\end{align}
	and
	\begin{align}\label{num:0003}
		\int \phi_{\rm sol} \cdot a \psi_{\rm sol} + \mu \phi_{\rm sol} \cdot \psi_{\rm sol} + \phi_{\rm sol}\cdot g = - \int \phi_{\rm sol} \cdot a \psi_{\rm pot}.
	\end{align}
	Inserting $\phi:=\psi$ in \eqref{num:001-1}, using Assumptions \ref{A-i} and  \ref{A-ii}, we get by the Cauchy-Schwarz inequality
	\begin{align*}
		\lambda^{-1}\int \psi^2 + \mu \int \psi_{\rm sol}^2 \leq \lt(\int g^2\rt)^{\frac{1}{2}} \lt(\int \psi^2\rt)^{\frac{1}{2}}.
	\end{align*}
	Hence, we obtain
	\begin{align}\label{num:0004}
		\|\psi\|_2 \leq \lambda \|g\|_2.
	\end{align}
	Similarly, using $\phi_{\rm sol}:=\psi_{\rm sol}$, in \eqref{num:0003} and using Assumption~\ref{A-i}, we get
	\begin{align*}
		(\mu+\lambda^{-1}) \int |\psi_{\rm sol}|^2 \leq \int |\psi_{\rm sol} \cdot g| + |\psi_{\rm sol} \cdot a \psi_{\rm pot}|,
	\end{align*}
	and hence, by the Cauchy-Schwarz inequality, and using once more Assumptions \ref{A-i} and \ref{A-ii}, we have
	\begin{align*}
		\|\psi_{\rm sol}\|_2
		\leq \frac{1}{\mu+\lambda^{-1}} \lt(\|g\|_2 + \Lambda \|\psi_{\rm pot}\|_2\rt) \overset{\eqref{num:0004}}{\leq}
		\frac{\lambda + \Lambda}{\mu+\lambda^{-1}} \|g\|_2.
	\end{align*}
	Then, comparing $\psi_{\rm pot}$ with the minimizer $\nabla u$ of $\mathcal{J}$, we have from \eqref{num:0002}
	\begin{align*}
		\int \phi_{\rm pot} a (\psi_{\rm pot} - \nabla u) = - \int \phi_{\rm pot} \cdot a \psi_{\rm sol}.
	\end{align*}
	Choosing $\phi_{\rm pot}:= \psi_{\rm pot} - \nabla u$ and appealing to Assumption \ref{A-i} yields
	\begin{align*}
		\|\psi_{\rm pot} - \nabla u\|_2
		\leq \lambda\Lambda \|\psi_{\rm sol}\|_2.
	\end{align*}
	As a consequence, we establish \eqref{num:0005}.

\subsection{Proof of Proposition \ref{Prop:disc}}

	In this proof, the symbol $\lesssim$ means ``$\leq C$'', where $C$ is a constant depending only on $p$, $\lambda$ and $\Lambda$ in \eqref{A-i-eq}, and $C_p$ in \eqref{As-0}.
	
	The proof is done in two steps. First, we assume that $\nu_a=\nu_P=\nu_g=0$,
	and we call $\bhpsi^\circ_{\mu, L}$ the associated minimizer of $\hat{\bm{\mathcal{J}}}_\mu$.
	We prove that:
	\begin{equation}\label{Disc:error}
		\lt\|\psi^\circ_{\mu, L} -\nabla u\rt\|_2 \lesssim
		(2^{L}\epsilon)^{-p} + \mu^{-1},
	\end{equation}
	Second, we compare $\bhpsi^\circ_{\mu, L}$ and $\bhpsi^*_{\mu, L}$ and establish \eqref{Disc:error2}.

	\paragraph{Step 1}
	Let $\bhP^\circ_k = \hat{p}(k)$ and $\bhg^\circ=\hat{g}(k)$, and let $\bhpsi^\circ_{\mu, L}$ be the solution to
	\begin{align}\label{psicirc}
		\bhA \bhpsi^\circ_{\mu, L} + \mu \bhP^\circ \odot \bhpsi^\circ_{\mu, L} = -\bhg^\circ.
	\end{align}
	Notice that, if $\mu=\infty$, $\psi^\circ_{\infty, L} \in \LLpot^2$.
	Hence, applying Céa's Lemma \cite[Lem.\ 2.28]{ern_theory_2004} in the space $\LLpot^2$, we have
	\begin{align*}
		\lt\|\psi^\circ_{\infty, L} -\nabla u\rt\|_2 \lesssim \inf_{\bhpsi \in \LLpot^2}
		\lt\| \sum_{k\in \hat{G}_L} \bhpsi_k \exp(-\ii k \cdot x) - \nabla u \rt\|_2.
	\end{align*}
	Moreover, since $a$ and $g$ are smooth from scale $\epsilon$ downward, by the properties of the Sobolev space $\HH^p$ and by the regularity properties of the elliptic equations, we have
	\begin{align*}
		\inf_{\bhpsi}
		\lt\| \sum_{k\in \hat{G}_L} \bhpsi_k \exp(-\ii k \cdot x) - \nabla u \rt\|_2
		\lesssim 2^{-pL} \|\nabla u\|_{\HH^p(\Q_{1,\per})}
		\lesssim (2^{L}\epsilon)^{-p} \|g\|_{\CC^p(\Q_{1,\per})}
		\lesssim (2^{L}\epsilon)^{-p}.
	\end{align*}
	Moreover, using the isometric properties of the Fourier transform, and adapting Lemma \ref{Prop:Mu} to the discrete case, we obtain
	\begin{equation*}
		2^{dL/2}\lt\|\hat\psi^\circ_{\mu, L} -\hat\psi^\circ_{\infty, L}\rt\|_{2}=
		\lt\|\bhpsi^\circ_{\mu, L} -\bhpsi^\circ_{\infty, L}\rt\|_{\ell^2} \lesssim C \mu^{-1}\|\bhg^\circ\|_{\ell^2} \leq  2^{dL/2}\|g\|_{2}.
	\end{equation*}
	Using  and the triangle inequality and the two previous estimates, we finally get \eqref{Disc:error}.
	
	\paragraph{Step 2}
	We define $\delta\bhpsi = \bhpsi^*_{\mu, L} - \bhpsi^\circ_{\mu, L}$.
	By definition of the minimization process, the Euler-Lagrange equation \eqref{E-lin} is satisfied:
	\begin{align}\label{E-lin*}
		(\bha * + \mu \bhP\odot) \bhpsi^*_{\mu, L} = - \bhg.
	\end{align}
	Comparing this with \eqref{psicirc}, we obtain
	\begin{equation*}
		(\bhA + \mu \bhP\odot) \delta \bhpsi
		= -\delta \bhg
		- (\delta \bhA + \mu \delta \bhP\odot) \bhpsi^\circ_{\mu, L}
		- (\delta \bhA + \mu \delta \bhP\odot) \delta \bhpsi
		.
	\end{equation*}
	Since the spectral norm of the operator $\delta\bhP\odot$ is equal to $\|\bhP\|_{\ell^\infty}$, we obtain from the previous identity and from \eqref{A-i-eq} that
	\begin{align*}
		\|\delta \bhpsi\|_{\ell^2} \leq \lambda \lt( \|\delta \bhg\|_{\ell^2} +
		\|\delta \bhA + \mu \delta \bhP\odot\|_{2, 2} \|\bhpsi^\circ_{\mu, L}\|_{\ell^2} + \|\delta \bhA + \mu \delta \bhP\odot\|_{2, 2} \|\delta \bhpsi\|_{\ell^2}
		\rt),
	\end{align*}
	that is
	\begin{align*}
		\|\delta \bhpsi\|_{\ell^2} \leq \lambda \lt( 2^{dL/2}\nu_g +
		(\nu_a+\mu \nu_P) \|\bhpsi^\circ_{\mu, L}\|_{\ell^2} + (\nu_a+\mu \nu_P)  \|\delta \bhpsi\|_{\ell^2}\rt).
	\end{align*}
	If $\lambda (\nu_a+\mu \nu_P) \leq 1/2$, we may absorb a last right-hand side term into the left-hand side, to the effect of
	\begin{align*}
		\|\delta \bhpsi\|_{\ell^2} \lesssim 2^{dL/2}\nu_g +
		(\nu_a+\mu \nu_P) \|\bhpsi^\circ_{\mu, L}\|_{\ell^2}.
	\end{align*}
	Using the energy estimate on \eqref{psicirc}, we obtain
	\begin{align*}
		\|\delta \bhpsi\|_{\ell^2} \lesssim 2^{dL/2}\nu_g +
		(\nu_a+\mu \nu_P) \|\bhg\|_{\ell^2}.
	\end{align*}
	Using the isometric properties of the Fourier transform and \eqref{As-0}, this turns into
	\begin{align*}
		\|\delta \psi\|_{2} \lesssim \nu_g +
		(\nu_a+\mu \nu_P) \|g\|_{2} 
		\lesssim\nu_g +
		\nu_a+\mu \nu_P.
	\end{align*}
	Combining this with \eqref{Disc:error}, we obtain by the triangle inequality the error \eqref{Disc:error2}.

\subsection{Proof of Proposition \ref{Prop:Aposteriori}}

	The proof is a simple variational argument.
	We expand
	\begin{align}\label{ExpandPsi}
		\phi = \nabla v + \phi_{\rm sol} \qquad \text{with} \quad
		\phi_{\rm sol} = P\phi.
	\end{align}
	Then, there holds
	\begin{align*}
		-\nabla\cdot (a\phi) = \nabla \cdot g + \nabla \cdot \nabla w,
		\qquad \text{with} \quad
		\nabla w = \Gamma (a\phi + g).
	\end{align*}
	As a consequence, we have
	\begin{align}\label{Num:010}
		-\nabla\cdot (a (\nabla v-\nabla u)) = \nabla \cdot (a\phi_{\rm sol}) + \nabla \cdot \nabla w,
	\end{align}
	so that, by the energy estimate and Assumption \ref{A-i}, we obtain
	\begin{align*}
		\|\nabla v-\nabla u \|_2 \leq \lambda \lt(\Lambda \|\phi_{\rm sol}\|_2 + \|\nabla w\|_2\rt).
	\end{align*}
	As a consequence, by the triangle inequality, we finally obtain \eqref{Error-total} in form of
	\begin{align*}
		\|\phi-\nabla u \|_2 \leq \lt((1+\lambda \Lambda) \|\phi_{\rm sol}\|_2 + \lambda\|\nabla w\|_2\rt)
		= \lt((1+\lambda \Lambda) E_P[\phi] + \lambda E_\Gamma[\phi]\rt) \|\nabla u\|_2.
	\end{align*}
	
	To prove the reverse inequality \eqref{Error-total2}, we use the $\LL^2$-orthogonality of \eqref{ExpandPsi} to the effect of
	\begin{align*}
		\|\phi-\nabla u \|_2^2 = \|\nabla v-\nabla u \|_2^2 + \|\phi_{\rm sol}\|_2^2.
	\end{align*}
	Using the energy estimate on $\nabla w$ in \eqref{Num:010}, we get
	\begin{align*}
		\|\nabla w\|_2 \leq \Lambda \lt(\|\phi_{\rm sol}\|_2 + \|\nabla v - \nabla u\|_2 \rt).
	\end{align*}
	Thus, we get \eqref{Error-total2} as a consequence of
	\begin{align*}
		\frac{\|\phi-\nabla u \|_2}{\|\nabla u\|_2} \geq 
		E_P[\phi]  \et 
		\frac{\|\phi-\nabla u \|_2}{\|\nabla u\|_2} \geq \frac{1}{\Lambda \sqrt{2}} E_\Gamma[\phi].
	\end{align*}
  
\end{document}